\documentstyle[twoside,12pt]{article}
\def\kies#1{#1}    

\newif\ifpctex
\pctexfalse
\ifpctex\hoffset=-lin\voffset=-lin\fi
\ifpctex\def\bb{bf}\else\newfont{\bb}{msbm10 scaled\magstep1}\fi

\newcommand{\E}{{\cal E}}

\newcommand{\R}{{\cal R}}

\newcommand{\W}{{\cal W}}

\def\nz{\mbox{\bb N}}
\def\rz{\mbox{\bb R}}

\def\cz{\mbox{\bb C}}

\newcommand{\qed}{\hfill$\Box$\linebreak\medskip\par}

\newcommand{\pr}{{\em Proof:\quad}}

\newcommand{\be}{\begin{enumerate}}
\newcommand{\ee}{\end{enumerate}}
\newcommand{\bdm}{\begin{displaymath}}
\newcommand{\edm}{\end{displaymath}}
\newcommand{\beq}{\begin{equation}}
\newcommand{\eeq}{\end{equation}}
\newcommand{\beqa}{\begin{eqnarray}}
\newcommand{\eeqa}{\end{eqnarray}}
\newcommand{\beqas}{\begin{eqnarray*}}
\newcommand{\eeqas}{\end{eqnarray*}}

\def\smalll{}
\def\text#1{\mbox{#1}}
\def\Box{\sqcup\llap{$\sqcap$}}

\kies{

\catcode`@=11 \catcode`!=11

\expandafter\ifx\csname fiverm\endcsname\relax
  \let\fiverm\fivrm
\fi
  
\let\!latexendpicture=\endpicture 
\let\!latexframe=\frame
\let\!latexlinethickness=\linethickness
\let\!latexmultiput=\multiput
\let\!latexput=\put
 
\def\@picture(#1,#2)(#3,#4){%
  \@picht #2\unitlength
  \setbox\@picbox\hbox to #1\unitlength\bgroup 
  \let\endpicture=\!latexendpicture
  \let\frame=\!latexframe
  \let\linethickness=\!latexlinethickness
  \let\multiput=\!latexmultiput
  \let\put=\!latexput
  \hskip -#3\unitlength \lower #4\unitlength \hbox\bgroup}

\catcode`@=12 \catcode`!=12
 \input pictex.tex 

\catcode`@=11 \catcode`!=11
  
\let\!pictexendpicture=\endpicture 
\let\!pictexframe=\frame
\let\!pictexlinethickness=\linethickness
\let\!pictexmultiput=\multiput
\let\!pictexput=\put

\def\beginpicture{%
  \setbox\!picbox=\hbox\bgroup%
  \let\endpicture=\!pictexendpicture
  \let\frame=\!pictexframe
  \let\linethickness=\!pictexlinethickness
  \let\multiput=\!pictexmultiput
  \let\put=\!pictexput
  \let\input=\@@input   
  \!xleft=\maxdimen  
  \!xright=-\maxdimen
  \!ybot=\maxdimen
  \!ytop=-\maxdimen}

\let\frame=\!latexframe

\let\pictexframe=\!pictexframe

\let\linethickness=\!latexlinethickness
\let\pictexlinethickness=\!pictexlinethickness

\let\\=\@normalcr
\catcode`@=12 \catcode`!=12
}

\newtheorem {lemma}{Lemma}[section]
\newtheorem {theo}{Theorem}[section]
\newtheorem {remark}{Remark}[section]
\newtheorem {example} {Example}[section]
\newtheorem {prop}{Proposition}[section]
\newtheorem {corr}{Corollary}[section]

\newcommand{\sect}[1]{\section{#1}\setcounter{equation}{0}}
\newcommand\Swia{\'Swi\accent'30atek }
\newcommand\st{\,;\,\,}
\newcommand\ov{\overline}
\newcommand\sma{<}

\begin{document}
\title
{Local connectivity of the Julia set of real polynomials}

\author{
Genadi Levin, Hebrew University, Israel
\thanks{e-mail: levin@math.huji.ac.il}\\
Sebastian van Strien, University of Amsterdam, the Netherlands
\thanks{e-mail: strien@fwi.uva.nl.}
}
\date{December 31, 1994 and extended January 27 and April 5, 1995}
\maketitle

\def\IMSmarkvadjust{0 pt}
\def\IMSmarkhadjust{0 pt}
\def\IMSmarkhpadding{0 pt}
\def\IMSpubltext{Published in modified form:}
\def\SBIMSMark#1#2#3{
 \font\SBF=cmss10 at 10 true pt
 \font\SBI=cmssi10 at 10 true pt
 \setbox0=\hbox{\SBF \hbox to \IMSmarkhpadding{\relax}
                Stony Brook IMS Preprint \##1}
 \setbox2=\hbox to \wd0{\hfil \SBI #2}
 \setbox4=\hbox to \wd0{\hfil \SBI #3}
 \setbox6=\hbox to \wd0{\hss
             \vbox{\hsize=\wd0 \parskip=0pt \baselineskip=10 true pt
                   \copy0 \break%
                   \copy2 \break%
                   \copy4 \break}}
 \dimen0=\ht6   \advance\dimen0 by \vsize \advance\dimen0 by 8 true pt
                \advance\dimen0 by -\pagetotal
	        \advance\dimen0 by \IMSmarkvadjust
 \dimen2=\hsize \advance\dimen2 by .25 true in
	        \advance\dimen2 by \IMSmarkhadjust

%
%
  \openin2=publishd.tex
  \ifeof2\setbox0=\hbox to 0pt{}
  \else 
     \setbox0=\hbox to 3.1 true in{
                \vbox to \ht6{\hsize=3 true in \parskip=0pt  \noindent  
                {\SBI \IMSpubltext}\hfil\break
                {\it Annals of Math}~{\bf 147} (1998), 471--541. 
                \vfill}}
  \fi
  \closein2
  \ht0=0pt \dp0=0pt
 \ht6=0pt \dp6=0pt
 \setbox8=\vbox to \dimen0{\vfill \hbox to \dimen2{\copy0 \hss \copy6}}
 \ht8=0pt \dp8=0pt \wd8=0pt
 \copy8
 \message{*** Stony Brook IMS Preprint #1, #2. #3 ***}
}

\SBIMSMark{1995/5}{April 1995}{}


%
\ifx\beginpic\undefined\else \fi

\chardef\oldatcat=\the\catcode`\@
\catcode`\@=11

\newskip\hsssglue \hsssglue=0pt plus 1fill minus 1fill \def\hsss{\hskip\hsssglue}

\newdimen\unitlength \newdimen\linethickness
\newdimen\@picheight \newdimen\@xdim \newdimen\@ydim \newdimen\@len \newdimen\@save
\newcount\@multicount \newcount\@xarg \newcount\@yarg
\newbox\@picbox \newbox\@mpbox

\font\tenln=line10     \font\tenlnw=linew10
\font\tencirc=lcircle10 \font\tencircw=lcirclew10
\font\smallln=linew10 scaled 483 

\def\thinlines{\let\linefont=\tenln \let\circlefont=\tencirc
  \linethickness=\fontdimen8\linefont}
\def\thicklines{\let\linefont=\tenlnw \let\circlefont=\tencircw
  \linethickness=\fontdimen8\linefont}
\thinlines

\def\beginpic(#1,#2)(#3,#4){\@picheight=#2\unitlength
  \setbox\@picbox=\hbox to#1\unitlength\bgroup\let\line=\@line
    \kern-#3\unitlength \lower#4\unitlength\hbox\bgroup\ignorespaces}
\def\endpic{\egroup\hss\egroup
  \ht\@picbox=\@picheight \dp\@picbox=\z@
  \leavevmode\box\@picbox}

\def\put(#1,#2)#3{\raise#2\unitlength\rlap{\kern#1\unitlength #3}\ignorespaces}

\def\multiput(#1,#2)(#3,#4)#5#6{\@multicount=#5
 \@xdim=#1\unitlength \@ydim=#2\unitlength \setbox\@mpbox=\hbox{#6}%
 \loop\ifnum\@multicount>0
   \raise\@ydim\rlap{\kern\@xdim \unhcopy\@mpbox}%
   \advance\@xdim#3\unitlength \advance\@ydim#4\unitlength
   \advance\@multicount\m@ne \repeat\ignorespaces}

\def\makebox(#1,#2)#3{\setbox\@picbox=\hbox to#1\unitlength{\hss#3\hss}%
  \@ydim=\ht\@picbox \advance\@ydim-\dp\@picbox
  \ht\@picbox=#2\unitlength \dp\@picbox=\z@
  \leavevmode\lower.5\@ydim\box\@picbox}

\newif\ifneg
\def\@line(#1,#2)#3{\@xarg=#1 \@yarg=#2 \@len=#3\unitlength \leavevmode
 \ifnum\@xarg<0 \reverseline \else \negfalse \@ydim=\z@\fi
 \ifnum\@xarg=0 \@vline
 \else\ifnum\@yarg=0 \@hline \else\@sline\fi\fi
 \ifneg\kern-\@len\else\@save=\@ydim\fi}
\def\reverseline{\negtrue \kern-\@len \@xarg=-\@xarg
 \@ydim=\@len \multiply\@ydim\@yarg \divide\@ydim\@xarg \@yarg=-\@yarg}

\def\@hline{\vrule height.5\linethickness depth.5\linethickness width\@len}
\def\@vline{\kern-.5\linethickness\vrule width\linethickness
  \ifnum\@yarg<0 height\z@ depth\else depth\z@ height\fi\@len
  \kern-.5\linethickness}

\def\@sline{\setbox\@picbox=\hbox{\linefont \count@=\@xarg \multiply\count@ 8
 \ifnum\@yarg>0 \advance\count@\@yarg \advance\count@-9
 \else \advance\count@-\@yarg \advance\count@ 55 \fi \char\count@}%
 \ifnum\@yarg<0 \@picheight=-\ht\@picbox \advance\@ydim\@picheight
 \else \@picheight=\ht\@picbox \fi
 \@xdim=\wd\@picbox \@save=\@ydim
 \loop\ifdim\@xdim<\@len \raise\@ydim\copy\@picbox
  \advance\@xdim\wd\@picbox \advance\@ydim\@picheight \repeat
 \advance\@xdim-\@len \kern-\@xdim
 \multiply\@xdim\@yarg \divide\@xdim\@xarg \advance\@ydim-\@xdim
 \raise\@ydim\box\@picbox}

\def\vector(#1,#2)#3{\@line(#1,#2){#3}%
 \ifnum\@xarg=0 \@vvector \else\ifnum\@yarg=0 \@hvector \else\@svector\fi\fi}
\def\@hvector{\ifneg\rlap{\linefont\char27}\else
 \smash{\llap{\linefont\char45}}\fi} 
\def\@vvector{\ifnum\@yarg<0 \raise-\@len\rlap{\linefont\char63}%
 \else\setbox\@picbox=\rlap{\linefont\char54}\advance\@len-\ht\@picbox
 \raise\@len\box\@picbox\fi}

\def\@svector{\setbox\@picbox=\hbox to\z@{\linefont
 \ifnum\@yarg<0 \count@=55 \@yarg=-\@yarg \else\count@=-9 \fi
 \ifneg\multiply\@xarg16 \multiply\@yarg2
 \else\hss 
  \ifnum\@xarg>2 \multiply\@xarg9 \multiply\@yarg2 \advance\count@29
  \else\ifnum\@yarg>2 \multiply\@xarg16 \multiply\@yarg9 \advance\count@-20
   \else\multiply\@xarg24 \multiply\@yarg3 \fi\fi\fi
  \advance\count@\@xarg \advance\count@\@yarg \char\count@
  \ifneg\hss\fi}
 \raise\@save\box\@picbox}

\def\disk#1{\@len=#1\unitlength \count@='160 \@diskcirc}
\def\circle#1{\@len=#1\unitlength \count@='140 \@diskcirc}
\def\@diskcirc{\setbox\@picbox=\hbox{\circlefont\char\count@}\@xdim=\wd\@picbox
 \leavevmode \ifdim\@len>15.499\@xdim \@bigdc \else \@smalldc\fi}
\def\@bigdc{\ifnum\count@<'160 \@bigcirc
 \else \@len=15\@xdim \@diskcirc\fi}
\def\@smalldc{{\advance\@len-.5\@xdim
 \loop\ifdim\@xdim<\@len \advance\count@\@ne \advance\@xdim\wd\@picbox\repeat
 \hbox{\circlefont\char\count@}}}
\def\@bigcirc{{\circlefont\count@=15
 \setbox\@picbox=\hbox{\char\count@}\@xdim=\wd\@picbox
 \ifdim\@len>2.5\@xdim \@len=2.5\@xdim\fi
 \advance\@len-.125\wd\@picbox
 \loop\ifdim\@xdim<\@len \advance\count@ 4 \advance\@xdim.25\wd\@picbox\repeat
 \@ydim=.5\@xdim \advance\@ydim.5\linethickness
 \setbox\@picbox=\vbox{\hbox{\char\count@\advance\count@-3\char\count@}%
  \nointerlineskip
  \hbox{\advance\count@\m@ne\char\count@\advance\count@\m@ne\char\count@}}%
 \kern-\@ydim\lower\@ydim\box\@picbox}}

\newif\ifovaltl \newif\ifovaltr \newif\ifovalbl \newif\ifovalbr
\ovaltltrue \ovaltrtrue \ovalbltrue \ovalbrtrue
\def\oval(#1,#2){\@xdim=#1\unitlength \@ydim=#2\unitlength
 {\circlefont \setbox\@picbox=\hbox{\char0}
 \ifdim\@xdim<\wd\@picbox \@xdim=\wd\@picbox\fi
 \ifdim\@ydim<\wd\@picbox \@ydim=\wd\@picbox\fi
 \@save=\@xdim \ifdim\@ydim<\@save \@save=\@ydim \fi
 \count@=39
 \loop \setbox\@picbox=\hbox{\char\count@}\ifdim\@save<\wd\@picbox
  \advance\count@-4 \repeat
 \setbox\strutbox=\hbox{\vrule height\ht\@picbox depth\dp\@picbox width\z@
   \kern\wd\@picbox}%
 \@save=.5\wd\@picbox \advance\@save-.5\linethickness
 \setbox0=\hbox to\@xdim{\ifovaltl\char\count@\else\strut\fi
  \kern-\@save\leaders\hrule height\ifovaltl\linethickness\else\z@\fi\hfil
  \leaders\hrule height\ifovaltr\linethickness\else\z@\fi\hfil\kern\@save
  \ifovaltr\advance\count@-3\char\count@\else\strut\fi\kern-\wd\@picbox}%
  \advance\count@\m@ne
 \setbox2=\hbox to\@xdim{\ifovalbl\char\count@\else\strut\fi
  \kern-\@save\leaders\hrule height\ifovalbl\linethickness\else\z@\fi\hfil
  \leaders\hrule height\ifovalbr\linethickness\else\z@\fi\hfil\kern\@save
  \ifovalbr\advance\count@\m@ne\char\count@\else\strut\fi\kern-\wd\@picbox}%
 \@save=\@ydim \advance\@save-\wd\@picbox \divide\@save 2
 \setbox\@picbox=\vbox{\box0\nointerlineskip
  \hbox to\@xdim{\vrule height\@save width\ifovaltl\linethickness\else\z@\fi
    \hfil\ifovaltr\vrule width\linethickness\kern-\linethickness\fi}%
  \nointerlineskip
  \hbox to\@xdim{\vrule height\@save width\ifovalbl\linethickness\else\z@\fi
    \hfil\ifovalbr\vrule width\linethickness\kern-\linethickness\fi}%
  \nointerlineskip\box2}%
  \@save=.5\@ydim \advance\@save.5\linethickness \leavevmode
  \kern-.5\@xdim \kern-.5\linethickness \lower\@save\box\@picbox}}

\def\cpic#1\endcpic{\vcenter{\hbox{\beginpic#1\endpic}}}


\newdimen\@xi \newdimen\@xii \newdimen\@xiii \newdimen\@xiv
\newdimen\@xpt \newdimen\@xoldpt
\newdimen\@yi \newdimen\@yii \newdimen\@yiii \newdimen\@yiv
\newdimen\@ypt \newdimen\@yoldpt
\def\squine(#1,#2,#3,#4,#5,#6){\setbox\@picbox\hbox{\tencirc q}%
 \global\@xoldpt=#1\unitlength \global\@yoldpt=#4\unitlength \kern\@xoldpt
 \@xi=\@xoldpt \@xii=#2\unitlength \@xiii=#3\unitlength
 \@yi=\@yoldpt \@yii=#5\unitlength \@yiii=#6\unitlength
 \squinerec
 \@xpt=#3\unitlength \@ypt=#6\unitlength \@addpoint
 \raise\@ypt\copy\@picbox}
\newif\iffar
\def\squinerec{\farfalse \testnear\@xi\@xiii \testnear\@yi\@yiii
 \iffar \decast \fi}
\def\testnear#1#2{\@save=#1\advance\@save-#2%
 \ifdim\@save<\z@ \@save=-\@save\fi \ifdim\@save>\p@ \fartrue \fi}
\def\decast{\@xpt=\@xi \advance\@xpt\@xii \divide\@xpt2
 \advance\@xii\@xiii \divide\@xii2
 \@xiv=\@xpt \advance\@xiv\@xii \divide\@xiv2
 \@ypt=\@yi \advance\@ypt\@yii \divide\@ypt2
 \advance\@yii\@yiii \divide\@yii2
 \@yiv=\@ypt \advance\@yiv\@yii \divide\@yiv2
 \begingroup\@xii=\@xpt \@xiii=\@xiv
  \@yii=\@ypt \@yiii=\@yiv \squinerec\endgroup
 \@xpt=\@xiv \@ypt=\@yiv \@addpoint
 \@xi=\@xiv \@yi=\@yiv \squinerec}
\def\@addpoint{
 \global\advance\@xoldpt-\@xpt \wd\@picbox=-\@xoldpt
 \raise\@yoldpt\copy\@picbox \global\@xoldpt=\@xpt \global\@yoldpt=\@ypt}

\catcode`\@=\oldatcat


\sect{Introduction and statements of theorems}
\thispagestyle{empty}

One of the main questions in the field of complex dynamics
is the question whether the Mandelbrot set is locally connected,
and related to this, for which maps the Julia set is
locally connected. In this paper we shall prove the following
\bigskip

\noindent
{\bf Main Theorem}
{\em Let $f$ be a polynomial of the form
$f(z)=z^\ell+c_1$ with $\ell$ an even integer
and $c_1$ real. Then the Julia set of $f$ is either
totally disconnected or locally connected.}

\bigskip

\bigskip
In particular, the Julia set of $z^2+c_1$ is locally connected
if $c_1\in [-2,1/4]$ and totally disconnected 
if $c_1\in \rz\setminus [-2,1/4]$ (note that $[-2,1/4]$ is equal to the set
of parameters $c_1\in \rz$ for which the critical point $c=0$
does not escape to infinity).
This answers a question posed by Milnor, see \cite{Mil1}.
We should emphasize that if the $\omega$-limit set $\omega(c)$
of the critical point $c=0$ is not minimal
then it very easy to see that the Julia set is
locally connected, see for example Section 10. 
Yoccoz \cite{Y} already had shown that each quadratic polynomial
which is only finitely often renormalizable
(with non-escaping critical point and no neutral periodic point)
has a locally connected Julia set.
Moreover, Douady and Hubbard \cite{DH1} already
had shown before  that each
polynomial of the form $z\mapsto z^\ell+c_1$ with an attracting or
neutral parabolic cycle has a locally connected Julia set.
As will become clear,
the difficult case is the infinitely renormalizable
case. In fact, using the reduction method
developed in Section 3 of this paper,
it turns out that in the non-renormalizable case the Main Theorem
follows from some results in \cite{L3} and \cite{L5}, see the final section
of this paper.

We should note that there are infinitely renormalizable
non-real quadratic maps with a non-locally connected
Julia set, see \cite{DH} and \cite{Mil}.
Hence, the results above really depend on the use of
real methods. On the other hand,
Petersen has shown that quadratic polynomials with a
Siegel disc such that the eigenvalue at the neutral fixed point
satisfies some Diophantine condition is locally connected,
see \cite{Pe}.

In principle, the methods of Yoccoz completely break down in the
infinitely renormalizable case and in the case
of polynomials with a degenerate critical point.
The purpose of Yoccoz's methods is to
solve the well-known conjecture about the
local connectedness of the Mandelbrot set
and therefore, some version of our ideas might be helpful
in proving this conjecture. For a survey of the results of
Yoccoz, see for example \cite{Mil} and also \cite{L4}.

We should note also that Hu and Jiang, see  \cite{HJ} and
\cite{Ji1} have shown that for
infinitely renormalizable quadratic maps
which are real and of so-called bounded type,
the Julia set is locally connected. Their result
is heavily based on the complex bounds which
Sullivan used in his renormalization results,
see \cite{Sul} and also the last chapter and in particular
Section VI.5 of \cite{MS} (cf. also \cite{Ji2}).

In fact, our methods enable us to extend Sullivan's result to
the class of all infinitely renormalizable unimodal polynomials
{\em independently of the combinatorial type}!
We should emphasize that these complex bounds form
the most essential ingredient for the renormalization results
of Sullivan \cite{Sul}; in fact in McMullen's approach to
renormalization, see \cite{McM}, these complex bounds play an even
more central role. In the previous proofs
of the complex bounds see \cite{Sul}, and also Section VI.5
of \cite{MS}, it is crucial that the renormalization is of bounded
type and, moreover, the proof is quite intricate.
Therefore we are very happy that our methods
give a fairly easy way to get complex bounds
independently of the combinatorial type of the map
(i.e., only dependent of the degree of the map):
 
\bigskip

\noindent
{\bf Theorem A}\,\,
{\em Let $f$ be a real unimodal polynomial infinitely renormalizable map.
Let $f^{s(n)}\colon V_n\to V_n$ be a renormalization
of this map. Then there exists a polynomial-like extension of this
map $f^{s(n)}\colon \Omega_n'\to \Omega_n$
such that the modulus of $\Omega_n\setminus \Omega_n'$
is bounded from below by a constant which only depends on
$\ell$ and such that the diameter of $\Omega_n$ is at most
a universally bounded constant times the diameter of $V_n$.}

\bigskip

The way we prove that such sets $\Omega_n$ exist is through
cross-ratio estimates. In fact, the estimates are
similar to those that were made previously in \cite{SN}.
In this way, we are able to get the `complex bounds'
of Theorem A similar to those used
by Sullivan in his renormalization result.
Note that our bounds
are completely independent of the combinatorial type of the map.
We should note that
Theorem A and its proof hold for
any renormalization $f^s$ of (a maybe only finitely renormalizable map)
$f$ provided $f^{2s}$ does not have an attracting or neutral
fixed point.

In the non-renormalizable case we also have complex bounds.
Firstly, for each level for which one has a high return one has
a polynomial-like mapping. 
(Our definition of high case also includes what is sometimes called a
central-high return, see the definition in the next section.) 

\bigskip

\noindent
{\bf Theorem B}\,\,
{\em Let $f(z)=z^\ell + c_1$ with $\ell$
an even integer and $c_1$ real be a non-renormalizable polynomial
so that $\omega(c)$ is minimal.
Assume $W$ is the real trace of a central Yoccoz puzzle piece
and $F\colon \cup V^i\to W$ is the corresponding first return map
(on the real line) and assume that this map
has a high return, i.e., assume that $F(V^0)\ni c$ where $V^0$ is the
central interval. Then there exist
topological discs $\Omega^i$ and $\Omega$ with
$\Omega^i\cap \rz=V^i$ and $\Omega\cap \rz=W$ and
a complex polynomial-like extension $G\colon \cup_i\Omega^i\to \Omega$
of $F$. The diameter of the disc $\Omega$ is comparable
to the size of $W$.}

\bigskip
Moreover, one has the following result which follows
from \cite{L3} and \cite{L5} (as was pointed out to
us in an e-mail by Lyubich). Graczyk and \Swia
informed us that they also have a proof of this Theorem C.

\bigskip

\noindent
{\bf Theorem C} \,\, [Lyubich]
{\em Let $f(z)=z^\ell + c_1$ with $\ell$
an even integer and $c_1$ real be a non-renormalizable polynomial
so that $\omega(c)$ is minimal.
If $W$ is the real trace of a central Yoccoz puzzle piece
and $F\colon \cup V^i\to W$ is the corresponding return map
(on the real line). Then after some `renormalizations'
one can obtain an iterate $\tilde F\colon \cup \tilde V^i\to \tilde W$
of $F$ with $\tilde W\subset W$ such that
there exist
topological discs $\tilde \Omega^i$ and $\tilde
\Omega$ with
$\tilde \Omega^i\cap \rz=\tilde V^i$ and $\tilde \Omega\cap\rz= \tilde W$ and
a complex polynomial-like extension $\tilde G\colon \cup_i\tilde
\Omega^i\to \tilde \Omega$
of $\tilde F$. The diameter of the disc $\tilde \Omega$ is comparable
to the size of $\tilde W$.}

\bigskip

Let us say a few words about our proofs.
The main idea behind our proof of the Main Theorem
is to construct generalized polynomial-like mappings
$F_n\colon \cup_i\Omega^i_n\to \Omega_n$ which
coincide on the real line with the first return maps
to certain Yoccoz puzzle-pieces. To do this we first
obtain real bounds to get Koebe space:
these are based on a sophisticated
version of the `smallest interval' argument. They are
a sharper version of those used before by Blokh, Lyubich, Martens,
de Melo, Sullivan, van Strien, \Swia and others.
Using those real bounds and the use of certain Poincar\'e domains
we construct these polynomial-like mappings
and show that the diameter of these domains
is comparable with that of the interval $\Omega_n\cap \rz$.
Next we compare these polynomial-like maps
with those from the Yoccoz puzzle because
the intersection of a Yoccoz puzzle-piece with the Julia
is connected. Next we show that the Julia set of the
polynomial-like mappings of the Yoccoz puzzle coincides
with the Julia set of the polynomial-like mappings $F_n$,
see Section 3.
Since these domains get small, we are able to conclude
local connectivity of the Julia set.

The paper is organized as follows. In Section 2 some background
information is given and in Sections 3 and 4 we give an abstract
description of our method for proving local connectivity
of the Julia set. In Section 5, 6 and 7 we develop real
bounds which will enable to estimate the shape
of the pullbacks of certain discs or other regions.
We should emphasize that the real bounds in these sections
hold for all unimodal maps with negative Schwarzian derivative.
In Sections 8 to 13 we apply these estimates to
several cases. The reader will observe that certain
cases are proved by several methods. For example,
in Section 8 the local connectivity
of the infinitely renormalizable case with $\ell\ge 4$ is
proved, while this case
also follows from the estimates (for a more general
case) in Section 12. However, the domains in Section 8
are discs and those in Section 12 are considerably more complicated.
We believe that for future purposes it might be important to
have good domains, and therefore
even if it was sometimes not necessary for the proofs of our theorems,
we have tried to treat each case in a fairly optimal way.
In the final six pages of this paper --
Section 14 -- we prove Theorem C and complete the
Main Theorem in the non-renormalizable case.

Finally, a short history of this paper since
several others have partial proofs of Theorem A and the Main Theorem
in the quadratic case. Firstly, we were inspired by the
papers of Hu and Jiang, see \cite{HJ} and \cite{Ji1}
where it is shown that infinitely renormalizable maps
of bounded type (where Sullivan's bounds hold)
have a locally connected Julia set. The first widely distributed
version of our paper (dated December 31, 1994)
included the proof of the Main Theorem in the quadratic case,
the infinitely renormalizable case, Theorem A (without doubling)
and also some non-renormalizable cases. 
Subsequently, Theorem B was included
in the version of this paper of January 27, 1995.
Graczyk and \Swia distributed a preprint with a
proof of Theorem A in the quadratic case
on February 3, 1995.
Lyubich and Yampolsky gave an alternative proof
of the Main Theorem and Theorem A in the quadratic case, in a
draft dated February 22, 1995.
The `quadratic' proofs of Graczyk, \'Swi\accent'30atek,
Lyubich and Yampolsky of Theorem A improve our estimates
in certain cases because it sometimes allows one to obtain
annuli with large moduli, but those proofs seem to
heavily rely on the map being quadratic.
(In view of the estimates in \cite{SN}
such large moduli cannot be expected to exist in the higher order
case.) After we told Lyubich about our methods to obtain local
connectivity, he realized the relevance of his
methods, see \cite{L3}, \cite{L5}, for proving local-connectivity
in the non-renormalizable case. In an e-mail dated February 10, 1995,
he told us how to prove Theorem C using these methods, thus
completing the proof of the Main Theorem in the non-renormalizable case.
To make this paper self-contained we added his proof in Section
14 in our paper, in the version of April 5, 1995.

The first author would like to thank the University of Amsterdam
where this work was started.
His research was partially supported by BSF Grant No. 92-00050,
Jerusalem, Israel. We thank Ben Hinkle
for a useful comment and sending us a very detailed list of typos.
We thank Misha Lyubich for telling us about
his results in \cite{L3} and \cite{L5} and
pointing out to us that they imply Theorem C.
We thank Edson Vargas for many discussions and
explanations about the ideas in Section 4 of
\cite{L3}.
Finally, we thank Curt McMullen, Mitsu Shishikura and Greg \Swia
for some very helpful remarks.

\sect{Some notation and some background}

Let $f$ be a real unimodal polynomial. For example,
$f(z)=z^\ell+c_1$ where $\ell$ is even. We find it convenient
to denote the critical point by $c$, i.e., $c=0$. The critical value
is therefore $c_1=f(c)$ and we shall write $c_s=f^s(c)$.
When $w\ne c$ then we shall define $\tau(w)$ to be the point
$\ne w$ so that $f(\tau( w))=f(w)$. For our specific map, we have
$\tau(z)=-z$ but since most results in this paper
do not rely on the specific form of the map $f$ we shall
write $\tau(z)$ rather than $-z$.
If $A,B$ are intervals then we shall write
$[A,B]$ for the smallest interval containing $A$ and $B$.
Furthermore, we shall use the following notation
$$(A,B]=[A,B]\setminus A\,\, ,
\quad [A,B)=[A,B]\setminus B\text{ and }
(A,B)=[A,B]\setminus (A\cup B).$$

As usual, if $J\subset T$ are two intervals and
$L,R$ are the components of $T\setminus J$ 
then we define
$C(T,J)$ to be the cross-ratio of this pair of intervals:
$$C(T,J)=\frac{|J||T|}{|L||R|}.$$
Here $|U|$ stands for the length of an interval $U$.
Cross-ratios play a crucial role in all recent results
in real interval dynamics. Often, it suffices to use
some qualitative estimates based on the
so-called Koebe Principle. In our analysis, we shall need
somewhat sharper estimates, which are based on
direct use of the cross-ratio.
For example, we shall often use the inequality that
$$|L|/|J|\ge C^{-1}(T,J).$$
If $g$ is a map which is monotone on $T$ and $Sg\sma 0$
then 
$$C(gT,gJ)\ge C(T,J)\text{, i.e., }
C^{-1}(T,J)\ge C^{-1}(g(T),g(J)).$$
In our case we shall apply this to maps $g$ of the form
$f^n$. Since $Sf\sma 0$ one has also that $Sf^n\sma 0$ so the
previous inequality applies when we take $g=f^n$ and $f^n|T$ is
monotone. The Koebe Principle states that if $Sg<0$ and
$J\subset T$ are intervals so that $g\colon T\to g(T)$
is a diffeomorphism and so that each component
of $g(T\setminus J)$ has size $\tau|g(J)|$ (i.e.,
$g(T)$ is a {\it $\tau$-scaled neighbourhood} of $g(J)$) then
$|Dg(x)|/|Dg(y)|\le (1+\tau)^2/\tau^2$ for each $x,y\in J$.
The intervals $g(T\setminus J)$ are referred to as `Koebe space'.
We shall also use the following fact:
if $Sf\sma 0$ and if $f^n|T$ is monotone and has a hyperbolic
repelling fixed point, then $f^n(T)\supset T$.

We say that $W$ is a {\it symmetric} interval if it is of the
form $W=[w,\tau(w)]$. The boundary point $w$ is called
{\it nice} if $f^i(w)\notin W$ for all $i>0$. Note that there
are plenty of nice points: each periodic orbit
contains a nice point. Also, if $f$ is not renormalizable,
preimages of the orientation
reversing fixed point of $f$ can
be used to find nice points. This is done in the Yoccoz puzzle,
see also the proof of Theorems B and C. Nice points are also
considered in, for example, the thesis of Martens \cite{Mar},
see also Section V.1 of \cite{MS}.

If $f$ is renormalizable, then we can take for $u_n$
the points which are in the boundary of an interval
$I_n\ni c$ which is mapped into itself in a unimodal way
after $q(n)$ iterates.

If $f$ is not renormalizable then we can
construct a sequence of nice points $u_n$ as follows.
Assume that $f$ has an orientation reversing fixed point
$u_0$. Then we define $u_n$ inductively as follows:
let $k(n)$ be the smallest integer such that
$$(u_{n-1},\tau(u_{n-1}))\bigcap \left(\cup_{i=0}^k f^{-i}(u_0)\right)
\ne \emptyset$$
and let $u_n,\tau(u_n)$ be the points
in this intersection which are nearest to
$c$. If $f$ has no periodic attractor, then $u_n$ is defined
for each $n$. It is easy to see that each $u_n$ is a nice point.

Let us explain why these
nice points play such an important role.
Let $W$ be a symmetric interval with nice boundary points.
Let $$D_W=\{x\st \text{ there exists }k>0\text{ such that }
f^k(x)\in W\}.$$
For $x\in D_W$ let $k(x)$ be the smallest integer
$k>0$ for which $f^k(x)\in W$
and define
$$R_W(x)=f^{k(x)}(x).$$
Let $V$ be the component of $D_W$ which contains $c$ and take $s'\in \nz$
be so that $R_W|V=f^{s'}$. Because $W$ has nice boundary points,
each component -- except the component $V$ -- of the domain of $D_W$
is mapped diffeomorphically by $D_W$ onto $W$.
Clearly, $V$ is symmetric and also has nice
boundary points. Similarly, let $U$ be the components
containing $c$ of the
domain $D_V$ of the first return map $R_V$. Take $s\in \nz$
so that $R_V|U=f^s$. Note that $f(U),\dots,f^s(U)$ are disjoint and
that similarly $f(V),\dots,f^{s'}(V)$ are also disjoint.

We say that $R_V$ has a {\it high return} if $R_V(U)\ni c$.
(We should emphasize that this situation also includes
the so-called central-high return case.)
This implies that $R_V(U)=f^{s}(U)$ contains a component of
$V\setminus \{c\}$ and therefore
$f^{s+i}(U)\supset f^i(U)$ for $i\ge 1$.

It is possible that $U=V=W$ is a periodic interval: in this
case $f$ is renormalizable and $s'=s$ is the period of this interval $V$.
In this case, we certainly can assume that $R_V\colon V\to V$
(which is equal to $f^s$ in this case and consists of one fold)
has a high return: otherwise this return map has a periodic attractor
and therefore we do not have to consider this case.

If $f$ is non-renormalizable and
the critical point of $f$ is recurrent, then
taking $\hat W_{n-1}=[u_{n-1},\tau(u_{n-1})]$
one gets as the domain of $R_{\hat W_{n-1}}$
containing the critical point
the interval $\hat W_n=[u_n,\tau(u_n)]$.
In Theorem B we demand that there are infinitely many
$n$'s for which $R_{\hat W_{n-1}}$ has a {\it high return}.

Finally, as in the complex bounds of Sullivan,
we shall use the Poincar\'e metric
on a slit region in the complex plane.
Given a real interval $T$ we shall write $D_*(T)$ for the
disc which is symmetric with respect to the real line
and which intersects the real line exactly in $T$.
More generally, if $T$ is a bounded real interval and $\alpha\in (0,\pi)$
then $D(T;\alpha)$ will denote the union of two discs
which are symmetric w.r.t. the real axis, intersect the
real line exactly in $T$ and which have an external angle
with the real line of angle $\alpha$.
The reason these sets play an important role, can be
explained as follows. Let $\cz_T=\cz\setminus (\rz\setminus T)$.
The set $\cz_T$ with two infinite slits, carries
a Poincar\'e metric, and with respect to this metric
the set $D(T;\alpha)$ consists of all points whose distance
to $T$ is at most equal to some constant $k(\alpha)$.
From this interpretation and the Schwarz contraction principle,
it follows that if $\phi\colon \cz_T\to \cz_{T'}$ is
a univalent conformal mapping sending $T$ diffeomorphically to $T'$,
then
\beq
\phi(D(T;\alpha))\subset D(T';\alpha).
\eeq
We shall apply this statement, in the following way:

\begin{lemma}
\label{scp}
Let $F\colon \cz\to \cz$ be a real polynomial whose critical
points are on the real line and which maps $T'$ diffeomorphically onto
$T$, then there exists a set $D\subset D(T';\alpha)$
with $D\cap \rz=T'$ which is mapped diffeomorphically onto
$D(T;\alpha)$ by $F$.
\end{lemma}

Often we shall use $\alpha=\pi/2$ and so we define
$$D_*(T)=D(T;\pi/2).$$

\sect{Method showing that the Julia set of two polynomial-like mappings
coincide}

We shall use the fundamental notion of polynomial-like mapping
\cite{DH} or more precisely, we need its
extension due to Lyubich and Milnor from \cite{LM}.
Let $D^0,D^1,\dots,D^i$, and $D$ be topological discs
bounded by piecewise
smooth curves and such that the closures $D^0,\dots,D^i$ are contained
in the interior of $D$, and such that
each the discs $D^0,\dots,D^i$ are pairwise disjoint. Then we call
$$R\colon D^0\bigcup D^1\bigcup \dots \bigcup D^i\to D$$
by $\ell$-{\it polynomial-like} if
$R|_{D^j}$ is a univalent map onto $D$ for each $j=1,\dots,i$
and $R|_{D^0}$ is a $\ell$-fold covering of $D^0$ onto $D$.
If $i=0$ in this definition, we obtain a polynomial-like map
in the original sense of Douady-Hubbard.
 
The {\it filled Julia set} of $R$ is said to be the set
$F_R\subset \cup_{j=0}^i D^j$ of the points $z$ such that $R^k(z)$ is
defined for all $k=i,2,\dots$. The {\it Julia set} $J_R=\partial F_R$.
An equivalent definition of the filled Julia set $F_R$ is:
$$F_R=\bigcap_{k=1}^\infty R^{-k}(D).$$

We shall use an extension
of the Straightening Theorem due to Douady
and Hubbard, \cite{DH}. This extension was also used
in Lemma 7.1 of \cite{LM}, for the case that $i=1$.

\begin{lemma}
\label{gnl1}
Let $R\colon D^0\cup \dots \cup D^i \to D$
be a $\ell$-polynomial-like map. Then $R$ is quasi-conformally
conjugate to a polynomial in neighborhoods of the filled Julia set $F_R$
and filled Julia set of the polynomial.
\end{lemma}
\pr Let us first pick a point $x_0\in D\setminus (D^0\cup \dots
D^i)$ and choose closed simple
curves $\gamma_0,\dots,\gamma_i\colon [0,2\pi]\to \cz$
such that $\gamma_i(0)=\gamma_i(2\pi)=x_0$,
the curves $\gamma_i$ only meet at $x_0$
and $\gamma_i$ surrounds $D^i$. Moreover, we choose the
function $\gamma_i$ to be  smooth and so that
$\frac{d}{dt}\gamma_i(0)$ and $\frac{d}{dt}\gamma_i(2\pi)$
are two vectors based at $x_0$ having an angle $\pi/(i+1)$.

If, for example, $i=1$ then $\gamma_0\cup \gamma_1$
is a figure eight.
Next pick a curve $\gamma$ in $\cz\setminus D$ and a point
$x_1\in \gamma$. Moreover, choose a smooth
function $\phi$ defined on a neighbourhood $N$ of $x_0$ such that
$\phi(x_0)=x_1$ and so that $\phi$ maps
$\gamma_i\cap N$ diffeomorphically to $\gamma\cap \phi(N)$ for each
$i=0,1,\dots,i$. In local coordinates this map
will have an expression of the form $z\mapsto z^{i+1}$
plus higher order terms, i.e., this map $\phi$ will have a critical
point of order $i+1$. Now let $A_j$ be the open annulus between
$\gamma_j$ and $D^j$ and let $A$ be the open annulus between
$\gamma$ and $D$. Moreover, find a smooth map
$\tilde R\colon A_j\to A$ which extends to the closure of these sets
so that it agrees with $R$ on $\partial D^j$
and with $\phi$ on the neighbourhood $N$ of $z_0$.
Choose this extension so that $\phi\colon A_0\to A$
is a $\ell$-covering and
$\phi\colon A_j\to A$ is a diffeomorphism for $j=1,\dots,i$.
This map $\tilde R$ becomes an extension of $R$
if we define it equal to $R$ on $D^0\cup \dots \cup D^i$.
Next choose $r>1$ so that the circle centered at
the origin with radius $r>1$ surrounds
$A_0\cup \dots \cup A_0$.
We can extend $\tilde R$ to a map $\hat R\colon \cz\to \cz$
so that $\hat R(z)=z^{\ell+i}$ for $|z|\ge r$
and so that $\hat R$ coincides with $\tilde R$ on
$A_0\cup \dots \cup A_i$. The map $\tilde R$
on the annulus $\{z\st |z| \sma r\}\setminus (A_0\cup \dots \cup A_i)$
is a $\ell$-covering map to the annulus
$\{z\st |z| \sma (\ell+i)r\}\setminus A$.

Now we use the standard trick from the Straightening Theorem.
Take a standard conformal structure (i.e.,
the Beltrami coefficient $\mu=0$) on the basin
of $\infty$ of $R$ and extend this structure to 
a $L^1$ function $\mu\colon \cz\to \{z\st |z| \sma 1\}$
which is invariant under $\hat R$.
Since $\hat R$ is conformal near infinity and
on $D^0\cup \dots \cup D^0$, there are only a bounded number
of points in each orbit of $\hat R$ where this map is
not conformal. It follows that the supremum of $|\mu(z)|$
is bounded away from one, and by the Measurable Riemann Mapping
Theorem, it follows that there exists a quasiconformal homeomorphism
$h\colon \bar{\cz}\to \bar{\cz}$ with $h(\infty)=\infty$
which has $\mu$ as its Beltrami coefficient.
Since $\mu$ is invariant under $\hat R$, it follows that
$$h\circ \hat R\circ h^{-1}$$ is an holomorphic
$(\ell+i)$-covering. Hence
$\hat R$ is quasiconformally conjugate to
a polynomial map $P$ (of degree $(\ell+i)$).
\qed

A corollary is:

\begin{corr}
\label{gnl2}
The Julia set $J_R$ is the limit set for the preimages
of any point $z\in D$ (except, in the case that $i=0$,
 for the point zero
where zero is the $\ell$-multiple fixed point of $R$).
\end{corr}

We can use all this to show that the Julia set of
two polynomial-like mappings coincide. In the applications
of this we shall later on use for one of these
the polynomial-like mapping of the Yoccoz puzzles.

\begin{prop} (cf. \cite{Ji1}, \cite{McM}.) \label{gnp1}
Let
$$R_1\colon D_1^0\bigcup D_1^1\bigcup \dots\bigcup D_1^i\to D_1,$$
$$R_2\colon D_2^0\bigcup D_2^1\bigcup \dots\bigcup D_2^i\to D_2$$
be two $\ell$-polynomial-like mappings,
such that the critical point $c$ of these maps coincide.
That is, $c\in D_1^0\cap D_2^0$ is the unique
and $\ell$-multiple critical point for both $R_1$ and for $R_2$.
Moreover, assume that the following conditions hold:
\begin{enumerate}
\item
$R_1(z)=R_2(z)$ whenever both sides are defined, so that $R_1$
and $R_2$ are extensions of the same map $R$.
\item
Let $C$ be the component of $D_1\cap D_2$ which contains $R(c)$.
Then also $c\in C$, and there exist precisely
$i$ other points $c^1,\dots,c^i$ so that
$c^j\in D_1^j\cap D_2^j$ and $R(c^j)=R(c)$, and, furthermore,
$c^1,\dots, c^i\in C$.
\end{enumerate}
Under these conditions, the Julia sets of $R_1$ and $R_2$ coincide:
$$J_{R_1}=J_{R_2}.$$
If, additionally, $c\in J_{R_1}$, (and, hence, $c\in J_{R_2}$), then
there exists a component of a preimage $R_2^{-n}(D_2)$, which contains $c$
and is contained in $D_1$.
\end{prop}
\pr For $k=1,2$, let $C^0_k,\dots,C^i_k$ be  the components of
$R^{-1}_k(C)$, such that $c^j\in C^j_k$
when $j\not= 0$, and $c\in C^0_k$.
Firstly, $R_k\colon  C^j_k\to C$ is a
covering, which is just one-to-one if $j\not=
0$, and $R_k\colon  C^0_k\to C$ is a $\ell$-branching covering.
In particular,
boundaries are mapped to boundaries. Since $R_1=R_2$ on the common domain
of definition, we get that, in fact, $C^j_1=C^j_2:=C^j$, for every $j$.
Secondly, because of 2), each component
$C^j$ has a point $c^j$ in common with the component $C$.
Since $C^j$ is connected and is contained in both $D_1$ and $D_2$,
it belongs to a component of $D_1\cap D_2$ containing $c^j$, i.e.,
$C^j\subset C$. Now consider a map
$R\colon  C^0\cup C^1\cup\dots\cup C^i\to C$,
which is one-to-one on every $C^j,j\not= 0$, and $\ell$-to-one on $C^0$.
Take a point $x\in C$. Then $R^{-1}(x)$ is a subset of $C$ and it
consists of $l+i$ points
(counting with multiplicities). That is,
\beq
\text{for any $x\in C$, the sets $R^{-1}_1(x)$ and $R^{-1}_2(x)$
coincide and belong to $C$.}
\label{ster}
\eeq
Starting with $x_0\in C$, we apply the corollary to
Lemma~\ref{gnl2} and (\ref{ster})
to get
$ J_{R_1}=J_{R_2}:=J$. If  $c\in J$, then consider a component $K$
of $J$ containing $c$. Since $K\subset D_1$, there exists a component of a
preimage $R_2^{-n}(D_2)$, which contains $K$ and is contained in $D_1$.
\qed

In the sequel we will use a particular case of Proposition~\ref{gnp1}.
Let us state it separately:

\begin{prop}
\label{gnp2}
Let
$$R_1\colon D_1^0\bigcup D_1^1\bigcup \dots\bigcup D_1^i\to D_1,$$
$$R_2\colon D_2^0\bigcup D_2^1\bigcup \dots\bigcup D_2^i\to D_2$$
be two $\ell$-polynomial-like mappings, such that
the critical points of $R_k$ coincide,
this point $c\in D_1^0\cap D_2^0$ and
is a $\ell$-multiple critical point of both $R_1$ and $R_2$.
Moreover, we assume that the following conditions hold:
\begin{enumerate}
\item
$R_1(z)=R_2(z)$ whenever the both parts are defined, so that $R_1$
and $R_2$ are extensions of a map $R$.
\item
For $k=1,2$, all topological discs $D_k,D_k^0,\dots,D_k^i$ are
symmetric w.r.t. the real line $\rz$ and satisfy
$R_k(\overline z)=\overline {R_k(z)}$.
\item Denoting $I_k=D_k\cap \rz$ and
$I^j_k=D_k^j\cap \rz$, one has
$I_2\subseteq I_1$, $I_2^j\subseteq I_1^j$,
and, for $j=1,\dots,i$, the (real) map $R_k\colon  I_k^j\to I_k$ is one-to-one.
\end{enumerate}
Under these conditions, the Julia sets of $R_1$ and $R_2$ coincide.
If, additionally, $c\in \rz$ lies in the Julia set of $R_1$
(and, hence of $R_2$), then
there exists a component of a preimage $R_2^{-n}(D_2)$, which contains $c$
and is contained in $D_1$.
\end{prop}

\sect{How to construct a polynomial-like mapping?}
Let $f\colon \cz\to \cz$ be a map of the
form $f(z)=z^\ell+c_1$ with $c_1$ real and $\ell$ an even positive integer.
Let $V$ be a (real) symmetric interval with nice boundary points.
Let $U$ be the component of the domain of the
first return map to $V$ containing $c$ and
let $\hat U$ be the component of this map
containing $f(U)\ni c_1$. Take $s$ so that
$R_V|U=f^s$.

\begin{prop}
\label{pqp}
Let $\hat U\supset f(U)$ be the interval
which is mapped diffeomorphically onto $V$
by $f^{s-1}$. Write $v^f=f(v)$, $V=[v,\tau(v)]$
and $\hat U=[\hat u^f,u^f]$. Here $u^f=f(u)$
and $\hat u^f$ is a point which is not the $f$-image of
some real point.
Assume
\beq
|\hat u^f-c_1|
\sma
|v^f-c_1|.
\label{fitin}
\eeq
Moreover, assume that the critical point $c=0$ of
$f$ is recurrent, i.e., all iterates of $c$ under $R_V\colon D_V\to V$
remain in $D_V$ and that $\omega(c)$
is minimal. Then there exists a
$\ell$-polynomial-like mapping
$$R\colon D^0\cup \dots \cup D^i\to D_*(V')$$
such that $c\in J_R$.
Here $V'=V$ if $U\not= V$ (i.e. $f^s\colon U\to V$
is not a renormalization), and
$V'$ is equal to some $\varepsilon$-neighbourhood
of $V$ with $\varepsilon>0$ small enough if
$U=V$ (i.e., when $f^s\colon U\to U$ is a renormalization).
The map $R$ is a real polynomial on each of its components
and $\rz\cap D^i$ are the components of $D_V\cap V$
intersecting points of $\omega(c)$.
\end{prop}
\pr
Since $f\colon \omega(c)\to \omega(c)$ is minimal,
each point $x\in \omega(c)$ is in the domain of
the map $R_V$. By compactness, there exists therefore
a finite covering of $\omega(c)$
of disjoint intervals $I^0,\dots,I^i$ consisting
of components of $R_V$ with $I^0\ni c$.
Let us first consider a component $I^j$ with
$j\ne 0$.
Since then $I^j\notin c$ we get
that $R_V$ maps $I^j$ diffeomorphically onto
$R_V(I^j)=V$ and it follows that there is a region
$D^j$ contained in $D_*(I^j)$ which is mapped
diffeomorphically onto $D_*(V)$ by $R_V$.
So consider $I^0=U=[u,\tau(u)]$. The map $f^{s-1}$
sends $\hat U\supset f(U)=f(I^0)\ni c^1$ diffeomorphically onto
$V$. Again there is a region $D'\ni c_1$ contained in
$D_*(\hat U)$ which is mapped diffeomorphically onto $D_*(V)$
by $f^{s-1}$. Because of (\ref{fitin}),
the $f$-inverse $D^0$ of $D'\subset D_*(\hat U)$
is contained in $D_*(V)$.

In the case of renormalization, we replace $V$ above
by its $\varepsilon$-neighborhood $V'$, with $\varepsilon>0$
so small that (\ref{fitin}) holds for the new points $\hat u^f, v^f$,
and so that the new interval $U$ is strictly inside $V'$
(this is possible since in this case the point $u$ is a repelling
periodic point of $f^s$).
\qed

\begin{remark}
As we will show in Section~\ref{sec6},
one can apply this proposition for any renormalizable unimodal
polynomial $f$ of degree $\ell\ge 4$. In addition, we shall
give a specific bound for the modulus
of the corresponding annuli in Section~\ref{sec6}.
For the degree $\ell=2$ we will need a modification
of the above domains: see Section~\ref{sec7}.
\end{remark}

\sect{Real bounds if $R_V$ has a high return}

As before, let $W$ be a symmetric interval with nice boundary points,
let $R_W$ be the first return map to $W$ and
let $V$ be the domain of $R_W$ containing $c$.
Similarly, let $R_V$ be the first return map to $V$
and $U$ the component of the domain of $R_V$ which contains
$c$. Let $\hat U$ (resp. $\hat V$) be the component of $R_V$
(resp. of $R_W$) containing
the critical value $c_1$. Let $s,s'$ be so that
$R_V|\hat U=f^{s-1}$ and $R_W|\hat V=f^{s'-1}$.

In this section
we will assume that $R_V$ has a high return
and derive some conditions which - when satisfied - will imply
that the component of the $f^{-s}(D_*(V))$
which contains $c$ is contained in $D_*(V)$.

Let $j$ be the component of $\hat U\setminus c_1$ which is outside
$[c_1,c_2]$, i.e., $j=\hat U\setminus f(U)$.
If $R_V$ has a high return, $f^{s-1}|f(U)$ contains $c$
and so we can define
$r$ to be the interval in $f(U)$ which contains $c_1$
and such that $f^{s-1}(r)\ni c$. Furthermore, let
$l$ be the maximal interval having a unique common point with
the boundary point of $\hat U$ outside $[c_1,c_2]$ on which
$f^{s-1}$ is monotone.
We also write,
$$U=[u,\tau(u)],\quad \hat U= [\hat u^f, f(u)]\ni c_1,\quad
t=l\cup j \cup r,\quad V=[v,\tau(v)],$$
$$l'=f^{s-1}(l),\quad j'=f^{s-1}(j), \quad r'=f^{s-1}(r),\quad
t'=l'\cup j' \cup r' $$
and
$$L=f^s(l),\quad J=f^s(j), \quad R=f^s(r)\text{ and }T=L\cup J \cup R.$$
Mark the typographical difference between the degree $\ell$
and the interval $l$.
The situation is drawn below. (The fat lines denote
the part near $c_1$ which is inside the interval
$[c_1,c_2]$; note that the map $f^s|t$ is orientation reversing.)

\hbox to \hsize{\hss\unitlength=1.3mm
\beginpic(70,60)(0,0) \let\ts\textstyle
\put(7,5){\line(1,0){60}}
\put(45,4.7){\line(1,0){22}}
\put(7,25){\line(1,0){60}}
\put(7,45){\line(1,0){60}}
\put(65,44.7){\line(-1,0){58}}
\put(24,-4){{\it The intervals of Lemma~\ref{str1}.}}
\put(10,6){\line(0,-1){2}}
\put(30,6){\line(0,-1){2}}\put(29,1){$\hat u^f$}
\put(45,6){\line(0,-1){2}}\put(44,1){$c_1$}
\put(65,6){\line(0,-1){2}} \put(64,1){$c_{-s+1}$}
\put(15,1){$b$}
\put(20,7){$l$} \put(37,7){$j$} \put(52,7){$r$}
\put(10,26){\line(0,-1){2}}
\put(30,26){\line(0,-1){2}} \put(29,21){$\tau(v)$}
\put(45,26){\line(0,-1){2}} \put(44,21){$c_s$}
\put(65,26){\line(0,-1){2}} \put(64,21){$c$}
\put(20,27){$l'$} \put(37,27){$j'$} \put(52,27){$r'$}
\put(15,10){\vector(0,1){9}}

\put(10,46){\line(0,-1){2}}
\put(15,41){$a$}
\put(30,46){\line(0,-1){2}}
\put(45,46){\line(0,-1){2}}
\put(65,46){\line(0,-1){2}}
\put(65,41){$c_1$}\put(45,41){$c_{s+1}$}\put(30,41){$v^f$}
\put(20,47){$L$} \put(37,47){$J$}\put(52,47){$R$}
\put(15,31){\vector(0,1){9}}

\put(3,15){$f^{s-1}$} \put(3,35){$f$}
\endpic\hss}
\vskip1.2cm

Given $V$ as above, take
$a\in L$ (including possibly $v^f$), choose
$b\in l$ so that $f^s(b)=a$
and define
$$K_\ell(a) =\frac{|b-c_1|}{|a-c_1|}.$$
In the case that $a=v^f$ and $b=\hat u^f$
this becomes
$$K_\ell(v^f) =\frac{|\hat u^f-c_1|}{|v^f-c_1|}=\frac{|j|}{|J\cup R|}.$$
This number is important for our question. Indeed,
if
\beq
K_\ell(v^f)\sma 1
\label{ks1}
\eeq
then, if $f(z)=z^\ell+c_1$ then we get that
$f^{-1}(D_*(\hat u^f,u^f))\subset D_*(v,\tau(v))$.
As in Lemma~\ref{scp} this allows us to get
a polynomial-like extension of
$R_V\colon D_V\to V$.
In this section we shall derive a condition for
(\ref{ks1}).
Define
$$t={|c_1-c_{s+1}|\over |T|}=\frac{|R|}{|T|},$$
$$y={|a-c_1|\over |T|}.$$
This last quantity measures the amount of `extendability'
around $[a,c_1]$. For example,
if $a=v^f$ then $y=\frac{1}{1+|L|/(|J\cup R|)}$
where $|L|/(|J\cup R|)$ is the `space' which exists
around $f(V)=J\cup R$.

\begin{lemma} (See also Proposition 3.2 in \cite{SN}.)
\label{str1}
Assume that $R_V$ has a high return and that
$f^{2s}$ has no neutral or attracting fixed point.
Then
$$K_\ell (a)\le
{t(y^{1/\ell }-t^{1/\ell })\over t^{1/\ell }y(1-y^{1/\ell })}.$$
\end{lemma}
\pr
Denote $\ov J=[a,c_{s+1}]$, $\ov L=T\setminus (R\cup \ov J)$.
Then instead of $j',l'$ we can choose some intervals $\ov j',\ov l'$
so that $f(\ov j')=\ov J$, $f(\ov l')=\ov L$), and
$\ov l,\ov j$ replace the
intervals $l,j$, i.e., $\ov l=f^{-(s-1)}(\ov l')$,
$\ov j=f^{-(s-1)}(\ov j')$).
If we do this, then the intervals $r',r,R$ do not change. Write
$$\alpha=|r'|,\ov \beta=|r'\cup \ov j'|,\gamma=|r'\cup \ov j'\cup
\ov l'|,$$
and
$$Q={|\ov J|\cdot |T|\over |\ov L|\cdot |R|} \,\, / \,\,
{|\ov j'|\cdot |t'|\over |\ov l'|\cdot |r'|}.$$
Then, from the expansion of the cross-ratio's
$${|\ov J|\cdot |T|\over |\ov L|\cdot |R|}
\ge Q\cdot {|\ov j|\cdot |t|\over |\ov l|\cdot |r|}\ge Q.{|\ov
j|\over |r|},$$
and using this inequality we get
$$1/K_l(a)=|\ov J\cup R|/|\ov j|\ge
(|\ov J\cup R|/|r|)\cdot Q\cdot (|\ov L|\cdot |R|)/(|\ov J|\cdot |T|)$$
$$\ge
{{|\ov J|\cdot |T|\over |\ov L|\cdot |R|}
\over {|\ov j'|\cdot |t'|\over |\ov
l'|\cdot |r'|}}
\cdot {|\ov J\cup R|\cdot |\ov L|\over |\ov J|\cdot |T|}
={|\ov J\cup R|\cdot |\ov
l'|\cdot |r'|\over |R|\cdot |\ov j'|\cdot |t'|}$$
$$={\ov \beta^l\cdot (\gamma-\ov \beta)\cdot \alpha\over
\gamma\cdot (\ov \beta-\alpha)\cdot \alpha^l}.$$
Here we have used in the last inequality that
$|r|\le |R|$ which holds
because $f^{2s}|r$ has no periodic attractor.
Now writing
$\ov \beta^l/\gamma^l=y$ and $\alpha^l/\gamma^l=t$
the lemma follows.
\qed

\begin{corr}
\label{cor51}
$$K_\ell(a)  \le K^\ast_\ell (y)={(1-1/\ell )^{\ell -1}\over \ell
\cdot (1-y^{1/\ell })},$$
so that 
$$K^\ast_\ell (y)\to 1/(e\cdot \log(1/y))$$
as $\ell \to \infty$.
\end{corr}

\noindent
\begin{example}.

\noindent
(a) If the extendability space is $0.6$, i.e.,
$y_1=1/(1+0.6)=0.625$ then
$$K^\ast_2(y_1)=1.19371...,\,\, K^\ast_4(y_1)=0.951366...\sma 1.$$

\noindent
(b) If the extendability space is $1/2$, i.e., if $y_2=1/(1+1/2)=2/3$
then 
$$K^\ast_2(y_2)=1.36237...,\,\, K^\ast_4(y_2)=1.0941..,\,\,
K^\ast_6(y_2)=1.02502..,\,\,K^\ast_8(y_2)=0.993...$$

\noindent
(c) If the extendability space is $1/3$, i.e.,
$y_3=1/(1+1/3)=3/4$ then
$$K^\ast_2(y_3)=1.8660...\text{ and }
\lim_{\ell\to \infty} K^\ast_\ell(y_3)\to 1.2788....$$

\noindent
This means that with these
estimates for the extendability space,
we can apply the method suggested by
Proposition~\ref{pqp} respectively for $\ell\ge 4$, $\ell\ge 8$ and not at all
in the last case. (In fact, if we can prove the space
is more than  $e^{1/e}-1=0.44466..$ then we could apply this method
for each $\ell$ sufficiently large.)
In the last section of this paper we shall use a slightly
different method (using different Poincar\'e neighbourhoods)
which also works when the space is equal to $1/3$.
\end{example}

In the next two section we shall derive estimates for
the number $y$ from above.

\sect{Lower bounds for `space' in the renormalizable case}

In this and the next section we shall derive lower bounds
for the number $y$, i.e., find lower bounds for
`space' by looking for a
`smallest' interval among a finite number of intervals.
This idea is used in a large number of
results in one-dimensional dynamics.
In particular we were inspired by
the thesis of Martens \cite{Mar}
or, specifically, by Lemma 1.2 in Section V.1 of \cite{MS}.
In this section we shall obtain quite sharp bounds,
which will enable to deal with all real
infinitely renormalizable maps
$z\mapsto z^\ell+c_1$ of degree $\ell\ge 2$.
Unfortunately, the proof splits in quite a few subcases.
The main result in the section is Lemma~\ref{str5}.
In the next section, we shall obtain weaker bounds which
work in a more general context; these weaker bounds
only apply to the case that $\ell\ge 4$.

Let $\hat V\supset f(V)$ be the interval which is
mapped monotonically onto $W$ by $f^{s'-1}$.

\begin{lemma}
\label{str2}
Let $2\le k\sma s'$ and let
$H_1$ be the maximal interval containing
$f(V)$ such that $f^k|H_1$ is monotone.
Then $f^k(H_1)$ contains $f^k(f(V))$ and
on each side of this interval also an interval
of the form $f^i(V)$ with $i\le k$.
\end{lemma}
\pr
Let $H_{1,-},H_{1,+}$ be the components of $H_1\setminus f(V)$.
From the maximality of $H_1$ it follows that there exists
$i'\sma k$ such that $f^{i'}(H_{1,+})$ contains $c$.
Since $f^{i'+1}(V)$ is outside $W$ it follows that
$f^{i'}(H_{1,+})$ contains one component $W_+$ of $W\setminus \{c\}$.
It follows that $f^k(H_{1,+})$ contains $f^{k-i'}(W)\supset
f^{k-i'}(V)$. Since the same holds for $H_{1,-}$, the lemma follows.
\qed

\hbox to \hsize{\hss\unitlength=1.3mm
\beginpic(70,60)(0,0) \let\ts\textstyle
\put(7,5){\line(1,0){60}}
\put(30,4.7){\line(1,0){37}}
\put(7,25){\line(1,0){60}}
\put(7,45){\line(1,0){60}}
\put(28,-3.5){{\it The proof of Lemma~\ref{str2}.}}
\put(10,6){\line(0,-1){2}}
\put(30,6){\line(0,-1){2}} \put(29,2){$c_1$}
\put(45,6){\line(0,-1){2}}
\put(65,6){\line(0,-1){2}}
\put(20,7){$H_{1,-}$} \put(37,7){$f(V)$} \put(52,7){$H_{1,+}$}
\put(10,26){\line(0,-1){2}} \put(9,22){$c$}
\put(20,26){\line(0,-1){2}}
\put(30,26){\line(0,-1){2}}
\put(45,26){\line(0,-1){2}}
\put(65,26){\line(0,-1){2}}
\put(13,27){$W_+$} \put(32,27){$f^{i'+1}(V)$} \put(48,27){$f^{i'}(H_{1,+})$}
\put(15,10){\vector(0,1){9}}

\put(10,46){\line(0,-1){2}}
\put(20,46){\line(0,-1){2}}
\put(30,46){\line(0,-1){2}}
\put(45,46){\line(0,-1){2}}
\put(55,46){\line(0,-1){2}}
\put(65,46){\line(0,-1){2}}
\put(10,47){$f^{k-i'-1}(W)$} \put(32,47){$f^{k}(f(V))$}
\put(55,47){$f^{k-i"-1}(W)$}
\put(15,31){\vector(0,1){9}}

\put(3,15){$f^{i'}$} \put(3,35){$f^{k-i'-1}$}
\endpic\hss}
\vskip1cm

Let $\hat U\supset f(U)$ be the interval which is mapped
monotonically onto $V$ by $f^{s-1}$.

\begin{lemma}
\label{str3}
Assume that $R_W$ has a high return
and let $l$ be one of the two maximal intervals
outside $\hat U$
for which $f^s|l$ is monotone and which has a unique common point
with $\hat U$. (If we take the interval which is outside $[c_1,c_2]$
then it is equal to the interval $l$ from Lemma~\ref{str1}.)
Then
$L:=f^{s}(l)$ contains an interval
of the form $f^i(V)$, $1\le i \le s'$.
If $U=V$ (so $f$ is renormalizable with period $s$)
and $f$ is not also renormalizable of period $s/2$ then
$L:=f^s(l)$ contains
{\it two distinct} intervals of the form $f^{i+2}(V), f^{s-i}(V)$,
with $1\le i+2, s-i \sma  s'$ and $i+2\ne s-i$.
\end{lemma}
\pr Let $H=l\cup \hat U$.
By maximality of $H$ there is $i$ with $0\sma i\sma s$ such that
$f^i(H)$ contains $c$ in its boundary. Choose $i$ maximal
with this property.
Since $f^i(\hat U)$ is outside
$V$ it follows that $f^i(H)$ contains one component $V_+$ of
$V\setminus \{c\}$.
Hence $f^{i+1}(l)$
contains $f(V)$ and therefore $f^{s-1}(H)$ contains $f^{s-i-1}(V)$
(and also a point in $f^{s-1}(\hat U)=V$).
Since $R_W$ has a high return and $f^{s'-1}(f(V))$ contains
$c$ in its interior,
and since by definition $f^s|l$ is monotone,
it follows that $s-i\sma s'$. Hence $f^{s-1}(H)$ contains one of the intervals
$f(V),\dots,f^{s'-1}(V)$.

Now assume that $U=V$ and take $\tilde H=f^i(H)=[c,f^i(\hat U)]$.
If $f^{s-i-1}(\tilde H)=f^{s-1}(H)$
only contains $f^{s-i-1}(V)$ from the collection
$f(V),\dots,f^{s'}(V)$,  then $f^{s-i-1}(V)=f^{s-i-1}(U)$
is contained in the interval $f^i(\hat U)$
(i.e., $f^{i+1}(U)=f^{s-i-1}(U)$).
Hence $s-i-1=i+1$, i.e., $i+1=s-i-1=s/2$.
It follows that $f^{s-i-1}=f^{s/2}$ maps $[\tilde H,\tau(\tilde H)]$
inside itself. Since, by assumption, this interval
only contains two of the intervals of the
orbit $V,\dots,f^{s-1}(V)$, it follows that
$f$ is also renormalizable of period $s/2$.
Therefore $f^{s-1}(H)$ contains $f^{s-i-1}(V)$
and $f^{i+1}(V)$.
\qed

\hbox to \hsize{\hss\unitlength=1.3mm
\beginpic(70,80)(0,0) \let\ts\textstyle
\put(7,5){\line(1,0){60}}
\put(43,4.7){\line(1,0){24}}
\put(7,25){\line(1,0){60}}
\put(7,45){\line(1,0){60}}
\put(7,65){\line(1,0){60}}
\put(7,64.7){\line(1,0){38}}
\put(28,-3.5){{\it The proof of Lemma~\ref{str3}.}}
\put(15,6){\line(0,-1){2}}
\put(30,6){\line(0,-1){2}}
\put(43,6){\line(0,-1){2}} \put(43,2){$c_1$}
\put(60,6){\line(0,-1){2}}
\put(22,7){$l$}\put(32,7){$\hat U\setminus f(U)$} \put(49,7){$f(U)$}
\put(45,11){$\hat U$}
\put(15,26){\line(0,-1){2}}\put(15,22){$c$}
\put(22,26){\line(0,-1){2}}
\put(30,26){\line(0,-1){2}}
\put(43,26){\line(0,-1){2}}
\put(60,26){\line(0,-1){2}}
\put(18,27){$V_+$}
\put(44,27){$f^{i}(f(U))$}
\put(10,10){\vector(0,1){9}}

\put(15,46){\line(0,-1){2}}
\put(22,46){\line(0,-1){2}}
\put(30,46){\line(0,-1){2}}
\put(45,46){\line(0,-1){2}}\put(45,42){$c$}
\put(60,46){\line(0,-1){2}}
\put(13,47){$f^{s-i-1}(V)$}
\put(40,48){$V=f^{s-1}(\hat U)$}
\put(10,31){\vector(0,1){9}}

\put(15,66){\line(0,-1){2}}
\put(22,66){\line(0,-1){2}}
\put(30,66){\line(0,-1){2}}
\put(45,66){\line(0,-1){2}}
\put(45,62){$c_1$}
\put(17,68){$L=f^s(l)$}
\put(14,60){$f^{s-i}(V)$}
\put(33,68){$f(V)$}
\put(10,51){\vector(0,1){9}}

\put(-5,55){$f$} \put(-5,15){$f^i$} \put(-8,35){$f^{s-i-1}$}
\endpic\hss}
\vskip0.7cm

\begin{lemma}
\label{str4}
Assume that $U=V$ has period $s$ and $f$ is not also renormalizable
of period $s/2$. Consider the disjoint intervals
$f(V),\dots,f^s(V)$ and assume that $f(V)$
and $f^2(V)$ are both smaller than their neighbours.
Then there exists an integer $k\ge 2$ such that $f^k(V)$ is shorter
than its two neighbours from the collection
$f(V),\dots,f^{k-1}(V)$. Take $k$ maximal with respect to this
property.
Let $1\le i_0, i_1\sma k$
be so that $f^{i_0}(V),f^{i_1}(V)$ are the neighbours
of $f^k(V)$ from the collection $f(V),\dots,f^{k-1}(V)$.
Let
$$Q_k=[f^{i_0}(V),f^{i_1}(V)]$$
and define
$H_1\supset f(V)$ to be the maximal interval
on which $f^{k-1}$ is monotone. Then
$H_k=f^{k-1}(H_1)\supset Q_k$. Let
$Z_k\subset H_k$ be the maximal interval
such that each component of $Z_k\setminus Q_k$ contains
at most one interval of the form $f^j(V)$ with
$k<j\le s$. If we define $Z_1\subset H_1$
so that $Z_k=f^{k-1}(Z_1)$ then
$$C^{-1}(Z_1,f(V))\ge 0.6.$$
\end{lemma}
\pr
Such an integer $k$ exists because
otherwise $\{2,\dots,s'\}\ni k\mapsto |f^k(V)|$
would be increasing, contradiction our assumption that
$f^2(V)$ is smaller than its neighbour.
Let $i_0,i_1$ be the intervals as in the statement of the lemma.
By the choice of $k$ these neighbours are longer
than $f^k(V)$.
Let
$$Q_k=[f^{i_0}(V),f^{i_1}(V)]\supset f^k(V).$$
Throughout the remainder of the proof
we shall consider the case that $f^{i_0}(V)$
lies to the left of $f^{i_1}(V)$.
Notice that the fact that $f^{i_0}(V)$ and $f^{i_1}(V)$
are neighbours implies that $Q_k$ only contains intervals
of the form $f^j(V)$ with $k\sma j\le s$. From the
maximality of $k$ this implies that
each such interval $f^j(V)\subset Q_k$ is longer than
the intervals $f^{i_0}(V)$, $f^{i_1}(V)$ and $f^k(V)$.
Lemma~\ref{str2} gives $f^{k-1}(H_1)\supset Q_k$.
Write
$$H_k=f^{k-1}(H_1).$$
Let $Z_1\subset H_1$ be as in the statement of the lemma.
Let $Q_1\supset f(V)$ be the
subset of $H_1$ for which $f^{k-1}(Q_1)=Q_k$
and let
$$Q_{i_j}=f^{i_j-1}(Q_1)\,,\, Z_{i_j}=f^{i_j-1}(Z_1)
\text{ and }
H_{i_j}=f^{i_j-1}(H_1)\text{ for }j=0,1.
$$
Since $f^{i_0}(V)$ and $f^{i_1}(V)$ are longer than
$f^k(V)$ we at least have
$$C^{-1}(Z_1,f(V))\ge
C^{-1}(Q_1,f(V))\ge C^{-1}(Q_k,f^k(V))\ge 1/3.$$
We shall now improve this estimate, by pulling back
the interval $Q=Q_k$ either to $Q_{i_0}$
or to $Q_{i_1}$.  In this way we shall
either find another interval $f^j(V)$ inside the interval
$Q_k$ or find a lower bound
for the space between the intervals $f^j(V)$ in $Q_k$.
For this we shall distinguish between several cases
depending on whether or not $Q_k=H_k$ and depending on
the position of $Q_{i_0}$ and of $Q_{i_1}$ relative to $Q_k$.
Often we shall even show that
$$C^{-1}(Q_k,f^k(V))\ge 0.6.$$
Since $Z_1\supset Q_1$ this suffices:
$$C^{-1}(Z_1,f(V))\ge C^{-1}(Q_1,f(V)).$$

\medskip
\noindent
{\bf Case I.} Assume that $Q_1=H_1$.
By maximality of $H_1$ this implies that there exist $i_0$ and
$i_1$ such that $f^{k-i_0}(H_1)$ and $f^{k-i_1}(H_1)$
contain $c$ in their boundary,
see the figure below.

\hbox to \hsize{\hss\unitlength=1.3mm
\beginpic(70,30)(-15,0) \let\ts\textstyle
\put(3,5){\line(1,0){6}}
\put(13,5){\line(1,0){6}}
\put(23,5){\line(1,0){6}}
\put(13,7){$k-i_1$}\put(25,7){$0$}
\put(40,5){$f^{k-i_1-1}(Q_1)$}
\put(3,15){\line(1,0){6}}
\put(13,15){\line(1,0){6}}
\put(23,15){\line(1,0){6}}
\put(5,17){$0$} \put(13,17){$k-i_0$}
\put(40,15){$f^{k-i_0-1}(Q_1)$}
\put(3,25){\line(1,0){6}}
\put(13,25){\line(1,0){6}}
\put(23,25){\line(1,0){6}}
\put(5,27){$i_0$} \put(15,27){$k$}\put(25,27){$i_1$}
\put(40,25){$f^{k-1}(Q_1)=Q_k$}
\put(8,-2.5){{\it Case I.}}
\endpic\hss}
\vskip0.7cm

If $f^{k-i_0}(V)$ lies closer to $c$ then
$f^{k-i_0+i_1}(V)$ lies between $k$ and $i_1$.
Similarly, if $f^{k-i_1}(V)$ lies closer to $c$ then
$f^{k-i_1+i_0}(V)$ lies between $k$ and $i_0$.
Therefore, since there is no $f^j(V)\subset Q_k$ with $j\sma k$
and $j\ne i_0,i_1$,
this implies that the first possibility occurs if $i_1>i_0$
and the second one if $i_1\sma i_0$.
In order to be definite, we shall assume (in this case) that
$i_0\sma i_1$. This implies that the situation inside $Q_k$
is as drawn below.

\hbox to \hsize{\hss\unitlength=1.3mm
\beginpic(70,10)(10,0) \let\ts\textstyle
\put(23,5){\line(1,0){6}}
\put(33,5){\line(1,0){6}}
\put(43,5){\line(1,0){6}}
\put(58,5){\line(1,0){6}}
\put(25,7){$i_0$} \put(35,7){$k$}\put(40,7){$k+i_1-i_0$}
\put(60,7){$i_1$}
\put(75,5){$Q_k$}
\put(18,-2.5){{\it Case I: The next interval in $Q_k$.}}
\endpic\hss}
\vskip0.7cm

Since each of these intervals $f^{i_0}(V)$,  $f^{i_1}(V)$
and $f^{k+i_1-i_0}(V)$ is at least as long as
$f^k(V)$, it follows that
$$C^{-1}(Q_k,f^k(V))=\frac{\text{left}\cdot \text{right}}
{\text{middle}\cdot \text{total}}\ge
 \frac{2\cdot 1}{1\cdot 4}=\frac{1}{2}.$$
Since
$$C^{-1}(Q_{i_0},f^{i_0}(V))\ge C^{-1}(Q_k,f^k(V))\ge \frac{1}{2}$$
this implies that each of the components of
$Q_{i_0}\setminus f^{i_0}(V)$ has length
$1/2$ times the length of $f^{i_0}(V)$.
Therefore, if $Q_{i_0}$ does not contain any points of
$f^k(V)$,
then one of these components of
$Q_{i_0}\setminus f^{i_0}(V)$
is contained in the gap between $f^{i_0}(V)$
and $f^k(V)$. Hence
$$C^{-1}(Q_k,f^k(V))\ge \frac{2\cdot 3/2}{3/2+1+2}=\frac{6}{9}>0.6\quad .$$
So we are finished in this case.
If there exists an interval $f^j(V)$ between $f^{i_0}(V)$
and $f^k(V)$ then we also are finished, because this interval
then has length $\ge |f^k(V)|$ and therefore
$$C^{-1}(Q_k,f^k(V))\ge \frac{2\cdot 2}{2+1+2}=\frac{4}{5}>0.6\quad .$$
So we shall consider the case that
$Q_{i_0}$ contains some points of $f^k(V)$
and that there exists no interval $f^j(V)$
between $f^{i_0}(V)$ and $f^k(V)$.
Therefore the map $f^{k-i_0}\colon Q_{i_0}\to Q_k$ is orientation preserving.
Indeed, otherwise the interval $[f^{i_0}(V),f^k(V)]$ would be mapped
inside itself by the map $f^{k-i_0}$. Since there is no
interval $f^j(V)$ contained in this interval
$[f^{i_0}(V),f^k(V)]$, this implies that $f$ is also renormalizable
with half the period $s/2$. By assumption this is not the case.
So $f^{k-i_0}\colon Q_{i_0}\to Q_k$ is orientation preserving.
If $Q_{i_0}$ contains (some points of) $f^k(V)$ and no points to the right
of $f^k(V)$ then
the gap between $f^{i_0}(V)$ and $f^k(V)$ is mapped
onto $f^{k+i_1-i_0}(V)$. So if we define
$W_1\subset H_1$ so that $W_k=f^{k-1}(W_1)=[f^{i_0}(V),f^{k+i_1-i_0}(V)]$
then
$$C^{-1}(W_{i_0},f^{i_0}(V))
\ge C^{-1}(W_k,f^k(V))\ge \frac{1}{3}.$$
Hence the component of $W_{i_0}\setminus f^{i_0}(V)$
which is between $f^{i_0}(V)$ and $f^k(V)$
has at least length $1/3$ times the length
of $|f^{i_0}(V)|$. This implies that
$$C^{-1}(Q_k,f^k(V))\ge \frac{4/3\cdot 2}{4/3+1+2}=\frac{8}{13}>0.6.$$
So we finally have to consider the case that
$Q_{i_0}$ strictly contains a neighbourhood of $f^k(V)$. Then
$Q_k$ contains $f^{k+k-i_0}(V)$ and therefore
this entire interval. Since $k+k-i_0>k+i_1-i_0>i_1$,
$Q_k$ contains three intervals of the form $f^j(V)$ to the right of $f^k(V)$
and therefore
$$C^{-1}(Q_k,f^k(V))\ge \frac{1\cdot 3}{1+1+3}=\frac{3}{5}=0.6.$$
This completes case I.

\medskip

\noindent
{\bf Case II and Case III.}
Assume that $Q_1$ is strictly contained in
$H_1$. In order to be specific, let us assume that
$f^{k-1}(H_1)$ contains a neighbourhood of $f^{i_1}(V)$.
This information is useful
since $f^j(V)$ cannot be mapped monotonically
onto $f^{i_i}(V)$ by an iterate of $f$ when $j>i_1$.
In particular, $f^{k-i_0}$ cannot map $f^k(V)$
to $f^{i_1}(V)$.
Therefore there are only two possibilities:
II) $Q_{i_0}$ lies to the left of 
$f^k(V)$ or III) $Q_{i_0}$ contains some points to
the right of $f^k(V)$. 
(Remember that we had assumed that $f^{i_0}(V)$ lies to the left
of $f^{i_1}(V)$.)
\medskip

\noindent
{\bf Case II.}  $Q_{i_0}$ lies to the left of 
$f^k(V)$. 

Now we shall analyze the situation near
$f^{i_1}(V)$. We shall subdivide several cases:

\medskip
\noindent
{\bf Case II.a.}
$Q_{i_1}$ lies to the right of $f^k(V)$.
Let $\alpha \cdot |f^{i_0}(V)|$
be the size of the gap between $f^{i_0}(V)$ and $f^k(V)$.
Similarly, let
$\beta \cdot |f^{i_1}(V)|$
be the size of the gap between $f^{i_1}(V)$ and $f^k(V)$.

\hbox to \hsize{\hss\unitlength=1.3mm
\beginpic(70,20)(10,0) \let\ts\textstyle
\put(23,14){\line(1,0){26}}
\put(13,5){\line(1,0){18}}
\put(38,5){$Q_{i_0}$}
\put(23,15){\line(1,0){6}}
\put(33,15){\line(1,0){6}}
\put(43,15){\line(1,0){6}}
\put(25,17){$i_0$} \put(35,17){$k$} \put(45,17){$i_1$}
\put(30,17){$\alpha$} \put(40,17){$\beta$} 
\put(65,15){$Q_k$}
\put(18,-3.5){{\it Case II: $Q_{i_0}$ lies to the left of $f^k(V)$.}}
\endpic\hss}
\vskip1cm

\noindent
Since $Q_{i_0}$ lies to the left of $f^k(V)$,
$$\alpha\ge
C^{-1}(Q_{i_0},f^{i_0}(V))\ge
C^{-1}(Q_k,f^k(V))\ge \frac{(1+\alpha)(1+\beta)}{3+\alpha+\beta}\ge
1/3.$$
In this case II.a, we also have a similar inequality as above
for $\beta$, i.e., we have
$$\alpha,\beta \ge
\frac{(1+\alpha)(1+\beta)}{3+\alpha+\beta}.$$
Since the right hand side of the above
inequalities is increasing in both $\alpha$ and in $\beta$,
it follows that $\alpha,\beta \ge \kappa$ where
$$\kappa=\frac{(1+\kappa)(1+\kappa)}{3+\kappa+\kappa}\text{, i.e., }
\kappa^2+\kappa= 1.$$
Hence
$$\alpha,\beta\ge \kappa=\sqrt{(5/4)}-1/2 > 0.6\quad .$$
This implies that 
$$C^{-1}(Q_k,f^k(V)) \ge \frac{(1+a)(1+b)}{1(3+a+b)}\ge \kappa \ge 0.6\quad .$$

\medskip
\noindent
{\bf Case II.b.}
$Q_{i_1}$ contains some points of
$f^k(V)$ but no point to the left of
$f^k(V)$. 
As we remarked above, the fact that we
are in Case II or Case III, implies that 
$f^k(V)$ cannot be mapped homeomorphically onto $f^{i_1}(V)$.
It follows that $f^{k-i_1}$ cannot map
$Q_{i_1}$ in an orientation reversing way homeomorphically onto
$Q_k$. Hence $f^{k-i_1}\colon Q_{i_1}\to Q_k$ is orientation preserving
and this map sends $f^k(V)$ to $f^{i_0}(V)$. Writing
$r=k-i_0$, this gives $k+(k-i_1)=i_0+s$, i.e.,
$k-i_1=s-(k-i_0)=s-r$.
So
\beq
k=i_0+r\mbox{ and }i_1=k+(r-s)=i_0+2r-s.
\label{ki}
\eeq
Let $\alpha$ be so that the gap
$(f^{i_0}(V),f^k(V))$ has length $\alpha|f^{i_0}(V)|$
and similarly, let
$\beta$ be so that the gap between $f^k(V)$
and $f^{i_1}(V)$ has size $\beta|f^{i_1}(V)|$.
Let $Z_k^r$ be the right component of $Z_k\setminus f^{i_1}(V)$
and define $\gamma$ so that the size of
$Z_k^r$ is equal to $\gamma|f^{i_1}(V)|$.
Similarly, let $H_k^r$
be the right component of $H_k\setminus f^{i_1}(V)$.

\medskip

\hbox to \hsize{\hss\unitlength=1.3mm
\beginpic(70,20)(0,0) \let\ts\textstyle
\put(24,5){\line(1,0){20}}
\put(60,5){$Q_{i_1}$}
\put(13,14){\line(1,0){26}}
\put(13,15){\line(1,0){6}}
\put(24,15){\line(1,0){5}}
\put(33,15){\line(1,0){6}}
\put(15,17){$i_0$} \put(26,17){$k$}\put(35,17){$i_1$}
\put(60,15){$Q_k$}
\put(8,-7.5){{\it Case II.b.  The map $f^{k-i_1}\colon Q_{i_1}\to Q_k$
is orientation preserving.}}
\endpic\hss}
\vskip1cm

Since we are in Case II, the interval $H_k^r$
contains at least some interval of the form
$f^m(V)$ (with in fact $m\sma k$).
If $Z_k^r\ne H_k^r$ then $Z_k^r$ contains an interval
$f^{j_1}(V)$ with $k \sma j_1 \sma s$. So in any case
$Z_k^r$ contains an interval $f^n(V)$.
If the right component of $Q_{i_1}\setminus f^{i_1}(V)$
is not contained in $Z_k^r$,
then $Q_{i_1}$ contains a neighbourhood of
$f^n(V)$ and therefore
then $Q_k$ contains an  interval
$f^j(V)$ between $f^k(V)$ and $f^{i_1}(V)$.
Since $\alpha\ge 1/3$, this implies that
$C^{-1}(Q_k,f^k(V))\ge (1+1/3)2/(4+1/3)=8/13> 0.6$.
Therefore, we may assume that
the right component of $Q_{i_1}\setminus f^{i_1}(V)$
is contained in $Z_k^r$.
Define $W_1\subset Q_1$ so that
$W_k=f^{k-1}(W_1)=(f^{i_0}(V),f^{i_1}(V)]$.
Since we are in Case II.b, we have that  $f^{k-i_1}$ maps
$f^k(V)$ to $f^{i_0}(V)$. Hence 
one component of $W_{i_1}\setminus f^{i_1}(V)$
is contained in the gap $(f^k(V),f^{i_1}(V))$ corresponding to 
$\beta$ and the other in the interval
$Z_k^r$ corresponding to $\gamma$. Therefore
$$\frac{\beta \gamma}{1+\beta+\gamma}\ge 
C^{-1}(W_{i_1},f^{i_1}(V)).$$
Since
$$C^{-1}(W_{i_1},f^{i_1}(V))\ge C^{-1}(W_k,f^k(V))\ge
\frac{\alpha (1+\beta)}{2+\alpha+\beta}$$
this gives,
$$\frac{\beta \gamma}{1+\beta+\gamma}
\ge
\frac{\alpha (1+\beta)}{2+\alpha+\beta}.$$
Hence
$$\beta(2+\alpha+\beta)
\ge \alpha(1+\beta) + \frac{\alpha(1+\beta)^2}{\gamma},$$
or
$$\beta^2+2\beta \ge \alpha + \frac{\alpha(1+\beta)^2}{\gamma},$$
i.e.,
$$\beta \ge \sqrt{1+\alpha + \frac{\alpha(1+\beta)^2}{\gamma}}-1.$$
Now we study the situation on the other side, around $f^{i_0}(V)$.
Define $\hat Z_k=[f^{i_0}(V),Z_k^r]$
and let $\hat Z_1$ be so that
$f^{k-1}(\hat Z_1)=\hat Z_k$.
Then $\hat Z_{i_0}$ is to the left of $f^k(V)$.
Indeed, $f^r=f^{k-i_0}$ maps $f^k(V)$ to $f^{k+(k-i_0)}(V)=f^{k+r}(V)$.
Because of (\ref{ki}) we have $k+r>s$,
that $f^r|f^k(V)$ is not monotone and $f^{k+r}(V)=f^{i_1}(V)$.  
Since $f^r$ maps $Z_{i_0}$ monotonically to $Z_k$ (which strictly contains
$f^{i_1}(V)$), we finally get that
$Z_{i_0}$ lies to the left of $f^k(V)$.
Hence,
$$\alpha\ge C^{-1}(\hat Z_{i_0},f^{i_0}(V))
\ge C^{-1}(\hat Z_k,f^k(V))\ge
\frac{(1+\alpha)(1+\beta+\gamma)}{3+\alpha+\beta+\gamma}$$
which gives that
$$\alpha^2+2\alpha\ge 1+\beta+\gamma\text{ i.e., }
\alpha \ge \sqrt{2+\beta+\gamma}-1.$$
If $\gamma\ge 0.56$ then
$\alpha\ge \sqrt{2.56}-1=0.6$
and so we have
$$C^{-1}(Z_k,f^k(V))\ge \frac{(1+\alpha)(1+\gamma)}{3+\alpha+\gamma}
\ge
\frac{1.6\cdot 1.56}{4.16}=0.6\quad .$$
If $\gamma\le 0.56$ then
we have that
$$\frac{\alpha}{\gamma}\ge\frac{\sqrt{2+\gamma}-1}{\gamma}
\ge \frac{\sqrt{2.56}-1}{0.56}> 1.$$
Therefore
$$\beta\ge \sqrt{1+\alpha+1}-1.$$
Since we also have $\alpha \ge \sqrt{2+\beta}-1$,
and since the function $(0,\infty)\ni x \mapsto \sqrt{2+x}-1$
has a unique attracting fixed point
$\sqrt{5/4}-1/2>0.61$ it follows that $\alpha,\beta>0.61$.
Again this is sufficient and this completes case II.b.

\medskip
\noindent
{\bf Case II.c.}
$Q_{i_1}$ contains $f^k(V)$ and also some point between $f^{i_0}(V)$
and $f^k(V)$. As before, $f^{k-i_1}\colon Q_{i_1}\to Q_k$
is orientation preserving in this case,
because otherwise $f^{k-i_1}$ maps $[f^k(V),f^{i_1}(V)]$
monotonically into itself and since $f$ has no periodic attractor,
this is impossible. Hence the interval $f^{k+(k-i_1)}(V)$ lies between
$f^{i_0}(V)$ and $f^k(V)$.

\vskip0.3cm

\hbox to \hsize{\hss\unitlength=1.3mm
\beginpic(70,25)(0,0) \let\ts\textstyle
\put(34,5){\line(1,0){25}}
\put(70,5){$Q_{i_1}$}
\put(8,14){\line(1,0){48}}
\put(8,15){\line(1,0){6}}
\put(19,15){\line(1,0){13}}
\put(40,15){\line(1,0){6}}
\put(50,15){\line(1,0){6}}
\put(10,18.5){$i_0$}
\put(19,18.5){$k+(k-i_1)$}
\put(42,18.5){$k$}\put(52,18.5){$i_1$}
\put(59,15.5){$\gamma$}
\put(70,15){$Q_k$}
\put(28,-3.5){{\it Case II.c.}}
\endpic\hss}
\vskip1cm

Note that $Z_k$ contains another interval $f^j(V)$ to the right of $f^{i_1}(V)$
because we have assumed in Cases II and III that $H_k$ contains a neighbourhood
of $f^{i_1}(V)$. Therefore we may assume that $Q_{i_1}$ is contained in
$Z_k$, because otherwise $Q_k$ contains another (i.e., a fifth)
interval $f^{j'}(V)$ inside $Q_k$ and so $C^{-1}(Q_k,f^k(V))$ is
at least $0.6$.
Now let $\gamma$ be so that the length of the component $Z_k\setminus Q_k$
to the right of $f^{i_1}(V)$ is equal to $\gamma|f^{i_1}(V)|$.
Since we have assumed that $Q_{i_1}$ is contained in $Z_k$,
$$\gamma\ge C^{-1}(Q_{i_1},f^{i_1}(V))
\ge C^{-1}(Q_k,f^k(V))\ge {2\cdot 1 \over 4}=0.5\quad .$$
Hence
$$C^{-1}(Z_k,f^k(V))\ge
{2(1+\gamma) \over 2+1+1+\gamma}\ge
{2(1+0.5) \over 4+0.5} = 2/3 > 0.6\quad .$$
This completes the proof of Case II.

\medskip
\noindent
{\bf Case III.}
$Q_{i_0}$ contains a neighbourhood of $f^k(V)$ and, moreover,
$H_k$ contains a neighbourhood of $f^{i_1}(V)$.
The first assumption implies as before that
$f^{k-i_0}\colon Q_{i_0}\to Q_k$ is orientation
preserving.
The last assumption implies that $Q_{i_0}$ cannot
have its right endpoint in some interval
$f^j(V)$ with $j>i_1$ since then $f^{k-i_0}$ cannot
map $f^j(V)$ monotonically onto $f^{i_1}(V)$.
Hence $f^{k+(k-i_0)}(V)$ is contained between $f^k(V)$ and $f^{i_1}(V)$.
If $Q_{i_0}$ contains points from $f^{k+(k-i_0)}(V)$, then
$Q_{i_0}$ contains this interval in its interior and therefore
$Q_k$ contains three intervals $f^j(V)$ to the right of $f^k(V)$
and so we get the required estimate.
Therefore we can (and will) assume in the remainder
of the proof of Case III that
$Q_{i_0}$ is to the left of $f^{k+(k-i_0)}(V)$.

\hbox to \hsize{\hss\unitlength=1.3mm
\beginpic(70,25)(0,0) \let\ts\textstyle
\put(3,5){\line(1,0){6}}
\put(13,5){\line(1,0){6}}
\put(24,5){\line(1,0){5}}
\put(30,5){\line(1,0){1}}
\put(13,14){\line(1,0){48}}
\put(70,15){$Q_k$}\put(70,5){$Q_{i_0}$}
\put(13,15){\line(1,0){6}}
\put(24,15){\line(1,0){5}}
\put(35,15){\line(1,0){12}}
\put(55,15){\line(1,0){6}}
\put(15,19){$i_0$}
\put(20,16){$\alpha$}
\put(25,19){$k$}
\put(31,16){$\beta$}
\put(35,19){$k+(k-i_0)$}\put(57,19){$i_1$}
\put(28,-3.5){{\it Case III.}}
\endpic\hss}
\vskip1cm

Let $\alpha>0$ be so that the gap $(f^{i_0}(V),f^k(V))$
has size $\alpha|f^{i_0}(V)|$.
Define $\beta$ so that 
$(f^k(V),f^{k+(k-i_0)}(V))$ has size $\beta|f^k(V)|$.
The gap corresponding to $\alpha$ is mapped
to the gap corresponding to $\beta$ by $f^{k-i_0}$.
Defining $W_k=[f^{i_0}(V),f^{k+(k-i_0)}(V))$
and $W_1\subset H_1$ so that $f^{k-1}(W_1)=W_k$, we have that
$$\alpha\ge C^{-1}(W_{i_0},f^{i_0}(V))
\ge C^{-1}(W_k,f^k(V))\ge
\frac{(1+\alpha)\beta}{2+\alpha+\beta}.$$
This means that
$$\alpha\ge \sqrt{1+\beta}-1.$$
Now let $\sigma$ be so that
$|f^{k+(k-i_0)}(V)|=\sigma|f^k(V)|$;
one has $\sigma\ge 1$.
Since one component of $Q_{i_0}\setminus f^k(V)$ is contained
in the gap corresponding to $\beta$, and using the definition
of $\sigma$,  we have that
$$\beta\ge C^{-1}(Q_{i_0},f^k(V))
\ge C^{-1}(Q_k,f^{k+(k-i_0)}(V))
\ge
\frac{(2+\alpha+\beta)1}{\sigma(3+\sigma+\alpha+\beta)}.$$
This means that
\beq
\beta^2+\beta(3+\sigma+\alpha-1/\sigma)\ge 2/\sigma+\alpha/\sigma.
\label{in1}
\eeq
Our aim is to prove that
$$C^{-1}(Q_k,f^k(V))
\ge \frac{(1+\alpha)(\beta+\sigma+1)}{3+\alpha+\beta+\sigma}
\ge 0.6\quad .$$
If $\sigma\ge 2$ or if $\alpha\ge 1/3$ then this holds.
So assume that $1\le \sigma\le 2$ and that $\alpha\le 1/3$.
If $\sigma\in [3/2,2]$ then
$3+\sigma+\alpha-1/\sigma\le 5$
and so (\ref{in1}) implies that $\beta^2+5\beta\ge 1$, i.e.,
$\beta\ge 0.19$. Hence $\alpha\ge \sqrt{1.19}-1\ge 0.09$
and hence
$$\frac{(1+\alpha)(\beta+\sigma+1)}{3+\alpha+\beta+\sigma}
\ge \frac{1.09 \times 2.69}{4.78}> 0.6\quad .$$
If $\sigma\in [5/4,3/2]$ then
$3+\sigma+\alpha-1/\sigma\le 4.2$ and therefore
$\beta^2+4.2\beta\ge 4/3$ and so $\beta\ge 0.29$.
Therefore, $\alpha\ge \sqrt{1.29}-1\ge 0.13$
and
$$\frac{(1+\alpha)(\beta+\sigma+1)}{3+\alpha+\beta+\sigma}
\ge \frac{1.13 \times 2.54}{4.67}> 0.6\quad .$$
If $\sigma\in [1,5/4]$ then
$\beta^2+ 3.8 \beta \ge 8/5$ and so $\beta\ge 0.38$.
Therefore, $\alpha\ge \sqrt{1.38}-1\ge 0.17$
and
$$\frac{(1+\alpha)(\beta+\sigma+1)}{3+\alpha+\beta+\sigma}
\ge \frac{1.17 \times 2.38}{4.55}> 0.6\quad .$$
Thus we get the required estimate in each case.
This completes the proof of this lemma.
\qed

The above lemma allows us to show that the interval
$L=f^s(l)$ from Lemma~\ref{str3} is not too short compared to
$f(V)$:

\begin{lemma}
Assume that $R_W$ has a high return.
\label{str5}
As in Lemma~\ref{str3}, let $l$ be one of the two maximal intervals
for which $f^s|l$ is monotone and which has
a unique common point with $\hat U$.
(If we take the interval which is outside $[c_1,c_2]$
then it is equal to the interval $l$ from Lemma~\ref{str1}.)
Write $L=f^s(l)$.
Assume that $U=V$ has period $s$ and assume that
$f$ is not renormalizable of period $s/2$. Then
$$|L|\ge 0.6 \cdot |f(V)|.$$
\end{lemma}
\pr
From the previous lemma,
$L=f^s(l)$ contains at least two intervals of the form
$f^2(V),\dots,f^{s'}(V)$.
We shall consider the disjoint intervals $f(V),\dots,f^{s'}(V)$.
Of course, 
$f(V)$ and $f^2(V)$ have just one neighbour in this collection,
and all other have two.

First consider the case that $f(V)$ is shorter than its neighbour. 
Because $f^s(l)$ contains at least this neighbour,
one gets $|f^s(l)|\ge |f(V)|$ which gives the required estimate.
So we may assume that $f(V)$ is longer than its neighbour.

Similarly, let us consider the
case that $f^2(V)$ is shorter than its nearest neighbour
$f^j(V)$. Then let $G$ be the interval containing $f(V)$
which is mapped diffeomorphically onto $[f^2(V),f^j(V)]$.
Because $z\mapsto |Df(z)|$ is monotone on $G$ and
takes its maximum on $f(V)$ it follows that $G\setminus f(V)$
is longer than $f(V)$. Now since $f^s(l)$ contains
a neighbour, it certainly contains $G\setminus f(V)$.
In particular, as in the previous case we get
$|f^s(l)|\ge |f(V)|$ and the proof is complete
in this situation.

So we may assume that $f(V)$ and $f^2(V)$ are both smaller
than their neighbours. This implies that, as in the previous lemma,
there exists a maximal integer $k\ge 2$ with $k\le s'$
and such that $f^k(V)$ is shorter
than its two neighbours from the collection
$f(V),\dots,f^{k-1}(V)$.
As before,
$Q_k=[f^{i_0}(V),f^{i_1}(V)]\supset f^k(V)$
contains only intervals of the form $f^j(V)$ with $k\sma j\le s$
which are all longer than the
intervals $f^{i_0}(V)$, $f^{i_1}(V)$ and $f^k(V)$.
For simplicity assume again that $f^{i_0}(V)$ lies to the left
of $f^{i_1}(V)$.
Let $H_1\supset f(V)$ be the maximal monotone interval
and, as in the previous lemma,
let $Z_1$ be the maximal interval in $H_1$ such that
$Z_k=f^{k-1}(Z_1)$ contains at most one interval
of the form $f^j(V)$ with $k\sma j \le s$ on each side of $f^k(V)$.
Then
$f^{k-1}(Z_1)\supset Q_k$ by Lemma~\ref{str2}.
Write
$$H_k=f^{k-1}(H_1).$$
Let $Z_{1,-}$ and $Z_{1,+}$
be the components of $Z_1\setminus f(V)$ and
for simplicity take $Z_{1,+}$ be the component which lies on
the same side of $f(V)$ as $L$ (so it lies in $[c_1,c_2]$).  
Label $i_0$ and $i_1$
so that $f^{k-1}(Z_{1,+})$ contains $f^{i_1}(V)$.
If $L\supset Z_{1,+}$ then we have that
$$\frac{|L|}{|f(V)|}
\ge \frac{|Z_{1,+}|}{|f(V)|}
\ge C^{-1}(Z_1,f(V))\ge
C^{-1}(Z_k,f^k(V))\ge 0.6,$$
where in the last inequality we used the previous lemma.

Therefore we may assume that $L$
is a proper subset of $Z_{1,+}\subset H_1$.
Hence $f^{k-1}|L$ is monotone
and therefore, because of Lemma~\ref{str3}, $f^{k-1}(L)$ contains at least
two neighbours $f^j(V)$ with $k\sma j \le s$
of $f^k(V)$ (on the same side as $f^{i_1}(V)$).
These intervals $f^j(V)$ are to the left of $f^{i_1}(V)$
because otherwise $f^{k-1}(L)\subset H_k$ would contain an interval
$f^j(V)$ to the right of $f^{i_1}(V)$ and so
$f^{k-1}(L)\supset Z_{k,+}$. Hence $L\supset Z_{1,+}$,
a contradiction.
If $f^{k-1}(L)$ contains three or more intervals $f^j(V)$
then
$C^{-1}(W_k,f^k(V))\ge 4/5>0.6$ and we are done.
Here $W_1=H_{1,-}\cup f(V) \cup L$ and $W_j=f^{j-1}(W)$.
So we may assume that $f^{k-1}(L)$ contains precisely two neighbours
$f^{j_1}(V)$ and $f^{j_2}(V)$ and assume for simplicity
that $f^{j_1}(V)$ is to the left of $f^{j_2}(V)$.
Of course, this implies that we may also assume that
$L$ contains precisely two intervals of the form $f^j(V)$.
Hence, it suffices to show that
$$C^{-1}(W_k,f^k(V))\ge 0.6.$$
We have that
$f^{i_0}(V)\subset f^{i_0}(H_{1,-})$ lies to the left
of $f^k(V)$.
There are two cases.

\noindent
{\bf Case 1.}
$W_{i_0}$ lies to the left of $f^k(V)$.
In this case choose $\alpha>0$ so that the gap between
$f^{i_0}(V)$ and $f^k(V)$ has size $\alpha|f^{i_0}(V)|$.
Then
$$\alpha \ge C^{-1}(W_{i_0},f^{i_0}(V))
\ge C^{-1}(W_k,f^k(V))
\ge \frac{(1+\alpha)2}{4+\alpha}$$
and this implies that $\alpha\ge 1/2$ and
that
$C^{-1}(W_k,f^k(V))\ge 3/5=0.6$.

\noindent
{\bf Case 2.}
$W_{i_0}$ contains some points of $f^k(V)$.
Then, as before, $f^{k-i_0}\colon W_{i_0}\to W_k$
is order preserving.
If the image of the gap between $f^{i_0}(V)$ and $f^k(V)$ under
$f^{k-i_0}$ contains
one of the intervals $f^j(V)$ in $f^{k-1}(L)\subset W_k$
then define $\hat W_k=[f^{i_0}(V),f^j(V)]$, $\hat W_1\subset H_1$ so that
$f^{k-1}(\hat W_1)=\hat W_k$ and  we get that
$$\alpha \ge C^{-1}(\hat W_{i_0},f^{i_0}(V))
\ge C^{-1}(\hat W_k,f^k(V))
\ge \frac{(1+\alpha)1}{3+\alpha},$$
i.e., $\alpha\ge 1/3$ and $C^{-1}(W_k,f^k(V))>0.6$.
Since, by assumption, $f^{k-1}(L)$ contains no more than two intervals 
$f^j(V)$
and  since $f^{k-i_0-1}(L)$ also contains two interval of the form
$f^m(V)$,  the only remaining possibility is that
$f^{k-i_0}$ maps $f^k(V)$ to $f^{j_1}(V)$
and $f^{j_1}(V)$ to $f^{j_2}(V)$.
(Here we use that $W_{i_0}$ cannot contain
$W_k$
because otherwise $f$ would have a periodic attractor.)
Hence
$$j_1=k+(k-i_0)\text{ and }j_2=k+2(k-i_0).$$
But now we use a more precise statement from Lemma~\ref{str3}:
the intervals $f^j(V)$ in $L$ are of the form
$f^{i+2}(V)$ and $f^{s-i}(V)$. Hence
$$j_1=(i+2)+(k-1)\text{ and }j_2=(s-i)+(k-1)$$
for some $i$.
Writing $r=k-i_0$,
combining all this gives
$$r=j_1-k=i+1\text{ and }2r=j_2-k=s-i-1.$$
Hence $3r=s$
and $L$ contains the intervals $f^{r+1}(V)$ and $f^{2r+1}(V)$.
But then the map $f^{k-1}=f^{r+i_0-1}$ cannot be monotone
on $f^{2r+1}(V)$ (because $3r=s$ and $f^s|f(V)$ is not monotone).
Therefore we get a contradiction with the assumption that
$f^{k-1}|L$ is monotone.
\qed

Finally, we shall also give in this section an estimate for the
case that a map is renormalizable of period $s$ and also of
period $s/2$ (this case was not covered by the previous lemma).

\begin{lemma}
\label{str6}
Assume that $R_W$ has a high return.
As in Lemma~\ref{str3}, let $l$ be one of the two maximal intervals
for which $f^s|l$ is monotone and which has
a unique common point with $\hat U$.
Write $L=f^s(l)$.
Assume that $U=V$ has period $s$ and assume that
$f$ {\bf is} renormalizable of period $s/2$. Then
$$|L|\ge (1/2) \cdot |f(V)|.$$
\end{lemma}
\pr
Let $R$ be the renormalizable interval of period $r:=s/2$
containing both $U$ and $f^{s/2}(U)$.
Let $k=0,1,\dots,r-1$ be so that $|f^k(R)|\le |f^i(R)|$
for each $i=0,1,\dots,r-1$. There are two cases:
$|f^k(U)|\le |f^{k+r}(U)|$ or $|f^{k+r}(U)|\le |f^k(U)|$.
In the former case, let $m=k$ and define $f^{m\pm r}(U)=
f^{m+r}(U)$ and in the latter case we take
$m=k+r$ and define $f^{m\pm r}(U)=f^{m-r}(U)$.
If $m=1,2$ then one has, just like in the
proof of Lemma~\ref{str5}, that
\beq|f^{r+1}(U)|\ge |f(U)|.
\label{str61}
\eeq
(Note that $f^{r+1}(U)$ is the nearest neighbour of $f(U)$ from
the collection of disjoint
intervals $U,f(U),\dots,f^{s-1}(U)$ because by assumption
$f$ is also renormalizable of period $r=s/2$.)
So assume that $m>2$. Then define
$Q_m$ to be the smallest interval containing
$f^{m \pm r}(U)$ on one side of $f^m(U)$ and
containing also the nearest neighbour from the collection
$R,\dots,f^{r-1}(R)$ on the other side of $f^m(U)$.
Let $H_1$ be the maximal interval containing
$\hat U$ so that $f^{m-1}|H_1$ is monotone.
We claim that $f^{m-1}(H_1)\supset f^m(U)$ contains $Q_m$.
Indeed, let $H_1^{1,2}$ be the components of
$H_1\setminus \hat U$. By maximality, there exist
$i_1,i_2<s$ with $i_1\ne i_2$ so that $f^{i_1}(H_1^1)$, $f^{i_2}(H_1^2)$
contains $c$ in its boundary. If $i_j\ne r-1$ then
$f^{i_j}(U)\cap R=\emptyset$ and therefore
$f^{i_j}(H_1^j)$ contains one component of $R\setminus \{c\}$.
Therefore $f^{m-1}(H_1^j)$ contains a neighbour of
$f^m(R)$ from the collection $R,\dots,f^{r-1}(R)$.
If $i_j=r-1$ then $f^{i_j}(H_1^j)$ merely contains one component
of $U\setminus \{c\}$ and then $f^{m-1}(H_1^j)$ contains $f^{m-r}(U)$.
Since either $i_1$ or $i_2$ is different from $r-1$
the claim follows.
By the choice of $m$ we therefore get
$$C^{-1}(Q_m,f^m(U))\ge \frac{1\cdot 2}{1 + 1 + 2}.$$
Hence the interval $Q_1\supset \hat U$ for which
$f^{m-1}(Q_1)=Q_m$ satisfies
$$C^{-1}(Q_1,f(U))\ge 1/2.$$
In particular,
\beq
\mbox{both components of }Q_1\setminus f(U)\mbox{ have length }\ge (1/2)|f(U)|.
\label{str63}
\eeq

From the first part of Lemma~\ref{str3}
it follows that $L=f^s(l)$ contains at least one of neighbours
of $f(U)$ from the collection $U,f(U),\dots,f^{s-1}(U)$.
The nearest neighbour of $f(U)$ is $f^{r+1}(U)$.
So if $m=1,2$ then from (\ref{str61}) it follows that
$|L|\ge |f(U)|$.
If $m>2$ then we have that
either $f^{m-1}(L)$ contains either at least two intervals
from the collection $U,f(U),\dots,f^{s-1}(U)$ or
it contains $f^{r+m}(U)$. Hence from the definition of $Q_m$
and since $f^{m-1}|Q_1$ is monotone
we get that
$L$ contains one component of $Q_1\setminus f(U)$.
In particular, from (\ref{str63}), $|L|\ge (1/2)|f(U)|$.
\qed

\sect{Lower bounds for `space' when $R_V$ has a high return}

Let $\hat V\supset f(V)$ be the interval which is
mapped monotonically onto $W$ by $f^{s'-1}$.
Let $s''$ be the smallest integer such that
$f^{s''-1}(V)\ni c$. In this
section we assume that $R_V$ has a high return.
Hence $f^{s-1}(\hat U)\ni c$ and so $s''\le s$.

\begin{prop} Assume that $R_V$ has a high return.
\label{sp13}
Let $T_0$ be the smallest interval containing
$f(V)$ and another disjoint interval
from the collection $f^2(V),\dots,f^{s''}(V)$.
Write $L_0=T_0\setminus f(V)$.
Then
$$|L_0|\ge \frac{1}{3}|f(V)|.$$
Moreover, if we define
$l$ to be one of the two maximal intervals
outside $\hat U$
for which $f^s|l$ is monotone and which has a unique common point
with $\hat U$, then $f^{s}(l)$ contains $L_0$.
\end{prop}

For the proof of this proposition we need two lemmas.

\begin{lemma} Assume that $R_V$ has a high return.
\label{str2n}
Let $2\le k\sma s''$ and let
$H_1$ be the maximal interval containing
$f(V)$ such that $f^k|H_1$ is monotone.
Then $f^k(H_1)$ contains $f^k(f(V))$ and
on each side of this interval also an interval
of the form $f^i(V)$ with $i\le k$.
\end{lemma}
\pr
Let $H_{1,-},H_{1,+}$ be the components of $H_1\setminus f(V)$.
From the maximality of $H_1$ it follows that there exists
$i'\sma k$ such that $f^{i'}(H_{1,+})$ contains $c$.
Since $i'\sma k \le s''-1 \le s-1$,
the interval $f^{i'+1}(U)$ (which is contained in
$f^{i'+1}(V)\subset f^{i'}(H)$) is outside $V$.
Since $f^{i'+1}(U)$ and $f^{i'+1}(V)$ have $f^{i'+1}(c)$ as one common
endpoint and the other endpoint of $f^{i'+1}(V)$
is certainly outside $V$ (because the endpoints of
$V$ are nice), it follows that $f^{i'+1}(V)$ is outside
$V$. Hence $f^{i'}(H_{1,+})$ contains one component of $V\setminus \{c\}$.
It follows that $f^k(H_{1,+})$ contains $f^{k-i'}(V)$.
\qed

Let $\hat U\supset f(U)$ the interval which is mapped
monotonically onto $V$ by $f^{s-1}$. In the next lemma
we prove the second part of the statement of Proposition~\ref{sp13}.

\begin{lemma}
\label{str3n}
Let $l$ be one of the two maximal intervals
outside $\hat U$
for which $f^s|l$ is monotone and which has a unique common point
with $\hat U$.
Then
$L:=f^{s}(l)$ contains at least one interval
of the form $f^i(V)$, $1\le i \le s''$ which is disjoint from
$f^s(\hat U)=V$.
\end{lemma}
\pr
By maximality of $l$ there is $i$ with $0\sma i\sma s$ such that
$f^i(l)$ contains $c$ in its boundary. Choose $i$ maximal
with this property.
Since $f^i(\hat U)$ is outside
$V$ it follows that $f^i(l)$ contains one component $V_+$ of
$V\setminus \{c\}$.
Hence $f^s(l)$ contains $f^{s-i}(V)$
(and also a point in $f^{s-1}(\hat U)=V$).
Since $f^{s''-1}(f(V))$ contains $c$ in its interior,
and since by definition $f^s|l$ is monotone,
it follows that $s-i< s''$.
Hence $f^{s-1}(H)$ contains one of the intervals
$f(V),\dots,f^{s''-1}(V)$.
\qed

\noindent
{\em Proof of Proposition~\ref{sp13}}
Let $L_0$ be defined as in the
proposition and consider the intervals $f(V),\dots,f^{s''}(V)$.
Since $s''$ might be larger than $s'$, we cannot
be sure that these intervals are disjoint.
Evenso, there exists $k$ such that
$$|f^k(V)|\le |f^i(V)|\mbox{ for each }i=1,2,\dots,s''.$$
If $k=1$ then in particular the length $f(V)$ is at most
equal to the length of its nearest disjoint neighbour
from the collection $f(V),\dots,f^{s''}(V)$.
Since $L_0$ contains an interval $f^i(V)$, $1\le i\le s''$ which is disjoint
from $f(V)$ it follows that $|L_0|\ge |f(V)|$.
If $k=2$ then a similar argument applies: again
the length $f^2(V)$ is at most
equal to the length of its nearest disjoint neighbour
from the collection $f(V),\dots,f^{s''}(V)$.
Since $L_0$ contains an interval $f^i(V)$, $1\le i\le s''$ which is disjoint
from $f(V)$, and since $z\mapsto |Df(z)|$ increases monotonically as
$z$ moves away from $c=0$,  it follows again that $|L_0|\ge |f(V)|$.
So we may assume that $k>2$.
Because of Lemma~\ref{str2n} we can
find an interval
$Z_1$ around $\hat V$ such that $f^{k-1}|Q_1$
is monotone and so that $Z_k=f^{k-1}(Z_1)$ contains on each
side of $f^k(V)$ an interval of the form $f^i(V)$
with $i<k$ (and which is disjoint with $f^k(V)$).
Let $Z_{1,\pm}$ be two components of $Z_1\setminus \hat U$ marked so that
$Z_{1,+}$ intersects $L_0$. If $Z_{1,+}\subset L_0$, then
$${|L_0|\over |f(V)|}\ge {|Z_{1,+}|\over |f(V)|}\ge C^{-1}(Z_1,f(V))\ge
C^{-1}(Z_k,f^k(V))\ge {1\over 3}$$
because $Z_k$ contains on each side of $f^k(V)$
an interval $f^i(V)$ with $i\le s''$ which is at least as long
as $f^k(V)$ and because of the choice of $k$.

It remains to consider the case when $L_0\subset Z_{1,+}$,
i.e., when $f^{k-1}$ is monotone
on $L_0$. As we have seen above,
$L_0$ contains some $f^j(V)$ with $j\le s''$. Hence, the interval
$f^{j+k-1}(V)$ lies in $Z_k$ and $k\le j+k-1\le s''$. By choice of $k$,
the interval $f^{j+k-1}(V)$ is longer than $f^k(V)$.
Hence $|f^{k-1}(L_0)|\ge |f^k(V)|$.
Moreover, again by Lemma~\ref{str2n},
$f^{k-1}(Z_{1,-})$ also contains an interval of the $f^i(V)$ with $i\le s''$
and again by the choice of $k$ this implies that
$|f^{k-1}(Z_{1,-})|\ge |f^k(V)|$.
Let $W_1=Q_{1,-}\cup \hat U\cup L_0$ and $W_k=f^{k-1}(W_1)$.
Then both components of $W_k\setminus f^k(V)$ are at least as long
as $f^k(V)$ and therefore
$${|L_0|\over |f(V)|}\ge C^{-1}(W_1,f(V))\ge
C^{-1}(f^{k-1}(W_1),f^k(V))\ge {1\over 3}.$$
This completes the proof of this proposition.
\qed

\sect{The proof of the Main Theorem in the infinitely renormalizable case
for $\ell>2$}
\label{sec6}
In this section we consider an infinitely renormalizable $f(z)=z^\ell+c_1$
map with $\ell\ge 4$. Such a map has renormalizations of period $q(n)$
where $q(n+1)=q(n)\cdot a(n)$ where $a(n)\ge 2$ is an integer.
If $a(n)=2$ for all $n$ larger than some $n_0$ than
some renormalization has Feigenbaum dynamics, and
then local connectedness immediately follows from
Hu and Jiang's result \cite{HJ}. In fact, we shall
prove the Main Theorem and Theorem A
for this case separately in Section~\ref{secll}
because the bounds obtained in Lemma~\ref{str5}
do not hold at the $n$-th renormalization if $a(n)=2$
(in that case the weaker bounds obtained in Lemma~\ref{str6}
will be used in Section~\ref{secll}.)
So assume in this section that $a(n)>2$ for infinitely many $n$.
Then we can find a sequence of $n$'s
tending to infinity and a sequence of
periodic intervals $U_n=V_n$ of period $s(n)$
such that $f$ is not also renormalizable of period $s(n)/2$.
Because $f$ has no wandering intervals \cite{MS},
it follows that $|U_n|\to 0$.

Let us pick such an $n$ and write $U_n=V_n=[-u_n,u_n]$
and let $s(n)$ be the period (note that for the map we consider
$\tau(z)=-z$). For convenience,
let us for the moment suppress the subscript $n$
and write $s$ for $s(n)$.
Let $\hat U\supset f(U)$ be the interval
from before and consider the diffeomorphism
$f^{s-1}\colon \hat U\to V$.
Let $F$ be the inverse of this map. Since $F$
is a diffeomorphism, it induces a univalent
map
$$F\colon \cz_{V} \to \cz_{\hat U}.$$
Hence $F(D_*(V))\subset D_*(\hat U)$.
Now we use that the space from Lemma~\ref{str1}
is at least $0.6$, because of the estimates of Lemma~\ref{str5}.
(In fact, in Section~\ref{secla}, another proof of the Main
Theorem is given in the infinitely renormalizable
case - this proof is based on the space $1/3$ but using domains
which are not Euclidean discs.)
Hence, see the estimates below Lemma~\ref{str1},
one gets that
$f^{-1}(D_*(\hat U))\subset D_*(V)$, i.e.,
$$f^{-1}\circ F(D_*(V))\subset D_*(V).$$
From this we get that
$$f^{s}\colon D'_*(V)\to D_*(V),$$
where $D'_*(V)=f^{-1}\circ F(D_*(V))$, is a proper degree $\ell$ map.
This is still not a polynomial-like mapping
since $\partial D'_*(V)\cap
\partial D_*(V)=\partial V$, so these regions
intersect in the repelling periodic point
$u$ (with $f^{s}(u)=u$) and its symmetric counterpart $-u$.
Of course, this problem can be easily
amended by adding to $D_*(V)$ some discs
containing $u$ because this point is a repelling.
In this way we will get a polynomial-like mapping
$f^s\colon \Omega_n'\to \Omega_n$.
In fact, we even want
a lower bound for the modulus of the annulus
$\Omega_n\setminus \Omega_n'$.
That such a lower bound exists, is not surprising since
$|Df^s(u)|-1$ is bounded from below,
see Theorem B in Chapter IV of \cite{MS}
or Theorem B in \cite{MMS}.
In the Lemma below we shall give use related estimates
to give lower bounds for the modulus
of the annulus $\Omega_n\setminus \Omega_n'$.

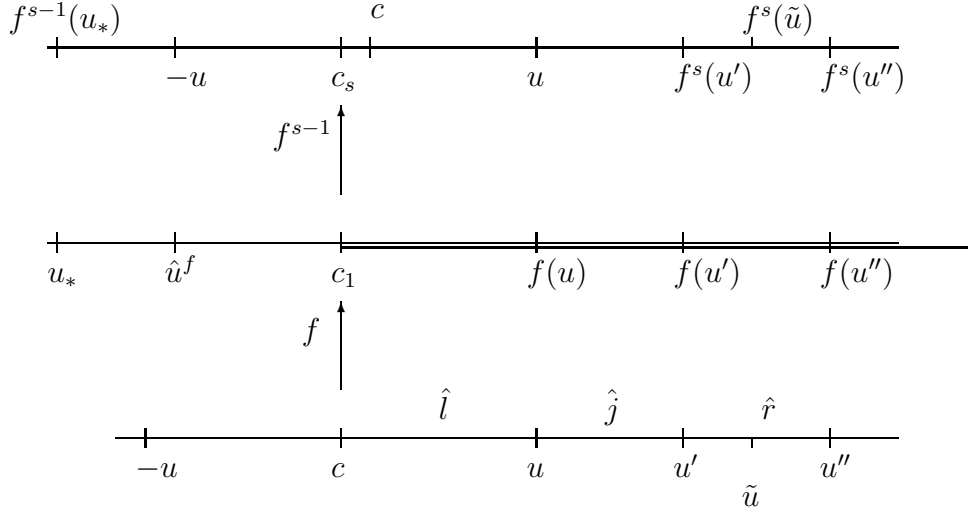
\hbox to \hsize{\hss\unitlength=1.3mm
\beginpic(70,55)(0,0) \let\ts\textstyle
\put(7,5){\line(1,0){80}}
\put(0,25){\line(1,0){87}}
\put(30,24.5){\line(1,0){64}}
\put(0,45){\line(1,0){87}}
%
%
\put(10,6){\line(0,-1){2}}\put(9,1){$-u$}
\put(30,6){\line(0,-1){2}}\put(29,1){$c$}
\put(50,6){\line(0,-1){2}}\put(49,1){$u$}
\put(65,6){\line(0,-1){2}} \put(64,1){$u'$}
\put(72,5){\line(0,-1){1}} \put(71,-2){$\tilde u$}
\put(80,6){\line(0,-1){2}} \put(79,1){$u''$}
\put(40,7){$\hat l$} \put(57,7){$\hat j$} \put(73,7){$\hat r$}
\put(1,26){\line(0,-1){2}}\put(0,21){$u_*$}
\put(13,26){\line(0,-1){2}} \put(12,21){$\hat u^f$}
\put(30,26){\line(0,-1){2}} \put(29,21){$c_1$}
\put(50,26){\line(0,-1){2}} \put(49,21){$f(u)$}
\put(65,26){\line(0,-1){2}} \put(64,21){$f(u')$}
\put(80,26){\line(0,-1){2}} \put(79,21){$f(u'')$}
\put(30,10){\vector(0,1){9}}\put(26,15){$f$}
 \put(1,46){\line(0,-1){2}}\put(-4,47){$f^{s-1}(u_*)$}
\put(13,46){\line(0,-1){2}} \put(12,41){$-u$}
\put(30,46){\line(0,-1){2}} \put(29,41){$c_s$}
\put(33,46){\line(0,-1){2}} \put(33,48){$c$}
\put(50,46){\line(0,-1){2}} \put(49,41){$u$}
\put(65,46){\line(0,-1){2}} \put(64,41){$f^s(u')$}
\put(72,45){\line(0,1){1}} \put(71,47){$f^s(\tilde u)$}
\put(80,46){\line(0,-1){2}} \put(79,41){$f^s(u'')$}
\put(30,30){\vector(0,1){9}}\put(23,35){$f^{s-1}$}
\endpic\hss}
\vskip1cm

\begin{lemma}
\label{lsul}
There are universal constants $C_0,C_1,C_2>0$
with the following property.
Let $T_1\supset f(U)$ be the maximal interval on which $f^{s-1}$ is monotone.
Then,
\beq
\text{ each component of $f^{s-1}(T_1)\setminus V$
has length $\ge \frac{C_0}{\ell}|u-c|$.}
\label{ua}
\eeq
Let $T$ be the component of $f^{-1}(T)\setminus \{c\}$
which contains $u$. Then
there exists $\tilde u\in T$ such that
\beq
\frac{C_1/2}{\ell}|u-c| \le |f^s(\tilde u)-u| \le \frac{C_1}{\ell}|u-c|
\label{ub}
\eeq
such that
\beq
|f^{s}(\tilde u)-c|\ge (1+C_2/\ell^2) \cdot
|\tilde u-c|.
\label{uc}
\eeq
When $\ell\ge 4$, there exists also $u_*\in T$ such that
\beq
f^{s-1}(u_*)=-f^s(\tilde u)\text{ and }
|u_*-c_1| \sma  |f(u)-c_1|.
\label{ud}
\eeq
\end{lemma}
\pr Let $L_1$ be a maximal interval on which $f^s$ is monotone
with a unique common point with $f(U)$.
By Lemma~\ref{str5}, one has
\beq
|f^s(L_1)\setminus f(V)|\ge 0.6 \cdot |f(V)|=0.6\cdot |c_1-f(u)|.
\eeq
Since, $f(z)=(z-c)^\ell+c_1$
(where we usually take $c=0$), the first inequality (\ref{ua})
follows
because there exists $C_0$ such that
\beq
0.6^{1/\ell}\ge 1-\frac{C_0}{\ell}.
\eeq
So we can take $u',u''\in T$ with
\beq
|f^s(u')-u|=\frac{C_1}{2\ell}|c-u|\text{ and }|f^s(u'')-u|=
\frac{C_1}{\ell}|c-u|
\label{sp0}
\eeq
provided $C_1 \sma C_0$.
We shall now show in the remainder of the proof
that there exists $\tilde u\in [u',u'']$ for which
(\ref{uc}) holds. Let us now show that
this would complete the proof, i.e. that
(\ref{ud}) automatically would also hold.
Indeed, because of (\ref{ua}) one can
take $u_*\in T_1$ so that $f^{s-1}(u_*)=-f^s(\tilde u)$.
By Lemma~\ref{str1} there exists a universal constant $K_* \sma 1$
such that provided $\ell\ge 4$,
\beq
|\hat u^f-c_1| \sma  K_*|f(u)-c_1|.
\label{kster}
\eeq
Let us show that, provided we choose $C_1$ sufficiently small,
in (\ref{sp0}), the inequality in (\ref{ud}) holds.
For this define $j'=[u_*,\hat u^f]$ and $t'$ the interval
between $f(u)$ and the endpoint of $T_1$ outside $[c_1,c_2]$.
Let $l',r'$ be the components of $t'\setminus j'$ (with
$l'$ the component outside $[c_1,c_2]$.
One has
\beq
\frac{|c_1-\hat u^f|}{|u_*-\hat u^f|} =
\frac{|l'|}{|j'|}
\ge C^{-1}(t',j')\ge C^{-1}(f^{s-1}(t'),f^{s-1}(j'))
\ge \frac{|f^{s-1}(l')|}{|f^{s-1}(j')|}
\frac{|f^{s-1}(r')|}{|f^{s-1}(t')|}.
\label{uster}
\eeq
The first ratio in the last term will tend to infinity
if we choose $C_1$ small (because $|f^{s-1}(l'\cup j')|$
has size $\ge 0.6^{1/\ell}\cdot |V|$ and
$|f^{s-1}(j')|$ has size $(C_1/\ell)\cdot |V|$).
The second factor in the last term of (\ref{uster}) is of order one.
Hence, (\ref{uster}) implies that $|u_*-\hat u^f|/|c_1-\hat u^f|$
goes to zero provided $C_1$ goes to zero. Combined with
(\ref{kster}) one has that
$|u_*-c_1| \sma |f(u)-c_1|$ for some universal choice of $C_1$.

Thus it remains to prove (\ref{uc}).
Let $\hat l=[c,u]$, $\hat j=[u,u']$, $\hat r=[u',u'']$
and $\hat t=\hat l\cup \hat j \cup \hat r$.
Let us compare the size of $\hat j\cup \hat r$
with that of $\hat l$. Define $\tau=|\hat j\cup \hat r|/|\hat
l|$. If $\tau\leq C_1/(2\cdot \ell)$, then
$|f^s(u'')-c|/|u''-c|=(1+C_1/\ell)/(1+\tau)\geq
(1+C_1/\ell)/(1+C_1/(2\ell))\geq 1+C_2/\ell$,
and (\ref{uc}) follows with $\tilde u=u''$.
So we may assume that
\beq
\tau=\frac{|\hat j\cup \hat r|}{|\hat l|}=\frac{|u-u''|}{|u-c|}\ge C_1/(2\ell).
\label{sp2}
\eeq

Let us show that it suffices to show that
there exists a universal constant
$C_3$ and some $\tilde u\in \hat r=[u',u'']$ for which
\beq
|f^s(u)-f^s(\tilde u)|/|u-\tilde u|\ge (1+C_3/\ell^2).
\label{sp3}
\eeq
Indeed, then
$${|c-f^s(\tilde u)|\over |c-\tilde u|}\geq
{1+\tau (1+C_3/l^2)\over 1+\tau}.$$
Because of (\ref{sp2}), $\tau\ge C_1/(2\ell)$ and therefore
the last expression is bounded from below
by $1+C_7/\ell^2$. Therefore (\ref{ud}) follows and the proof
of the lemma is complete once we have shown that
(\ref{sp3}) holds.

In fact, we may also assume that there exists $C_3'\in (0,1)$ with
\beq
\frac{|f^s(\hat r)|}{|\hat r|}\le
1+C_3'/\ell^2.
\label{sp5}
\eeq
To see this, assume that (\ref{sp5}) fails.
Then, by the Mean Value Theorem,
$$
|Df^s(\xi)|\ge 1+C_3'/\ell^2\text{ for some }
\xi\in \hat r=[u',u''].
$$
Hence we have the following
cross-ratio inequality
$$\frac{(|f^s(u)-f^s(\xi)|/|u-\xi|)^2}{|Df^s(u)||Df^s(\xi)|}
\ge 1$$
(this is just the cross-ratio inequality
$C(f^st,f^sj)\ge C(t,j)$ where we let $t$ shrink
to $j=[u,\xi]$).
Hence, since $u$ is a repelling periodic point,
$$\frac{|f^s(u)-f^s(\xi)|}{|u-\xi|}
\ge \sqrt{|Df^s(u)||Df^s(\xi)|}\ge \sqrt{1\cdot (1+C_3'/\ell^2)}$$
which shows that
(\ref{sp3}) holds with $\tilde u=\xi$
and where we take $C_3=C_3'/4$.
Hence we are also finished if (\ref{sp5}) does not hold.
Of course, for the same reason we may assume
$\frac{|f^s(\hat r)|}{|\hat r|}\le 2$
and with (\ref{sp0}) we get
\beq
\frac{|\hat r|}{|\hat l|}\ge C_4/\ell
\label{spx}
\eeq
where $C_4=C_1/4>0$ is again a universal constant.
Moreover, we have
\beq
\frac{|\hat l|}{|\hat j|}\ge
C^{-1}(\hat t, \hat j)
\ge C^{-1}(f^s(\hat t),f^s(\hat j))
\ge \frac{1\cdot (C_1/(2\ell))}
        {(1+C_1/\ell)(C_1/(2\ell))}
\ge C_4'.
\label{spxa}
\eeq

Now, given intervals $j\subset t$ for which $t\setminus j$
has components $l,r$, define the cross-ratio operator
$$B(t,j)=\frac{|t||j|}{|l\cup j||r\cup j|}.$$
As with the cross-ratio operator $C$ we defined before, one has
$B(ft,fj)\ge B(t,j)$
if $f|t$ is monotone and has negative Schwarzian derivative.
Because of (\ref{spx}) and (\ref{spxa}), and using the expression
$f(z)=(z-\ell)^\ell+c_1$, one can easily check that
$$\frac{B(f(\hat t),f(\hat j))}{B(\hat t,\hat j)}
\ge 1+C_5/\ell^2$$
where $C_5>0$ is a universal constant.
This implies that
\beq
\frac{B(f^s(\hat t),f^s(\hat j))}{B(\hat t,\hat j)}
\ge 1+C_5/\ell^2.
\label{low}
\eeq

Therefore, suppose by contradiction that (\ref{sp3})
is false for $C_3\le \min(C_5/10,1)$.
Since $f^s$ is a repelling fixed point at $u$
and $\hat j$ has $u$ in its boundary,
we have that
$|f^s(\hat j\cup \hat r)|/|\hat j\cup \hat r|\ge 1$.
Because (\ref{sp3}) is false, this implies
\beq
\frac
{|f^s(\hat j)|/|\hat j|}
{|f^s(\hat j\cup \hat r)|/|\hat j\cup \hat r|}
\le 1+C_3/\ell^2.
\label{fac1}
\eeq
Now we also have
\beq
1\le \frac{|f^s(\hat t)|}{|\hat t|}, \frac{|f^s(\hat j\cup \hat l)|}
{|\hat j\cup \hat l|}
\label{sp5a}
\eeq
and
\beq
\frac{|f^s(\hat t)|}{|\hat t|}
\le \max\left(
\frac{|f^s(\hat j\cup \hat l)|}{|\hat j\cup \hat l|}
\quad , \quad
\frac{|f^s(\hat r)|}{|\hat r|}\right).
\label{sp6}
\eeq
Here (\ref{sp5a}) follows from the fact that $f^s|\hat t$
has no attracting fixed point (and $u$ is its fixed point),
and (\ref{sp6})
is a general fact about mean slopes of a monotone map.
Because we assumed that (\ref{sp5}) holds
we get from (\ref{sp5a}) and (\ref{sp6}) that
\beq
\frac
{|f^s(\hat t)|/|\hat t|}
{|f^s(\hat j\cup \hat l)|/|\hat j\cup \hat l|}
\le (1+C_3'/\ell^2).
\label{fac2}
\eeq
Combining (\ref{fac1}) and (\ref{fac2}) and using $C_3'=4C_3$ gives,
$$\frac{B(f^s(\hat t),f^s(\hat j))}{B(\hat t,\hat j)}
\le (1+6C_3/\ell^2).$$
On the other hand, by (\ref{low})
we get that the left hand side of this expression is bounded
from below by $1+C_5/\ell^2$. From this we get
a contradiction since $C_3\le \min(C_5/10,1)$.
\qed

Consider domain
$$\Omega_n=D((-f^s(\tilde u),f^s(\tilde u))=
D((f^{s-1}(u_*),f^s(\tilde u)).
$$
The inverse $F$ of $f^{s-1}$ extends in a univalent way to
$\cz_{[-f^s(\tilde u),f^s(\tilde u)]}$ and therefore
$$F(D(-f^s(\tilde u),f^s(\tilde u)))\subset D_*(u_*,f(\tilde u))$$
where we use that $f^{s-1}(u_*)=-f^s(\tilde u)$.
Because of (\ref{uc}) and (\ref{ud})
we get that
$$f^{-1}(D_*(u_*,f(\tilde u)))\subset
D_*(-f^s(\tilde u),f^s(\tilde u))=\Omega_n$$
and that, moreover,
the difference set is an annulus with
a modulus which is bounded away from zero
by a constant which only depends on $\ell$.
In particular, writing $\Omega'_n=f^{-1}\circ F(\Omega_n)$,
the map $f^{s(n)}\colon \Omega_n'\to \Omega_n$ is a polynomial-like
mapping. Moreover, since $f$ has no wandering intervals
and $s(n)\to \infty$ we have that
$|\Omega_n\cap \rz|\to 0$. By construction the diameter of $\Omega_n$
is at most twice that of $V_n=[-u_n,u_n]$.
Hence, Theorem A in the introduction follows immediately
from all this.

\medskip

\noindent
{\em The conclusion of the proof of the Main Theorem for the
infinitely renormalizable case when $\ell\ge 4$.}
For simplicity, let $T_n$ be the image under $f^{s(n)-1}$
of the maximal interval of monotonicity around $c_1$.
By (\ref{ua}), there exists a constant $C=C(\ell)>0$
such that $C(T_n, \Omega_n\cap \rz)\le C$.
To end the proof, we follow an argument similar to \cite{Ji1}.
Given $n$, consider so-called the maximal renormalization,
$g_n\colon W'_n\to W_n$, where $W_n$ is a slitted complex plane $\cz_{T_n}$
without a fixed neighborhood of infinity, so that
$\Omega_n\subset W_n$ and $f_n=g_n=f^s$ on $\Omega'_n$.
By Proposition~\ref{gnp2}, the Julia set $J_n$ of
$g_n$ is contained in $\Omega_n$. On the other
hand, $J_n$ is equal to the intersection of critical Yoccoz pieces started
from a piece based on the points $u,\hat u$, so there is a piece $P_n$
in $\Omega_n$ containing $c$. Since $P_n\cap J(f)$ is connected,
we have proved
the local connectivity of $J(f)$ at the critical point $c$.

Let $z\in J(f)$ be any
other point. If the forward
orbit of $z$ avoids some neighborhood of $c$, then $f$ is expanding along
this orbit, and the local connectivity at $z$ follows.

So assume that the orbit of $z$
hits any neighborhood $P_n$ of $c$. Then we use the fact that
the renormalizations $f_n$ are `unbranched' (\cite{McM}).
More precisely,
there exists a domain $M_n$, such that
$M_n\subset W_n$, the annulus $A_n=M_n\setminus \overline {P_n}$
has a modulus $\ge m=m(l)>0$, and $A_n$ does not contain any iteration
of $c$. To see this, denote by $t_n$ and $T'_n$, where $t_n\subset
T'_n\subset T_n$, the intervals, such that $g_n(t_n)=\Omega_n\cap
\rz$ and $g_n(T'_n)=T_n$. Since $C(T_n, \Omega_n\cap \rz)\le C$, we have
$C(T'_n, t_n)\le C'$, where $C'>0$ depends only on $\ell$. Furthermore,
the interval $T'_n$ is disjoint from the all iterations of $c$ which are
outside of the renormalization interval $V$ (we use that $f$ is not
$s/2$-renormalizable). On the other hand, by the 
construction, the domain $\Omega'$
is inside of a domain of a definite shape based on the interval $t_n$
(if $\ell\ge 4$, this is just a disc). This implies the existence of $M_n$
as above.
Let some iterate $f^k(z)$, $k=k_n$, of the orbit of $z$ hit $P_n$ the first
time. 
We can pull back the domain $P_n$ along $z,f(z),\dots,f^{k-1}(z)$
by a branch $G_n$ of $f^{-k}$, since $P_n$
contains only iterations of $c$ corresponding to the renormalization.
Indeed, otherwise some $P'=f^{-i}(P_n)$ covers $c$ for the 
first time. Hence, 
$f^i\colon P'\to P_n$ is an iteration of the renormalization $f_n$, i.e.,
$P'\subset P_n$ contradicting the choice of the iterate $f^k(z)$. By the
unbranching property, the pullback $G_n$ extends to the domain $M_n$.
Let $P_n(z)=G_n(P_n)$. We want to show, of course,
that the Euclidean diameters
of $P_n(z)$ tend to zero. For this, let us consider a domain $M_n'\subset M_n$
bounded by a core curve of the annulus $A_n$. Then
$\max_{y\in \partial M'_n}|f^k(z)-y|/\min_{y\in \partial M'_n}|f^k(z)-y|\le
C(m)$,$n=1,2,\dots$ (see e.g. \cite{McM}). Introduce domains $E_n=G_n(M_n')$.
Since the modulus of annulus $M_n\setminus M_n'$ is $m/2$, by Koebe
distortion theorem, $\max_{y\in \partial E_n}|z-y|/\min_{y\in
\partial E_n}|z-y|\le C_1(m)$, $n=1,2,\dots$. If we assume that
$\mbox{diam} P_n\ge d>0$ for $n=1,2,\dots$,
then $\min_{y\in \partial E_n}|z-y|\ge d/2C_1=r>0$,
i.e., the disc $D_z(r)\subset E_n$. Hence, $f^{k_n}(D_z(r))\subset M_n'$,
for $n\to \infty$ and $k_n\to \infty$. This is a contradiction with the
non-normality of the family $f^n$ at $z\in J(f)$. Thus, $\bigcap_{n>0}
P_n(z)=\{z\}$ and so $J(f)$ is again locally connected at $z$.

\sect{The proof of the Main Theorem in the infinitely renormalizable
case when $\ell=2$}
\label{sec7}
In this section we consider a renormalizable $f(z)=z^2+c_1$
with a periodic interval $V$ of the period $s$.
As before, we may assume that $R_V$ has a high return
because otherwise $f$ has a periodic attractor.
Again we shall delay dealing
with the case that $f$ is also renormalizable of period $s/2$
until Section~\ref{secll}.
Hence, because of Lemma~\ref{str5},
the space from Lemma~\ref{str1} is at least $0.6$.
Note however, that the bound we obtain for
$K^\ast_2(0.6)$ is equal to $1.19371...>1$ and
we cannot use the method of the previous section in the
quadratic case.
Therefore, we construct a domain of
renormalization $\Omega$ for $f^s$ which is
different from $D_\ast(V)$, but with a diameter depending on $|V|$
only, so that the renormalizations shrink to zero together with $V$.

\kies{
\begin{figure}[htp]
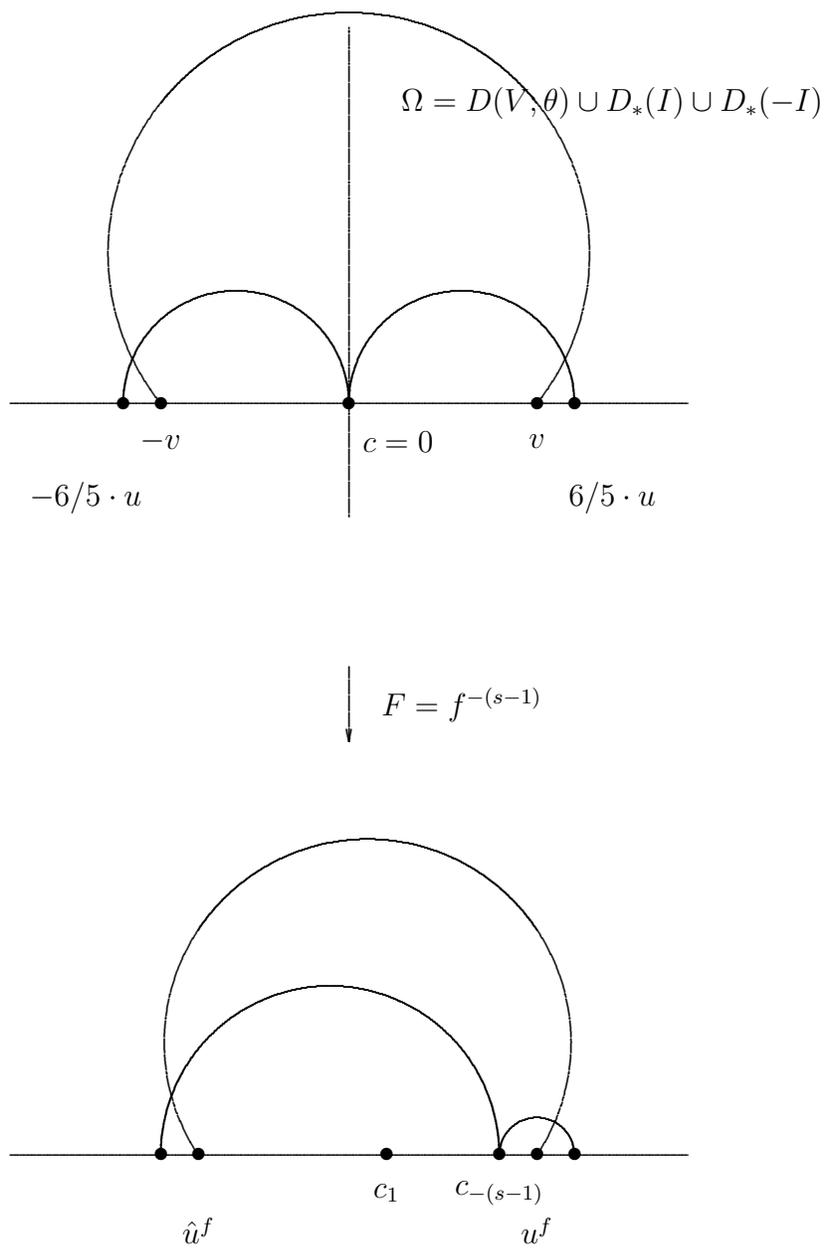
 \hfil
\beginpicture
\dimen0=0.5cm
\setcoordinatesystem units <\dimen0,\dimen0>
\setplotarea x from -9 to 9, y from -3 to 9
\setlinear
\plot -9 0 9 0 /
\plot 0 -3  0 10 /
\setsolid
\circulararc 180 degrees from 5 0  center at  0 4
\circulararc -180 degrees from -5 0  center at  0 4
\circulararc 180 degrees from 0 0  center at  -3 0
\circulararc -180 degrees from -6 0  center at  -3 0
\circulararc -180 degrees from   0 0  center at  3 0
\circulararc  180 degrees from   6 0  center at  3 0
\put {\smalll $c=0$} at 1.3 -1
\put {\smalll $v$}   at 5 -1
\put {\smalll $-v$} at -5 -1
\multiput {\smalll $\bullet$} at 6 0  5 0  0 0  -5 0 -6 0 /
\put {\smalll $6/5\cdot u$} at 7 -2.5
\put {\smalll $-6/5\cdot u$} at -7 -2.5
\put {\smalll $\Omega=D(V;\theta)\cup D_*(I)\cup D_*(-I)$} at 7 8
\arrow <5pt> [0.2,0.4] from 0 -7 to 0 -9
\put {\smalll $F=f^{-(s-1)}$} at 3 -8
\setcoordinatesystem units <\dimen0,\dimen0> point at 0 20
\setplotarea x from -9 to 9, y from -7 to 7
\setlinear
\plot -9 0 9 0 /
\circulararc 180 degrees from 5 0  center at  0.5  3
\circulararc -180 degrees from -4 0  center at  0.5 3
\circulararc 180 degrees from 6 0  center at  5  0
\circulararc -180 degrees from 4 0  center at  5  0
\circulararc 180 degrees from 4 0  center at  -0.5  0
\circulararc -180 degrees from -5 0  center at  -0.5  0
\put {\smalll $c_1$} at 1 -1
\multiput {\smalll $\bullet$} at 6 0 5 0 4 0 1 0 -4 0 -5 0 /
\put {\smalll $c_{-(s-1)}$}   at 4 -1
\put {\smalll $u^f$}   at 5 -2
\put {\smalll $\hat u^f$} at -4 -2
\endpicture
\caption[ ]{{\smalll The region $\Omega$ and its preimage.}}
\end{figure}


}

As before, denote $D(J;\theta)$ the Poincar\'e
neighbourhood of a real interval $J$,
with external angle $\theta\in (0,\pi/2]$.
We define the domain $\Omega=\Omega(\theta)$
as the Poincare neighborhood $D(V;\theta)$
of the periodic interval $V$
(with the angle $\theta$ to be specified later on),
united with two discs: $D_\ast(I)$ and $D_\ast(-I)$ with
$I=(0,6/5\cdot u)$, and $u$ is the periodic endpoint
of $V$, (i.e., $f^s(u)=u$).
We may and shall assume that $u>0$.
Since $R_V$ has a high return, one has
that $c_s\in [-u,0]$.

\hbox to \hsize{\hss\unitlength=1.3mm
\beginpic(70,60)(0,0) \let\ts\textstyle
\put(0,5){\line(1,0){67}}
\put(30,4.7){\line(1,0){37}}
\put(0,25){\line(1,0){67}}
\put(0,45){\line(1,0){67}}
\put(0,44.7){\line(1,0){40}}
\put(25,-7){{\it The intervals $I,-I$.}}
\put(10,6){\line(0,-1){2}} \put(9,1){$b$}
\put(20,6){\line(0,-1){2}}\put(19,1){$\hat u^f$}
\put(30,6){\line(0,-1){2}}\put(29,1){$c_1$}
\put(40,6){\line(0,-1){2}} \put(39,1){$c_{-(s-1)}$}
\put(50,6){\line(0,-1){2}}
\put(60,6){\line(0,-1){2}} \put(59,1){$u^f$}
\put(67,6){\line(0,-1){2}}
\put(10,26){\line(0,-1){2}} \put(9,21){$a'$}
\put(20,26){\line(0,-1){2}} \put(19,21){$-u$}
\put(30,26){\line(0,-1){2}} \put(29,21){$c_s$}
\put(40,26){\line(0,-1){2}} \put(39,21){$c=0$}
\put(50,26){\line(0,-1){2}}
\put(60,26){\line(0,-1){2}} \put(59,21){$u$}
\put(67,26){\line(0,-1){2}}\put(65,21){$6/5\cdot u$}
\put(19,29){$-I$} \put(43,29){$I=(0,6/5\cdot u)$}
\put(0,10){\vector(0,1){9}}

\put(10,46){\line(0,-1){2}} \put(9,41){$a$}
\put(20,46){\line(0,-1){2}} \put(19,41){$u^f$}
\put(30,46){\line(0,-1){2}} \put(29,41){$c_{s+1}$}
\put(40,46){\line(0,-1){2}} \put(39,41){$c_1$}
\put(0,31){\vector(0,1){9}}
\put(-8,15){$f^{s-1}$} \put(-8,35){$f$}
\endpic\hss}
\vskip1.6cm

\begin{remark}
\label{rem91}
Consider the map $f^{s-1}\colon \hat U\to V$.
As in Lemma~\ref{str5}, let $l$ be the maximal interval with
a unique common endpoint with $\hat U$ and outside $[c_1,c_2]$
such that $f^{s-1}|l$ is monotone. From Lemma~\ref{str5}
we have that $|f^{s-1}(l)|\ge \sqrt{1+0.6}\cdot |V_+|$
where $V_+$ is one component of $V\setminus \{c\}$.
The number $6/5$ is chosen for the following reason:
$6/5>K^\ast_2(0.6)=1.19371..$, but $6/5\sma \sqrt{1+0.6}$.
Because of Lemma~\ref{str5}, the latter
inequality shows that the pullback $F$ of $f^{s-1}\colon \hat U\to V$
has a monotone extension to $\widetilde I=I\cup V\cup (-I)$.
Let $F$ also denote the extension of this pullback to
$\cz_{\widetilde I}=(\cz \setminus \rz)\cup \widetilde I$.
The first inequality we shall need in the
proof of Corollary~\ref{corr91} below.
\end{remark}

Our aim is to prove that $f^{-1}\circ F(\Omega)$ is a proper domain inside
of $\Omega$. For this, it is enough to prove

\begin{prop}
\label{gnlq2}
$f^{-1}\circ F(\partial \Omega)\subset \Omega$.
\end{prop}
\pr
First of all, we observe that
$f^{-1}\circ F(D_\ast(I))$ consists of two components $D^1$ and $-D^1$,
where $D^1$ is a proper domain in $D_\ast(I)$. 
Note that there exists an interval
containing $c_1$ which is mapped diffeomorphically by $f^{s-1}$
onto $\widetilde I$. Therefore,
$f^{-1}\circ F$ maps $I$ homeomorphically
into itself, because $u$ is a repelling fixed point of $f^s$
(here we use that $f^s$ has negative Schwarzian derivative).
Moreover, $c_s=f^{s-1}(c_1)$ is not
in $I$ and therefore $F(D_*(I))$ does not contain $c_1$.
Hence
$f^{-1}\circ F(D_\ast(I))$ consists of two components $D^1$ and $-D^1$,
where $D^1$ is a proper domain in $D_\ast(I)$.
Thus, the corresponding branch of $f^{-1}\circ F(\Omega)$
maps $\cz_I$ univalently into itself. By symmetry, we get
also $-D^1\subset D_\ast(-I)$.
Thus we have shown that
\beq
f^{-1}\circ F(D_\ast(I))\subset \Omega.
\label{sub1}
\eeq

Next we want to show that
the pullback of $D(V;\theta)$ is inside $\Omega(\theta)$ provided
we choose $\theta$ conveniently. For this we shall
first consider the pullbacks through the map $P(z)=z^2$.
Fix $K\ge 1$. In the next lemma
we are going to compare
$P(D((-1,1);\theta))$
with the Poincar\'e disc $D((-K,1);\theta)$.
The latter disc is a `scaled-up' version
of $D((\hat u^f,u^f);\theta)\ni c_1$.

\begin{lemma}
\label{lem71}
Let $K>1$. There exists $\theta_0=\theta_0(K)>0$ such that,
for all
$\theta\in (0,\theta_0)$, the boundaries of $P(D((-1,1);\theta))$ and
$D((-K,1);\theta)$ intersect each other in
$Z(K,\theta)$ and its complex conjugate.
Furthermore,
$$Z(K,\theta)\to K^2\in \rz,$$
as $\theta\to 0$.
Hence, the difference $\Delta(K,\theta)=D((-K,1);\theta)
\setminus P(D((-1,1);\theta))$ tends to the interval $[1,K^2]$, as $\theta\to
0$.
\end{lemma}
\pr Consider a possible intersection point
$$z_1\in
P(\partial D((-1,1);\theta))\cap \partial D((-K,1);\theta).$$
That is, $z_1=P(z_2)=z_2^2$ with
$$z_2\in \partial D((-1,1);\theta)).$$
Since these sets are symmetric with respect to the real axis
and since $D((-1,1);\theta)$ is also symmetric with respect
to the imaginary axis, we may consider the case
that $z_1$ is in the upper half plane,
and that $z_2$ is in the first quadrant.

Since $z_2\in \partial D((-1,1);\theta)$,
\beq
z_2=1+{i\exp(i\theta)\over \sin(\theta)}(1-\exp(i\alpha))
\label{B1}
\eeq
where $\alpha\in (0,2\pi-\theta)$ is the angle between the vectors
$z_2-C,1-C$, with $C$ the centre of the circle $D((-1,1);\theta)$.
Similarly, since $z_1\in \partial D((-K,1);\theta)$,
\beq
z_1=1+{K+1\over 2}\cdot {i\cdot \exp(i\theta)\over \sin(\theta)}\cdot
(1-\exp(i\phi)).
\label{B2}
\eeq
Taking the square of (\ref{B1}),
$$z_2^2=1+{2\exp(i\theta)\over \sin^2(\theta)}(1-\exp(i\alpha ))
\left\{i\sin(\theta)-{1\over 2}\exp(i\theta)(1-\exp(i\alpha ))\right\}.$$
We have in the brackets $\{...\}$ term:
$$i\sin(\theta)-{1\over 2}\exp(i\theta)(1-\exp(i\alpha ))=$$
$${1\over 2}(\exp(i\theta)-\exp(-i\theta)-\exp(i\theta)+\exp(i(\theta+\alpha )))=$$
$${\exp(i\alpha )\over 2}(\exp(i\theta)-\exp(-i(\theta+\alpha ))).$$
So,
$$z_2^2=1+{2\exp(i\theta)\over \sin^2(\theta)}(1-\exp(i\alpha )){\exp(i\alpha )\over 2}
(\exp(i\theta)-\exp(-i(\theta+\alpha )))=$$
$$1+{2\exp(i\theta)\over 2\sin^2(\theta)}\exp(i\alpha )
\{\exp(i\theta)-\exp(-i(\theta+\alpha ))-\exp(i(\theta+\alpha ))
+\exp(-i\theta)\}=$$
$$1+{\exp(i(\theta+\alpha ))\over \sin^2(\theta)}
\{2\cos(\theta)-2\cos(\theta+\alpha )\}.$$

If $z_2^2=z_1$, then we compare the last expression with
(\ref{B2}), cancel $1$ in both hand-sides
and then divide them by $\exp(i\theta)$, and multiply by $\sin(\theta)$,
and after that separate Re and Im parts:
$${2(\cos(\theta)-\cos(\theta+\alpha ))\over \sin(\theta)}\cos(\alpha )=
{K+1\over 2}\sin(\phi),$$
$${2(\cos(\theta)-\cos(\theta+\alpha ))\over \sin(\theta)}\sin(\alpha )=
{K+1\over 2}(1-\cos(\phi)).$$
Now divide the second equality to the first one:
$$\tan(\alpha )=\tan(\phi/2).$$
Here $0\sma \phi/2\sma \pi$, so either $\alpha =\phi/2$, or $\alpha =\phi/2+\pi$.
The latter case is impossible, since $z_2$ is in the first quarter, i.e.
$0\sma \alpha +\theta\sma \pi/2$.

Thus, 
$$\alpha =\phi/2.$$
Now we substitute  $\phi=2\alpha $ in the equality, say, for the real parts:
$${2(\cos(\theta)-\cos(\theta+\alpha ))\over \sin(\theta)}\cos(\alpha )=
{K+1\over 2}\sin(2\alpha ),$$
or 
$$4\sin(\alpha /2)\sin(\theta+\alpha /2)\cos(\alpha )=
{K+1\over 2}2\sin(\alpha )\cos(\alpha )\sin(\theta),$$
or
$$\sin(\theta+\alpha /2)={K+1\over 2}\cos(\alpha /2)\sin(\theta),$$
or
$$\cos(\theta)\sin(\alpha /2)={K-1\over 2}\cos(\alpha /2)\sin(\theta),$$
or, at last,
$$\tan(\alpha /2)={K-1\over 2}\tan(\theta).$$

If $\theta$ here is small, we find a unique $\alpha$ in the admissible
interval $(0,\pi-\theta/2)$ for the angle $\alpha$ (this is because
$\alpha=\phi/2\in (0,\pi-\theta/2)$). So there is the unique point of
intersection of the curves $\partial P(D((-1,1);\theta))$ and
$\partial D((-K,1);\theta)$ in the upper halfplane. If $\theta \to 0$, then
$\phi\sim 2(K-1)\theta$ and $z_1\sim K^2$. The rest of the lemma follows
easily.
\qed

\begin{corr}
\label{corr91}
If $\theta$ is small enough, then
\beq
f^{-1}\circ F(\partial D(V;\theta)) \subset \Omega(\theta).
\label{sub2}
\eeq
\end{corr}
\pr
Take $K=K^\ast_2(0.6)=1.19371...$. By
the Schwarz contraction principle, see Lemma~\ref{scp},
$F(D(V;\theta))$ is contained in
$D(\hat I;\theta)$, where the interval $\hat I=(\hat u^f,u^f)$
is around $c_1$, with $|c_1-\hat u^f|\le K|u^f-c_1|$. By rescaling and
the previous lemma, for all $\theta$ less than some positive absolute constant
$\theta_1$ the closure of the domain $F(D(V;\theta))$ is
inside $f(\Omega(\theta))$.
Here we have used that the previous lemma implies that
$f^{-1}(F(D(V;\theta)))\setminus D(V;\theta)$ converges
to the interval $\pm (u,K\cdot u)$ as $\theta\to 0$.
Since $K<6/5$, see Remark~\ref{rem91}, this implies that this difference
set is contained in $D_*(I)\cup D_*(-I)$.
 
\qed

\medskip

\noindent
{\it Conclusion of the proof of Proposition~\ref{gnlq2}.}
To complete the proof, it remains to show that
for all angles $\theta$ less than some other positive constant $\theta_2$
the domain $F(D_\ast(-I))$ is contained properly inside $P(\Omega(\theta))$.
Note that $P(-I)=P(I)=(c_1,a)$, where
$|c_1-a|=(6/5)^2|c_1-u^f|\sma (1+0.6)\cdot |c_1-u^f|$, i.e. $a\in L$.
Hence, by Lemma~\ref{str1},
the interval $F(-I)$ is inside an interval $(b,u^f)$,
where 
$$|b-c_1|/|u^f-c_1|
=(|b-c_1|/|a-c_1|)\cdot (5/6)^2\sma K^\ast_2(y_0)\cdot (5/6)^2:=K_1\sma 1$$
where, using the notation of Lemma~\ref{str1},
$y_0=|a-c_1|/|T|=(6/5)^2/(1+0.6)\in (0.6,1)$.
It follows, $F(D_*(-I))$ is inside of
the ball
$D_\ast(b,u^f)$. By rescaling, we need to show that $P(D((-1,1);\theta))$
contains a fixed ball $D_\ast(-K_1,1)$. If $\theta$ is small, it is a not
difficult exercise. Hence
$$f^{-1}\circ F(D_\ast(-I)) \subset f^{-1}(D_*(b,u^f))
\subset D_*(U;\theta).$$
Together with (\ref{sub1}) and (\ref{sub2}) this implies
Proposition~\ref{gnlq2}.
\qed

With the constructed sequence of the domains of renormalizations $\{\Omega\}$
we end the proof of the Main Theorem in the infinitely renormalizable
case for degree two, simply repeating the proof
of this theorem for the larger degrees (see the end of the
previous section).

Let us now prove Theorem A for $\ell=2$. For this we make use of
Lemma~\ref{lsul}
(or rather its proof): for every small enough $\varepsilon>0$
there exists a constant $C, 0<C<1$ (depending only on $\varepsilon$), 
and a point $\widetilde u$ in the interval $I\setminus V$
such that the image $f^s(\widetilde u)$
lies in the $\varepsilon$-neighbourhood
of the point ${6\over 5}u$ and
\beq
|\widetilde u-c|\le C|f^s(\widetilde u)-c|.
\label{75}
\eeq
With this $\varepsilon$ small enough (but fixed) and the corresponding
point $f^s(\widetilde u)$, which replaces the previous point ${6\over 5}u$,
we can construct the domain $\Omega$ (see Remark~\ref{rem91}), with the same
angle $\theta_0$. Then the modulus of the annulus
$\Omega\setminus f^{-1}\circ F(D(V;\theta_0)\cup D_*(-I))$ is bounded from
below by
a positive absolute constant. On the other hand, by (\ref{75}),
two preimages $f^{-1}\circ F(D_*(I))$ are also on a proportionally
definite distance from the boundary of $\Omega$, and we obtain Theorem A.
\qed

 %


\sect{The proof of the Main Theorem when $\omega(c)$ is not minimal
and in the Fibonacci case}

In the remainder of the paper we shall deal with
the non-renormalizable case (except in Section 13, where we shall
finish the proof of Theorem A
in the infinitely renormalizable case when period doubling
occurs). Firstly, if the $\omega$-limit set $\omega(c)$
of the critical point $c=0$ is not minimal
then it very easy to see that the Julia set is
locally connected. 
To see this, note that if $\omega(c)$ is not minimal
then it contains a point $x$ whose forward orbit stays away from
the critical point $c$. Hence this forward orbit lies in
a hyperbolic set. Therefore the
Yoccoz puzzle-pieces $P_n(x)$ containing $x$ shrink down in diameter to zero.
Since $x\in \omega(c)$,
the forward orbit of $c$ enters these pieces and it follows that
all the puzzle-pieces tend to zero in diameter. Since
the intersection of the Julia set with puzzle-pieces is connected
the result follows.

Now we shall prove that
the Julia set of a Fibonacci map
of the form $f(z)=z^\ell+c_1$ with $c_1$ is real
is locally connected. We should emphasize that this
result also follows from Theorem B. However, since the Fibonacci
map is often thought of as the `bad case', we want
to show explicitly that the careful estimates obtained
in \cite{SN} imply that the proof of local connectivity in this case
is in fact very easy. Let us write as before
$\tau(z)=-z$. Let us remind that
a Fibonacci map is a map defined by the
following property:
For $i\ge 0$ and $x\in \rz$, let $x_i=f^i(x)$ and
choose $x_{-i}\in f^{-i}(x)$ so that
the interval connecting this point to $c=0$ contains no other points
in the set $f^{-i}(x)$. Note that if $c$ is not a periodic point
there are always precisely two such points $c_{-i}$
(which are symmetric with respect to each other).
Let $S_0=1$ and define $S_i$ inductively by
$$S_i=\min\{k\ge S_{i-1};\, c_{-k}\in (c_{-S_{i-1}},\hat c_{-S_{i-1}})\}.$$
$f$ is called a {\it Fibonacci} map
if the sequence $S_i$ coincides with
the Fibonacci numbers: $S_0=1$, $S_1=2$ and $S_{k+1}=S_k+S_{k-1}$, i.e.,
the sequence $1,2,3,5,8,\dots$.
Moreover,
let us define inductively a sequence of points $u_n$ as follows.
Let $u_0$ be the orientation reversing fixed point
$q$ of $f$ and let us define $u_{n+1}$ to be the nearest point
to $c$ with
$$u_{n+1}\in f^{-S_n}(u_n)$$
so that $u_{n+1}$ is on the same side of $c$ as $c_{S_{n+1}}$.
In particular, $u_1=\hat u_0=\hat q$.
We shall use

\begin{prop} \cite{SN}
For each even integer $\ell\ge 2$,
there exists a sequence of standard discs $D_n$
centred at the critical point $0$ and
relatively compact topological discs $D^0_n, D^1_n$ in $D_n$, such that:
\begin{itemize}
\item
A trace of $D_n$ on the real line is ended by two symmetric preimages
$u_{n-1}, \tau(u_{n-1})$ of
an orientation reversing fixed point $q$ of $f$.
\item
The sequence of discs $D_n$ shrinks to zero.
\item
For each $n$ big enough, the map
$$R_n:(D^0_n\bigcup D^1_n)\to D_n$$ defined by
$$R_n(z)=\{f^k(z)\st k>0 \ \hbox{is minimal with} f^k(z)\in D_n\}$$
is $l$-polynomial-like with $l$-multiple critical point zero.
\item
Moreover, the critical point of $R_n$ lies
in the Julia set of $R_n$.
\end{itemize}
\end{prop}
\pr For the proof see \cite{SN} (this theorem
is proved there for each even
$\ell\ge 2$). For $\ell=2$ this result
is also proved in \cite{LM}.
\qed

\noindent
{\it Proof of the Main Theorem in the Fibonacci case.}
There is a critical Yoccoz piece $P_n$ containing
$u_{n-1},\tau(u_{n-1})$ in its boundary.
Let $Q^{-1}$ be the extension of the
$R^{-1}_n$ from $D_n\bigcap P_n$ to $P_n$.
Let $P^0_n$, $P^1_n$ be the images of $P_n$
under this map $Q^{-1}$.  Then $Q\colon  (P^0_n\bigcup
P^1_n)\to P_n$ is again
$l$-polynomial-like map. Observe that $R_n=Q$ on the real line. Hence,
the critical point zero belongs to both Julia sets $J_{R_n}$ and $J_Q$.
As we have proved in Proposition~\ref{gnp2}, we have
$F_{R_n}=F_Q$ and that there exists a component of
some preimage of $Q^{-i}(P_n)$ which contains $c$
and which is contained in $D_n$. By construction,
$J(f)\cap P_n$ is connected. Hence $J(f)\cap Q^{-i}(P_n)$
is connected and for some $i$ and some component $C_n\ni c$
of $Q^{-1}(P_n)$ is contained in $D_n$.
Hence $C_n$ is an open neighbourhood of $c$
which is contained in $D_n$ and such that
$C_n\cap J(f)$ is connected.
Since the sequence $u_n$ tends to $c=0$,
the diameter of $D_n$ tends to zero and we get
that the diameter of $C_n$ tends to zero also.
Hence $J(f)$ is locally connected in $c=0$.

To prove the local connectivity at any other point $z\in J_f$, we can repeat
arguments for quadratic case (see \cite{Mil}),
which work in our case as well.
If the orbit of $z$ does
not hit a critical piece, then the pieces around $z$ shrink to $z$ by
contraction principle. Let now zero be an accumulation point of the orbit
of $z$. Consider the annuli given by the Yoccoz pieces.
First, note that the sum over the all depths $d=1,2,\dots$
of the moduli of the annuli $A_d(0)$ around zero is infinite, just because
the critical pieces tend to the point. Fix $d$ and find the first iterate
$z_j$ of $z$ that hits the critical piece at depth $d+1$. Then an annulus
of the puzzle around $z$ at depth $d+j$ is isomorphic to $A_d(0)$. Furthermore,
distinct values of $d$ give distinct values of $d+j$. Hence, the sum of the
moduli of annuli around $z$ is infinite, as we needed.
\qed

We should note that the Julia set of the Fibonacci polynomial
$z\mapsto z^\ell+c_1$ with $c_1$ real has positive Lebesgue measure
when $\ell$ is large, see \cite{SN}. It follows that there exists
Julia sets which are locally connected but have positive Lebesgue
measure.

\sect{The proof of the Theorem B for $\ell$ large}
\label{sec9}

In this section we prove Theorem B and the Main Theorem
in some non-renormalizable cases when $\ell$ is large.

\begin{theo}
There exists $\ell_0$ as follows.
Let $f(z)=z^\ell + c_1$ with $\ell$ an even integer and $c_1$ real be a
non-renormalizable polynomial such that the limit set $\omega(c)\ni c$
is minimal and $f$ has
infinitely many times a high return in the partition
given by the Yoccoz puzzle on the real line.
Then the Julia set
of $f$ is locally connected provided $\ell>\ell_0$.
\end{theo}

For a definition of the notion of a high return,
see the end of the introduction.
We should note that the proof of this theorem
also holds for every infinitely renormalizable $f$
with $\ell$ is large enough (and thus giving an alternative
proof of Main Theorem in the infinitely renormalizable
case when $\ell$ is large).

As before, let $W$ be a symmetric interval with nice boundary points,
let $R_W$ be the first return map to $W$ and
let $V=[v,\tau(v)]$ be the domain of $R_W$ containing $c$.
Similarly, let $R_V$ be the first return map to $V$
and $U=[u,\tau(u)]$ the component of the domain of $R_V$ which contains
$c$. Let $\hat U=[\hat u^f,u^f]$ be the component of $R_V$ containing
the critical value $c_1$. Let $s$ be so that
$R_V|\hat U=f^{s-1}$. In this section
we will assume that $R_V$ has a high return,
i.e., $R_V(U)=f^s(U)\ni c$.

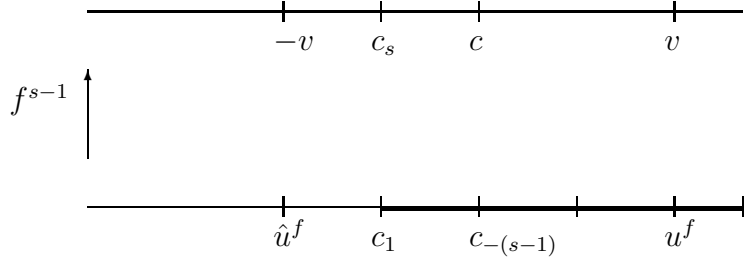
\hbox to \hsize{\hss\unitlength=1.3mm
\beginpic(70,35)(0,0) \let\ts\textstyle
\put(0,5){\line(1,0){67}}
\put(30,4.7){\line(1,0){37}}
\put(0,25){\line(1,0){67}}
\put(20,6){\line(0,-1){2}}\put(19,1){$\hat u^f$}
\put(30,6){\line(0,-1){2}}\put(29,1){$c_1$}
\put(40,6){\line(0,-1){2}} \put(39,1){$c_{-(s-1)}$}
\put(50,6){\line(0,-1){2}}
\put(60,6){\line(0,-1){2}} \put(59,1){$u^f$}
\put(67,6){\line(0,-1){2}}
\put(20,26){\line(0,-1){2}} \put(19,21){$-v$}
\put(30,26){\line(0,-1){2}} \put(29,21){$c_s$}
\put(40,26){\line(0,-1){2}} \put(39,21){$c$}
\put(60,26){\line(0,-1){2}} \put(59,21){$v$}
\put(0,10){\vector(0,1){9}}
\put(-8,15){$f^{s-1}$}
\endpic\hss}
\vskip1.6cm

We cannot use Lemmas~\ref{str4}-\ref{str5}
since $f$ above is not renormalizable, however we can Proposition~\ref{sp13}
(which also holds for renormalizable $f$).
Let us state it quickly again.
Let $T_0$ be a minimal interval containing $f(V)$ and its immediate
neighbour among the disjoint intervals $f^2(V),\dots,f^{s'}(V)$. Write
$L_0=T_0\setminus f(V)$. By Lemma~\ref{str3}, there exists an interval
$l$ on either side of $\hat U$ which has a unique common point with $\hat U$
such that $f^s\colon l\to L_0$ is one-to-one and by Proposition~\ref{sp13}:

\begin{lemma}
\label{lemmaC1}
$$|L_0|\ge {1\over 3}|f(V)|.$$
\end{lemma}

Given the interval $V=(-v,v)$ we construct an $\ell$-polynomial-like map sitting
inside the domain $\Omega=\Omega(\ell,V)$ defined as
$$\Omega=D(V;\theta)\bigcup D(I;\theta)\bigcup D(-I;\theta),$$
where
$\theta=\theta_0$ is some absolute constant (angle) to be determined later on,
the interval
$$I=(0,v+{\log(11/10)\over l}v).$$

\kies{
\begin{figure}[htp]
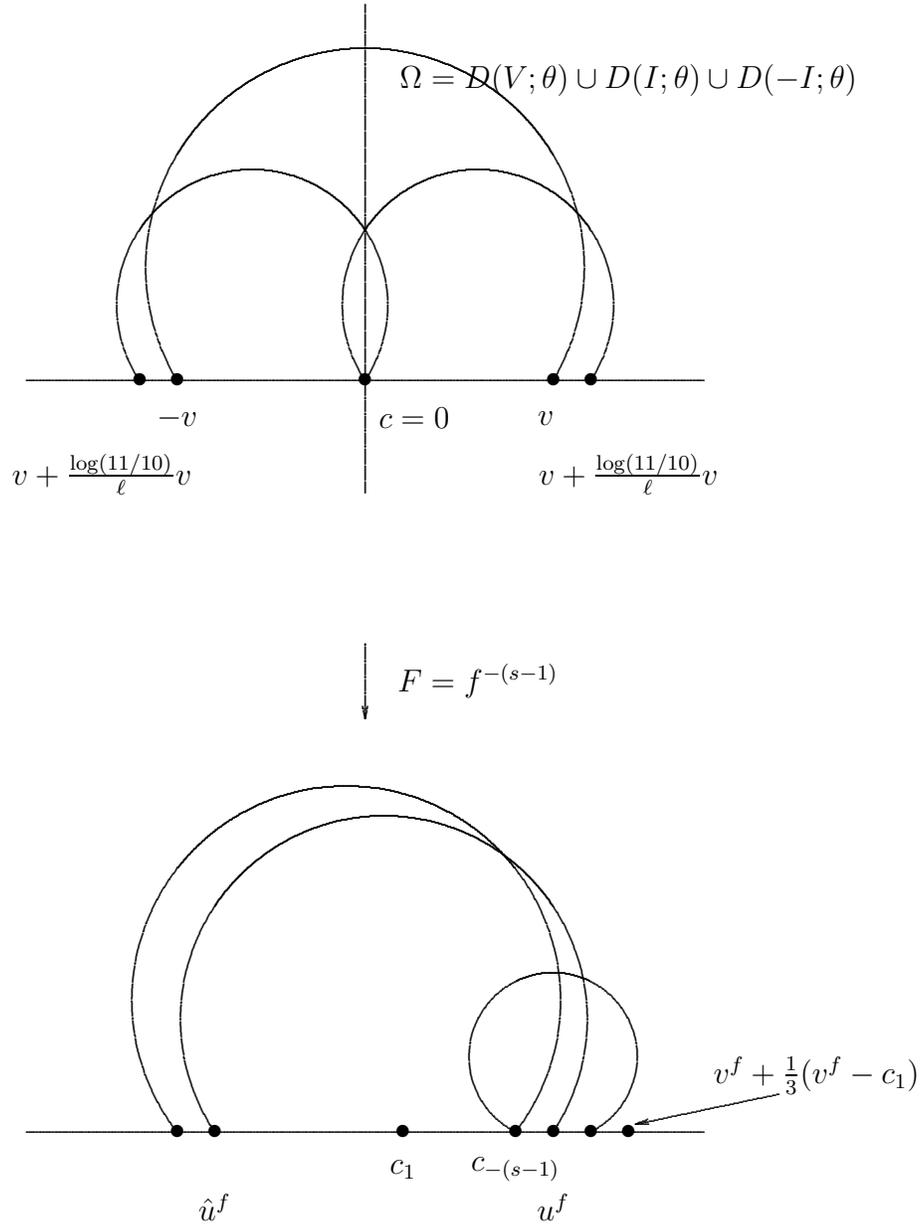
 \hfil
\beginpicture
\dimen0=0.5cm
\setcoordinatesystem units <\dimen0,\dimen0>
\setplotarea x from -9 to 9, y from -3 to 9
\setlinear
\plot -9 0 9 0 /
\plot 0 -3 0 10 /
\setsolid
\circulararc 180 degrees from 5 0  center at  0 3
\circulararc -180 degrees from -5 0  center at  0 3
\circulararc 180 degrees from 0 0  center at  -3 2
\circulararc -180 degrees from -6 0  center at  -3 2
\circulararc -180 degrees from   0 0  center at  3 2
\circulararc  180 degrees from   6 0  center at  3 2
\put {\smalll $c=0$} at 1.3 -1
\put {\smalll $v$}   at 4.8 -1
\put {\smalll $-v$} at -5 -1
\multiput {\smalll $\bullet$} at 6 0  5 0  0 0  -5 0 -6 0 /
\put {\smalll $v+\frac{\log(11/10)}{\ell}v$} at 7 -2.5
\put {\smalll $v+\frac{\log(11/10)}{\ell}v$} at -7 -2.5
\put {\smalll $\Omega=D(V;\theta)\cup D(I;\theta)\cup D(-I;\theta)$} at 7 8
\arrow <5pt> [0.2,0.4] from 0 -7 to 0 -9
\put {\smalll $F=f^{-(s-1)}$} at 3 -8
\setcoordinatesystem units <\dimen0,\dimen0> point at 0 20
\setplotarea x from -9 to 9, y from -7 to 7
\setlinear
\plot -9 0 9 0 /
\circulararc 180 degrees from 5 0  center at  0.5  3
\circulararc -180 degrees from -4 0  center at  0.5 3
\circulararc 180 degrees from 6 0  center at  5  2
\circulararc -180 degrees from 4 0  center at  5  2
\circulararc 180 degrees from 4 0  center at  -0.5  3.5
\circulararc -180 degrees from -5 0  center at  -0.5  3.5
\put {\smalll $c_1$} at 1 -1
\multiput {\smalll $\bullet$} at 7 0 6 0 5 0 4 0 1 0 -4 0 -5 0 /
\put {\smalll $v^f+\frac{1}{3}(v^f-c_1)$}   at 12 1.5
\arrow <5pt> [0.2,0.4] from 11 1 to 7.2 0.2
\put {\smalll $c_{-(s-1)}$}   at 4 -1
\put {\smalll $u^f$}   at 5 -2
\put {\smalll $\hat u^f$} at -4 -2
\endpicture
\caption[ ]{{\smalll $\frac{c_1-\hat u^f}{v^f-c_1}<1.3\text{ when }\ell
>>1$.}}
\end{figure}


}

Let $F\colon V\to \hat U$ be the inverse to
the map $f^{s-1}\colon \hat U\to V$.
Then $F$ extends to a univalent map on a maximal interval $T$ containing $V$
and then to domain $\cz_T$. The first observation is that
for all sufficiently large degrees $\ell$,
the interval 
$$\widetilde I=I\cup V\cup (-I)=\Omega\cap \rz$$
is inside the interval $T$. This is because
$$(1+{\log(11/10)\over \ell})^\ell\sim {11\over 10}<1+1/3,$$
that is $f(\widetilde I)\subset (c_1,v^f)\cup L_0$, by Lemma~\ref{lemmaC1}.
In the formula above the symbol $\sim$ means that the
left hand side converges to the right hand side as $\ell\to \infty$.
We shall use this convention throughout this section.

We are going to prove

\begin{prop}
\label{propC1}
There exists $\theta_0>0$, such that for all
sufficiently big $\ell$,
$$f^{-1}\circ F(\partial \Omega)\subset \Omega.$$
\end{prop}

Before proving this proposition we show

\begin{prop}
Proposition~\ref{propC1} implies the Main Theorem in the non-renormalizable
high case when $\ell$ is sufficiently large.
\end{prop}
\pr It is enough to construct the $\ell$-polynomial-like mapping
inside the domain $\Omega$.
Since $f\colon \omega(c)\to \omega(c)$ is minimal,
each point in $\omega(c)$ eventually is mapped into
$V$ and therefore is in the domain of definition
of the map $R_V$. 
There exists a finite collection of disjoint
intervals $I^0=U$ and $I^1,\dots,I^i$ in $V$,
which form the domain of definition
of the map $R_V$ (see Proposition~\ref{pqp} and its proof).
More precisely,
every $x\in \omega(c)\cap V$ belongs to some $I^j$ or to $U$,
the map $R_V\colon I^j\to V$
is a diffeomorphism for $j=1,2,\dots,i$,
and $U\in c$ with $R_V|_{U}(U)=f^s(U)\ni
c$ (the high return). Given $j=1,\dots,i$, there exists an interval $\hat I^j$
containing $I^j$ such that $R_V\colon I^j\to V$
extends to a diffeomorphism from
$\hat I^j$ onto the interval $\widetilde I$ (from the definition of the domain
$\Omega$). Indeed, if $l^j\supset I^j$ is a maximal interval on which $R_V$
is monotone, then $R_V(l^j)$ contains $V$ and its immediate neighbour (from
the collection of intervals $f^j(V)$)
on either side of $V$. Hence, $R_V(l^j)\subset f^{-1}(T_0)\cap
\rz$. Thus, $\hat I^j\subset l^j$. Since $f$ has no attracting periodic
orbit,
$\hat I^j$ is a subset of either the right part $I$ of $\widetilde I$ or its
left part $-I$. Let, for example, $\hat I^j\subset I$.
By Lemma~\ref{scp}, there exists a domain $\Omega_j$ inside
$D(I^j;\theta)\subset D(I;\theta)\subset \Omega$ which is mapped
diffeomorphically by a map $R^j$ (an iteration of $f$)
onto $\Omega$. We have constructed $\Omega_j$
for each $j=1,2,..,i$, and each $\theta\in (0/\pi/2)$. Assuming 
Proposition~\ref{propC1},
we fix the angle $\theta_0$ and find a domain
$\Omega_0\subset \Omega$ such that $R^0=f^s\colon \Omega_0\to
\Omega$ is a proper $\ell$-cover. Since the domains $\Omega_j, j=0,1,\dots,i$
may intersect each other, we modify so that they become
$\ell$-polynomial-like.
To do this, let us consider a Yoccoz piece $P_V$ containing the ends of
the interval $V$ on its boundary. Let $D$ be a component of $P_V\cap \Omega$
containing $V$. For every $j=0,1,\dots,i$, there exists a domain $D_j\subset
\Omega_j$ such that $R_j\colon D_j\to D$ is a diffeomorphism if $j>0$, and
an $\ell$-cover if $j=0$. Since all $D_j$ lie in different Yoccoz pieces,
we obtain the $\ell$-polynomial-like map.
\qed

\noindent
{\em Proof of Proposition~\ref{propC1}:}
We prove this proposition using three lemmas and their corollaries.

\begin{lemma}
\label{lemmaC2}
If $a\in L_0$ such that $0<a-c_1\le (11/10)(v^f-c_1)$, then
$K_\ell(a)<2$, for all $l$ big enough.
\end{lemma}
\pr
If $y\le (3/4)\cdot (11/10)$, then, by Corollary~\ref{cor51},
$K^\ast_\ell\sim 1/(e \cdot \log(1/y))\le 1.9123...<2$.
\qed

Let $P_\ell(z)=z^\ell$.
The next lemma gives information about the asymptotic shape
of $P_\ell(D((-1,1);\theta))$ as $\ell\to \infty$. We need even something
more general. Fix $K$ between $+\infty$ and $-1$, some $C>0$, and
$\theta\in (0,\pi/2)$, and consider a Poincar\'e neighbourhood
$D=D((-(1+{C\over \ell})K,1+{C\over \ell});\theta)$ of the interval
$(-(1+{C\over \ell})K,1+{C\over \ell})$
(here $K\ge -1$).
Take a real $\Lambda>0$ and consider a point $z_\ell(\Lambda)$ of
the boundary
of the above $D$ with $\arg(z)=\Lambda/\ell$.

\kies{
\begin{figure}[htp]
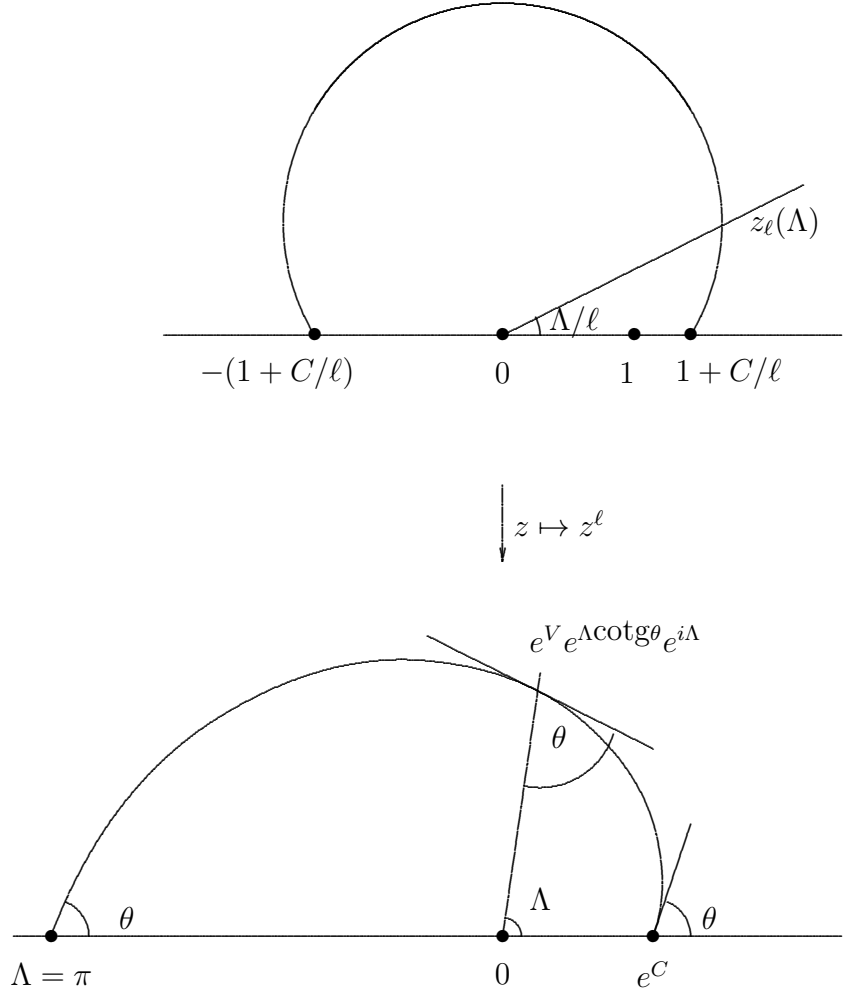
 \hfil
\beginpicture
\dimen0=0.5cm
\setcoordinatesystem units <\dimen0,\dimen0>
\setplotarea x from -9 to 9, y from -2 to 9
\setlinear
\plot -9 0 9 0 /
\plot 0 0 8 4 /
\circulararc 180 degrees from 5 0  center at  0 3
\circulararc -180 degrees from -5 0  center at  0 3
\multiput {\smalll $\bullet$} at 5 0  3.5 0  0 0  -5 0 /
\put {\smalll $0$} at 0 -1
\put {\smalll $1$}   at 3.3 -1
\put {\smalll $1+C/\ell$} at 6 -1
\put {\smalll $-(1+C/\ell)$} at -6 -1
\put {\smalll $z_\ell(\Lambda)$} at 7.5 3
\put {\smalll $\Lambda/\ell$} at 1.9 0.35
\circulararc 30 degrees from 1 0  center at  0 0
\arrow <5pt> [0.2,0.4] from 0 -4 to 0 -6
\put {\smalll $z\mapsto z^\ell$} at 1.5 -5
\setcoordinatesystem units <\dimen0,\dimen0> point at 0 16
\setplotarea x from -13 to 9, y from -2 to 7
\setlinear
\plot -13 0 9 0 /
\setquadratic
\plot
   4. 0
4.246574129 1.379795577
4.032389007 2.929702101
3.270362207 4.501267417
1.919251673 5.906849274
0           6.933012072
-2.391534848 7.360387428
-5.077915007 6.989150415
-7.801835337 5.668365165
-10.23806432 3.326548752
-12.01666410 0 /
\setlinear
\plot 0 0 1 7 /
\plot -2 8 4 5 /
\plot 4 0 5 3 /
\multiput {\smalll $\bullet$} at 4 0  0 0  -12 0  /
\put {\smalll $e^Ve^{\Lambda\text{cotg}\theta}e^{i\Lambda}$} at 3 8
\put {\smalll $\theta$} at 1.5 5.3
\put {\smalll $\theta$} at 5.5 0.5
\put {\smalll $\theta$} at -10 0.5
\put {\smalll $\Lambda=\pi$} at -12 -1
\put {\smalll $\Lambda$} at 1 1
\put {\smalll $0$} at 0 -1
\put {\smalll $e^C$} at 4 -1
\circulararc 80 degrees from 0.5 0  center at  0 0
\circulararc 68 degrees from 5 0  center at  4 0
\circulararc 84 degrees from 0.6 4  center at  1 6
\circulararc 68 degrees from -11 0  center at  -12 0
\endpicture
\caption[ ]{{\smalll The image of a Poinacar\'e disc under the map
$z\mapsto z^\ell$.}}
\end{figure}


}

\begin{lemma}
\label{lemmaC3}
$P_\ell(z_\ell(\Lambda))$ tends, as $\ell\to\infty$, to
a point 
$$\exp(C)\exp\left(\Lambda{\cos(\theta)\over \sin(\theta)}\right)
\exp(i\Lambda)$$
of a logarithmic
spiral, and the convergence is uniform in $\Lambda$ on every compact
of $(0,+\infty)$.
\end{lemma}
\pr
 If $z\in \partial D$, then, by \ref{B1},
$$z=(1+{C\over \ell})\cdot \left[1+{K+1\over 2}\cdot {i\cdot \exp(i\theta)\over
\sin(\theta)}\cdot (1-\exp(i\phi))\right],$$
where $\phi$ is the angle between the vectors
$z-z_0,1-z_0$, with $z_0$ the centre of the circle $D$.
As $\arg(z)=\Lambda/\ell\to 0$ with $\ell\to \infty$, then $\phi\to 0$
uniformly in $\Lambda$ on every compact
of $(0,+\infty)$ (remember that $\theta$ is fixed). So,
$$z^\ell\sim \exp(C)\exp\left\{{K+1\over
2\sin(\theta)}\exp(i\theta)\ell\phi\right\},$$
as $\ell\to \infty$. Let us prove that $\ell\phi\to \Lambda{2\over K+1}$.
Indeed, by the expression for $z$ and using the notation
$\alpha=\Lambda/\ell$,
$$\tan(\alpha)={{K+1\over
2\sin(\theta)}(\cos(\theta)-\cos(\theta+\phi))\over
1+{K+1\over
2\sin(\theta)}(\sin(\theta+\phi)-\sin(\theta))}.$$
Then
$$\cos(\theta+\phi-\alpha)=\cos(\theta)-2\sin(\alpha/2)
(\cos(\theta)\sin(\alpha/2)-{K-1\over K+1}\sin(\theta)\cos(\alpha/2)).$$
It follows, as $\alpha\to 0$ (and $\theta=$const),
$$\theta+\phi-\alpha\sim \theta+{\alpha(-{K-1\over K+1}\sin(\theta))\over
\sin(\theta)},$$
i.e., $\phi\sim 2\alpha/(K+1)$.
The uniform convergence also follows, and the statement is proved.
\qed

Given $A\ge 1$, $\theta\in (0,\pi/2)$, and $0\le \Lambda_1<\Lambda_2\le \infty$, 
we denote
$$\Gamma(A,\theta;\Lambda_1,\Lambda_2)=\{z=A\exp\left\{\Lambda{\cos(\theta)\over
\sin(\theta)}\right\}\exp(i\Lambda)\st\Lambda_1<\Lambda<\Lambda_2\}$$
a part of the logarithmic spiral. We have proved in the lemma above
that an arc
of $\partial D((-1,1);\theta)$ of the points $z$ with $0<\arg(z)<\Lambda/\ell$,
where $0<\Lambda<\infty$,
is mapped by $P_\ell$ asymptotically onto $\Gamma(1,\theta;0,\Lambda)$, and an
arc of $\partial D((0,1+(1/\ell)\log(11/10));\theta)$ of the points with
$0<\arg(z)<\Lambda/\ell$
is mapped by $P_\ell$ asymptotically onto
$\Gamma(11/10,\theta;0,\Lambda)$.
Let us note that the arc $\Gamma(1,\theta;0,\Lambda)$ is inside
$D((-K,1);\theta)$, $K>1$, if $\Lambda>0$ is small enough, because these curves
are
tangent at $1$, but the curvature of the logarithmic spiral at $1$ is less
than the curvature of the curve (a part of a circle) 
$\partial D((-K,1);\theta)$ at $1$. On the other
hand, it is clear that for given $K$ and for $A$ big enough the spiral
$\Gamma(A,\theta;0,\infty)$ is already outside of $D((-K,1);\theta)$.
We are going to find a lower bound for $A$.

Fix $K>1$, and $A>1$.

\kies{

\begin{figure}[htp]
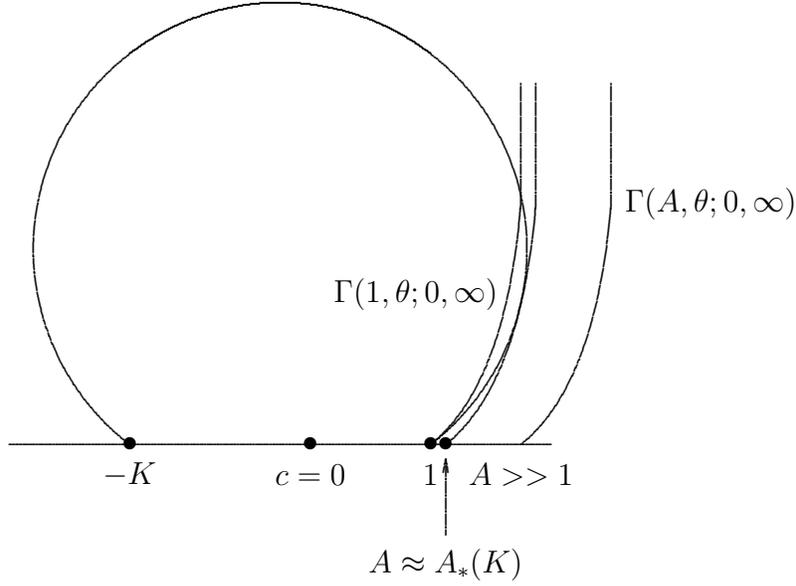
 \hfil
\beginpicture
\dimen0=0.4cm
\setcoordinatesystem units <\dimen0,\dimen0>
\setplotarea x from -9 to 9, y from -5 to 15
\setlinear
\plot -9 0 9 0 /
\circulararc 180 degrees from 5 0    center at  0 6.5
\circulararc -180 degrees from -5 0  center at  0 6.5
\multiput {\smalll $\bullet$} at 1 0  5 0  5.5 0  -5 0 /
\put {\smalll $c=0$} at 1 -1
\put {\smalll $-K$}   at -5 -1
\put {\smalll $1$}   at 5 -1
\setquadratic
\plot 5   0   7   3   8     8 /
\plot 5.5 0   7.5 3   8.5   8 /
\plot 8   0   10   3   11     8 /
\setlinear
\plot 8 8 8 12 /
\plot 8.5 8 8.5 12 /
\plot 11 8 11 12 /
\put {\smalll $\Gamma(1,\theta;0,\infty)$} at 4.5 5
\put {\smalll $\Gamma(A,\theta;0,\infty)$} at 14.3 8
\put {\smalll $A>>1$} at 8 -1
\arrow <5pt> [0.2,0.4] from 5.5 -3 to 5.5 -0.5
\put {\smalll $A\approx A_*(K)$} at 5.5 -4
\endpicture
\caption[ ]{{\smalll $D((-K,1);\theta)$ when $\theta<<1$}}
\end{figure}


}

\begin{lemma}
\label{lemmaC4}
If for some sequence of $\theta$ tending to zero,
the (open) curve $\Gamma(A,\theta;0,\pi)$ intersects the curve
$\partial D((-K,1);\theta)$ at a point $Z(\theta)$, then $Z(\theta)$
tends to $1+x$, where $x$ is a positive (real) solution of the equation
\beq
A\exp\left\{{x(1+{x\over K+1})\over 1+x}\right\}-(1+x)=0.
\label{C1}
\eeq
\end{lemma}

\begin{remark}
It can be seen from the proof below that this condition
is also `only if'.
\end{remark}

\noindent
\pr We have at the point of intersection $Z(\theta)$
(of argument $\alpha$):
\beq
A\exp\left\{\alpha{\cos(\theta)\over \sin(\theta)}\right\}\exp(i\alpha)=
1+{K+1\over 2}\cdot {i.\exp(i\theta)\over
\sin(\theta)}\cdot (1-\exp(i\phi)),
\label{C2}
\eeq
and a consequence is the equality for arguments:
\beq
\tan(\alpha)={{K+1\over
\sin(\theta)}\sin({\phi\over 2})\sin(\theta+{\phi\over 2})\over
1+{K+1\over
\sin(\theta)}\sin({\phi\over 2})\cos(\theta+{\phi\over 2})}.
\label{C3}
\eeq
Remember that some sequence of $\theta\to 0$. A priori the following cases are
possible for some subsequence:

\noindent
{\bf I.} $\phi/\theta\to \infty$. We are going to prove that this case
is, in fact, impossible.

\noindent
{\bf II.} $\phi/\theta\to t<\infty$.

\noindent
{\bf Case I} is divided into three subcases.

Ia.\,\, $\phi\to \pi$. Then (\ref{C3}) and $\phi/\theta\to \infty$
gives
$\tan(\alpha)\to \infty$, i.e., $\alpha\to \pi/2$. Now we compare
the left-hand side (LHS) and the
right-hand side (RHS) of (\ref{C2}).
The modulus of the LHS is equal to
$$A\exp\{\alpha{\cos(\theta)\over \sin(\theta)}\}\sim
A\exp({\pi\over 2\theta})$$
while for the modulus of the RHS we can write
$$|1+{K+1\over 2}\cdot {i\cdot \exp(i\theta)\over
\sin(\theta)}\cdot (1-\exp(i\phi))|\sim |1+i(K+1){1\over \theta}|,$$
so the equality (\ref{C2}) cannot hold in this case.

Ib. \,\, $\phi\to 2\pi$. Then $\alpha>\pi/2$, i.e.,
$$A\exp\{\alpha{\cos(\theta)\over \sin(\theta)}\}>
A\exp({\pi\cos(\theta)\over 2\sin(\theta)})\sim
A\exp({\pi\over 2\theta})$$  
while the modulus of the RHS of (\ref{C2}) can be at most
$$1+{K+1\over \theta},$$
as $\theta\to 0$.
This is again a contradiction.

Ic.\,\, $\phi$ tends neither to $\pi/2$ nor to $\pi$
(but $\phi/\theta\to \infty$). Then, from (\ref{C3}),
$\alpha\sim \theta+\phi/2$, and the modulus of the
LHS of (\ref{C2}) is at least
$$A\exp{3\over 4}(1+{\phi\over 2\theta})$$
while the modulus of the RHS of (\ref{C2}) is less than
$$1+(K+1){\phi\over 2\theta}.$$
Since $\phi/\theta\to \infty$, this is impossible again.

\noindent
{\bf Case II:} $\phi/\theta\to t<\infty$ (as $\theta\to 0$
along a sequence). Then 
$$(K+1){\sin(\phi/2)\over \sin(\theta)}\to x=(K+1)t/2<\infty$$
and $\tan(\alpha)\to 0$. If $t=0$, then the RHS tends to $1$,
but $|LHS|>A>1$. Thus, $t$ and $x$ are not zero, and, from (\ref{C3}),
$$\alpha\sim {x(1+\phi/(2\theta))\over 1+x}\theta
\sim {x(1+{x\over K+1})\over 1+x}\theta.$$
Then the RHS of (\ref{C2}) tends to
$$1+x.$$
Substituting these in (\ref{C2}), we obtain the equation (\ref{C1})
for $x$. Moreover, the point of intersection tends
to $1+x\in \rz$.
The lemma is proved.
\qed

\begin{corr}
\label{corC1}
Given $K>1$, if
\beq
A>A_\ast(K)={K\over \exp(2{K-1\over K+1})},
\label{C4}
\eeq
then, for all $\theta$ close
enough to zero, the arc $\Gamma(A,\theta;0,\pi)$ of the
logarithmic spiral does not intersect the domain $D((-K,1);\theta)$.
\end{corr}

\begin{remark}
The condition (\ref{C4}) is also `only if'.
\end{remark}

\noindent
\pr It is enough to prove only that, with this particular choice
of $A$,  equation
(\ref{C1}) does not have positive solutions. If $A>1$ is close enough to $1$,
equation (\ref{C1}) has at least two positive solutions. On the other hand,
since the second derivative of the left-hand side of (\ref{C1}),
$$A\exp\left\{{x(1+{x\over K+1})\over 1+x}\right\}{1\over (K+1)^2}
\left[1+{2K\over (1+x)^2}-{2K(K+1)\over (1+x)^3}+{K^2\over (1+x)^4}\right],$$
has exactly one positive root, the number of positive roots of (\ref{C1})
is at most two. If $A$ is large, there are no positive roots at all.
Hence, there exists some $A_\ast>1$, such that with $A<A_\ast$ there are
two roots, and with $A>A_\ast$ there are no roots, and $A_\ast$ can be
defined by a condition that with $A=A_\ast$, the equation has one multiple
positive root $x$. So, 
$$A_\ast\exp\{{x(1+{x\over K+1})\over 1+x}\}=1+x,$$
and
$$A_\ast\left[{1\over K+1}+{K\over (K+1)(1+x)^2}\right]
\exp\left\{{x(1+{x\over K+1})\over 1+x}\right\}=1.$$
Then 
$${1\over K+1}+{K\over (K+1)(1+x)^2}={1\over 1+x},$$
i.e., $x$ is either $0$, or $K-1$. The zero corresponds to
the trivial value $A_\ast=1$. Substituting $x=K-1$, we 
obtain the formula (\ref{C4}).
\qed

\begin{corr}
\label{corC2}
Given an arbitrary interval $J\subset \rz$, for all
$\theta$ small enough, the curve $\Gamma(1,\theta;\pi,2\pi)$
is outside the closure of the domain $D(J;\theta)$.
\end{corr}
\pr Obviously, the curve $\Gamma(1,\theta;\pi,2\pi)$
is a curve $\Gamma(A_\theta,\theta;0,\pi)$ rotated by the angle
$\pi$, where $A_\theta=\exp(\pi{\cos(\theta)\over \sin(\theta)})$
tends to $\infty$ as $\theta$ tends to the zero.
\qed

\noindent
{\em Conclusion of the proof of the Proposition~\ref{propC1}.}
Remember that we chose the domain of renormalization
$$\Omega=D(V;\theta)\bigcup D(I;\theta)\bigcup D(-I;\theta),$$
where $V=(-v,v)$ and $v$ is the boundary point of $V$
so that $F(v)=u^f\in (c_1,v^f)$
(we assume that $v>0$), $\theta$ is some
absolute constant (angle) to be determined later on,
and $I=(0,v+{\log(11/10)\over \ell}v)$.

We need to find $\theta=\theta_0>0$ and $\ell_0$, such that for every
$\ell>\ell_0$,
$$f^{-1}\circ F(\partial \Omega)\subset \Omega.$$
To do this, consider the pullback $F(\Omega)$. Let us rescale $\Omega$
such that the interval $V=(-v,v)$ turns into the interval $(-1,1)$.
We call the obtained domains by $\Omega^{\ast}$.
Let us rescale also $F(\Omega)$ by shifting first by $-c_1$,
and then rescaling it so that 
the interval $(c_1,v^f)$ turns into the interval $(0,1)$.
We call the obtained domains by 
$\widetilde \Omega^{\ast}$. It is convenient to introduce
also the scaled map $F^\ast$ corresponding to the map $F$
(the pullback of $f^{s-1}$):
$$F^\ast(z)=\frac{F(vz)-c_1}{v^f-c_1}.$$
So $\widetilde \Omega^{\ast}=F^\ast(\Omega^{\ast})$.
It is enough to find
$\theta=\theta_0$ and $\ell_0$, such that for every $\ell>\ell_0$ we have
that the closure of $P_\ell^{-1}(\widetilde \Omega^{\ast})$ is inside
$\Omega^{\ast}$. (As above, $P_\ell(z)=z^\ell$.)
By our choice,
$$\Omega^{\ast}=D((-1,1);\theta)\cup D(I^\ast;\theta)\cup D(-I^\ast;\theta),$$
where $I^\ast=(0,1+{\log(11/10)\over \ell})$.

Let us look at the all parts of $P_\ell^{-1}\circ F^\ast$.
Given $z\not= 0$, we let $E_i(z)$ be a unique point $w$ such that
$P_\ell(w)=w^\ell=z$ and
$\arg(w)\in [(2i-1)\pi/\ell,(2i+1)\pi/\ell)$, $i=0,1,\dots,l-1$.
Because of Lemma~\ref{lemmaC1}, the restriction of $F^\ast$ to the real axis
is defined on the interval $(-({4\over 3})^{1/\ell},({4\over 3})^{1/\ell})$.
Moreover, the real map $E_0\circ F^\ast \colon
(0,({4\over 3})^{1/\ell})\to  \rz$
is a homeomorphism and it sends the interval
$(0,({4\over 3})^{1/\ell})$ into itself (since $R_V$
has a high return and $f$ has no attracting periodic orbit).
It follows that
\beq
E_0\circ F^\ast
(D(I^\ast;\theta))\text{ is a proper subset of }D(I^\ast;\theta),
\label{C5}
\eeq
for any angle $\theta\in (0,\pi/2]$.

Let us consider the rest of $E_0\circ F^\ast(\Omega^\ast)$,
i.e., the set $E_0(W)$, where
$$W=F^\ast\left(D((-1,1);\theta)\cup D(-I^\ast;\theta)\right).$$
The trace of $W$ on
the real axis is contained in the
interval $(-K_1,1)$, where for $K_1>0$ we have
a bound controlled by Lemma~\ref{lemmaC2}:
$${K_1\over 11/10}<2,$$
for all big $\ell$.
It follows, that the set $W$ is covered by $D((-K_1,1);\theta)$,
for any $\theta$.

Observe that by (\ref{C4}),
\beq
A_\ast(K_1)<A_\ast(2.2)=1.04...<11/10.
\label{C6}
\eeq
Applying Lemma~\ref{lemmaC3} and
Corollary~\ref{corC1}, we find an angle $\theta_1>0$
and a degree $\ell_1$, such that, for all $\theta\le \theta_1$,
and for all $\ell\ge \ell_1$, the set
$P_\ell(\Omega^\ast\cap \{z\st\arg(z)\in [-\pi/\ell,\pi/\ell]\})$ contains $W$,
that is $E_0(W)$ is inside $\Omega^\ast$. Therefore,
we have proved that
$E_0\circ F^\ast(\Omega^\ast)$ is inside $\Omega^\ast$.
By the symmetry of $\Omega^\ast$, 
$E_{\ell/2}\circ F^\ast(\Omega^\ast)$ is inside
$\Omega^\ast$ too, for the same $\theta$ and $\ell$.

Now we consider $E_1\circ F^\ast(\Omega^\ast)$. The domain $F^\ast(\Omega^\ast)$
is contained in a domain $D((-K_1,4/3);\theta)$. So, we can apply
Corollary~\ref{corC2} (together with Lemma~\ref{lemmaC2}) to conclude that,
for some $\theta_2>0$ and $\ell_2$, if
$\theta\le \theta_2$ and $\ell\ge \ell_2$, then 
$E_1\circ F^\ast(\Omega^\ast)$ is
contained in the Poincar\'e neighbourhood $D((-1,1);\theta)$.
Essentially, this is the end of the proof. Indeed, each other
$E_i\circ F^\ast(D((-1,1);\theta))$, ($i\not=0,\ell/2$), is contained in
$D((-1,1);\theta)$, since $D((-1,1);\theta)$
is invariant under the rotation $z\mapsto \exp(i\cdot 2\pi/\ell)z$,
for $z$ in the first quarter.

Thus, for $\theta=\theta_0=\min\{\theta_1,\theta_2\}$, and for every
$\ell>\ell_0=\max\{\ell_1,\ell_2\}$, we
have that $\widetilde \Omega^\ast$ is inside
of $\Omega^\ast$.
\qed

Thus we have completed the proof of Proposition~\ref{propC1}
and of Theorem B for the case when $\ell$ is large.

\sect{The proof of Theorem B for all degrees}
\label{secla}

In this section we complete the proof of Theorem B and
of the Main theorem in some non-renormalizable cases:

\begin{theo}
Let $f(z)=z^\ell + c_1$ with $\ell$ an even integer and $c_1$ real be a
non-renormalizable polynomial such that the limit set $\omega(c)\ni c$
is minimal and $f$ has
infinitely many high returns in the partition
given by the Yoccoz puzzle on the real line.
Then the Julia set
of $f$ is locally connected.
\end{theo}

Of course, we may assume in the proof
below that $\ell\ge 4$ because when $\ell=2$ then the result
holds (even without the assumption about high returns)
by the result of Yoccoz \cite{Y}, see the Introduction.
So let us fix $\ell\ge 4$. In the previous section we have
proved the above result already for $\ell$ sufficiently large.
Since the estimates in this section for $\ell$ `small'
are more delicate and since
the proof in the previous section shows that the shape of the
domain (i.e., $\theta$) can be chosen uniformly in $\ell$,
we have dealt with the asymptotic case
separately in the previous section.
We should note that the proof of this theorem
also holds for every infinitely renormalizable $f$
with $\ell\ge 4$. Since we use only the `easy space' $1/3$,
one might hope to extend this result to certain non-real polynomials.

Given the interval $V=(-v,v)$ such that $R_V$ has a high return
we construct an
$\ell$-polynomial-like map sitting
inside the domain $\Omega=\Omega(\ell,V)$, where $\Omega$ is either
the disc $D_*(V)$ based on the diameter $V$, or
$$\Omega=D(V;\theta)\bigcup D(I;\theta)\bigcup D(-I;\theta),$$
where
$\theta=\theta_0$ is some absolute constant (angle) to be determined later on,
and 
$$I=(0,1.07^{1/\ell}v).$$
Here $F(v)=u^f\in (c_1,v^f)$ and we may assume that $v>0$.

As before, it is enough to prove

\begin{prop}
\label{propD1}
Given $\ell\ge 4$, there exists $\theta_0>0$, such that 
$$f^{-1}\circ F(\partial \Omega)\subset \Omega.$$
\end{prop}

The proof of the Proposition~\ref{propD1} is
somewhat similar to the proof of the Main
Theorem in the infinitely renormalizable
case for  degree $2$ and the proof of the Theorem B for
sufficiently large degrees.
The main new ingredient is contained in the following lemma:

\begin{lemma}
\label{lemmaD1}
Either the disc $D_*(V)$ is a domain of the $\ell$-polynomial-like
mapping, i.e., $f^{-1}\circ F(D_*(V))\subset D_*(V)$, or otherwise
$F(D(-I;\theta))$ lies inside $D(I^f;\theta)$, where
$I^f$ is an interval around $c_1$:
$$I^f=(c_1-2.12|v^f-c_1|, c_1+0.68|v^f-c_1|).$$
\end{lemma}
\pr
Remember that the constant $K_{\ell}(v^f)=|\hat u^f-c_1|/|v^f-c_1|$ depends
not only on the extendability space (which is $1/3$), but on the parameter
$t=|c-c_{s+1}|/|T|$ as well (see Lemma~\ref{str1}). If $t\ge 0.51$, then
$$K_{\ell}(v^f)\le {0.51((3/4)^{1/\ell}-0.51^{1/\ell})\over
0.51^{1/\ell}(3/4)(1-(3/4)^{1/\ell})}\le {0.51((3/4)^{1/4}-0.51^{1/4})\over
0.51^{1/4}(3/4)(1-(3/4)^{1/4})}\le 0.991818...<1,$$
so that we apply Proposition~\ref{pqp}. Thus,
we can assume $t<0.51$. The right end
of the interval $F(-I)$ is just the point $c_{-s+1}$,
which belongs to the interval $(c_1,c_{s+1})$ (since we have
a high return). Hence,
$${|c_1-c_{-s+1}|\over |c_1-v^f|}\le {|c_1-c_{s+1}|\over |c_1-v^f|}\le
0.51\times (4/3)=0.68.$$
The left end  $b$ of the interval $F(-I)$ is obtained from Corollary 5.1,
where we put $y=1.07/(4/3)=0.8025$, so that $K^*_{\ell}(y)\le
K^*_4(y)=1.97063...<1.98$. Since $a=c_1+1.07|v^f-c_1|$ and
$K_{\ell}(a)=|b-c_1|/|a-c_1|<1.971$, indeed, $|b-c_1|/|v^f-c_1|<1.98\times
1.07<2.12$. 
\qed

If the first alternative in the lemma holds then
Proposition~\ref{propD1} holds. So we will assume in the
remainder of the proof that the second alternative holds.
The next lemma will allow us to apply Lemma~\ref{lem71} for any
degree $\ell\ge 4$ (and not just for $\ell=2$).
Let us denote for simplicity $D(\theta)=D((-1,1);\theta)$. Set
$$\Pi_\ell(\theta)=\{z\in D(\theta): 0\le \arg z\le \pi/\ell\}.$$
As before, $P_\ell(z)=z^\ell$.

\begin{lemma}
\label{lemmaD2}
Fix $0<\theta<\pi/2$. Then
$$P_\ell(\Pi_\ell(\theta))\subset P_{\ell+2}(\Pi_{\ell+2}(\theta)),\quad
\ell=2,4,... .$$ 
\end{lemma}
\pr
Assume the contrary. Then the boundaries of
$P_\ell(\Pi_\ell(\theta))$
and $P_{\ell+2}(\Pi_{\ell+2}(\theta))$ have a common non-real point, i.e.,
$z^\ell=u^{\ell+2}$, for some $z\in \partial D(\theta)$,
$0\le \arg z\le \pi/\ell$, and $u\in \partial D(\theta)$,
$0\le \arg z\le \pi/(\ell+2)$. Hence, $u=z^t$, with $t=\ell/(\ell+2)$ between
zero and $1$. The point $z^t$ belongs to an arc $\Gamma(1,\theta_1;0,\Lambda_0)$
of a logarithmic spiral starting at the points $1$ and  ending at $z\in
\partial D(\theta)$ and crossing the circle $\partial D(\theta)$ at 
the other point $u$. If $\theta_1\le \theta$, it is clearly impossible
(see Section~\ref{sec9}). Consider the case $\theta_1>\theta$. Then
$\Gamma(1,\theta_1;0,\Lambda)$ is inside  $D(\theta)$, if $\Lambda$ is small.
Hence, there are two points $z_1,z_2$ of the intersection of the arc
$\Gamma(1,\theta_1;0,\Lambda_0)$ with $\partial D(\theta)$, such that
$\arg z_1<\arg z_2$, and
this arc leaves the disc $D(\theta)$ at $z_1$ and again enters it at $z_2$.
By the geometry of the logarithmic spiral, the angle between the vector
$z_1$ and the circle $D(\theta)$ is at least $\theta_1$,
and the angle
between the vector $z_2$ and the circle $D(\theta)$ is at most
$\theta_1$. This is a contradiction with the fact that the angle
between a vector $w\in \partial D(\theta)$ and the tangent
to $\partial D(\theta)$ at $w$ is
increasing as $w\in \partial D(\theta)$ moves
from $1$ to $-1$ (in fact, it increases from $\theta$ to $2\pi -\theta$).
\qed
 
In order to prove Proposition~\ref{propD1},
we need to find for any $\ell\ge 4$
some $\theta=\theta_0>0$ such that
$f^{-1}\circ F(\partial \Omega)\subset \Omega$.
To do this, consider the pullback $F(\Omega)$. Let us rescale $\Omega$
such that the interval $V=(-v,v)$ turns into the interval $(-1,1)$.
We call the obtained domain $\Omega^{\ast}$.
Let us rescale also $F(\Omega)$ by shifting first by $-c_1$, 
and then rescaling it so that 
the interval $(c_1,v^f)$ turns into the interval $(0,1)$.
We call the obtained domain
$\widetilde \Omega^{\ast}$. It is convenient to introduce
also the scaled map $F^\ast$ corresponding to the map $F$
(the pullback of $f^{s-1}$):
$$F^\ast(z)=\frac{F(vz)-c_1}{v^f-c_1}.$$
So $\widetilde \Omega^{\ast}=F^\ast(\Omega^{\ast})$.
It is enough to find $\theta=\theta_0$ such that
the closure of $P_\ell^{-1}(\widetilde \Omega^{\ast})$ is inside
$\Omega^{\ast}$. (As above, $P_\ell(z)=z^\ell$.)
By our choice,
$$\Omega^{\ast}=D((-1,1);\theta)\cup D(I^\ast;\theta)\cup D(-I^\ast;\theta),$$
where $I^\ast=(0, 1.07^{1/\ell})$.
Let us look at each piece of $P_\ell^{-1}\circ F$.
Given $z\not= 0$, define $E_i(z)$ to be the unique point $w$ such that
$P_\ell(w)=w^\ell=z$ and
$\arg(w)\in [(2i-1)\pi/\ell,(2i+1)\pi/\ell)$, $i=0,1,\dots,l-1$.
Because of Lemma~\ref{lemmaC1}, 
the restriction of $F^\ast$ to the real axis
is defined on the interval $(-({4\over 3})^{1/\ell},({4\over 3})^{1/\ell})$.
Moreover, the real map $E_0\circ F^\ast \colon
(0,({4\over 3})^{1/\ell})\to  \rz$
is a homeomorphism and it sends the interval
$(0,({4\over 3})^{1/\ell})$ into itself (since $R_V$
has a high return and $f$ has no attracting periodic orbit).
It follows that
\beq
E_0\circ F^\ast(D(I^\ast;\theta))\text{ is a proper subset of }D(I^\ast;\theta),
\label{C5n}
 \eeq
for any angle $\theta\in (0,\pi/2]$.

Let us now show that
\beq
E_0(F^\ast(D((-1,1);\theta)))\subset \Omega,
\label{C5n1}
\eeq
in other words, that $F^\ast((D((-1,1);\theta)))$ is covered by
the set $P_\ell(\Omega\cap \{z\st -\pi/\ell<\arg z<\pi/\ell\})$.
To see this, first note that
by Corollary~\ref{cor51}, $F^\ast((D((-1,1);\theta)))\subset D((-K_0,1);\theta)$
with the constant 
$$K_0=K^*_\ell(4/3)\le K^*_4(4/3)=1.51983...<1.52.$$
By Lemma~\ref{lemmaD2}, the difference
$\Delta_\ell(\theta)=D((-K_0,1);\theta)\setminus P_\ell(\Pi_\ell(\theta))$
is contained in $\Delta_2(\theta)$.
Hence (\ref{C5n1}) follows from:

\begin{lemma}
\label{lemmaD3}
For all $\theta$ small enough,
$$\Delta_2(\theta)\subset P_\ell\left( D(I^\ast;\theta)\cap
\{z\st 0<\arg z<\pi/\ell\}\right), \ell=4,6,\dots \quad .$$
\end{lemma}
\pr
Assume that this is not the case for some sequence $\theta\to 0$.
Then $z_1=z_2^\ell$, for some $z_1\in
\partial D((-K_0,1);\theta)$ and $z_2\in \partial D(I^\ast;\theta)$.
Moreover, $\arg z_2<\pi/\ell$, and, what is crucial, since
we were able to apply Lemma~\ref{lem71},
$z_1$ tends to a point of the real interval $[1,K_0^2]$
(for some subsequence) as $\theta\to 0$.
We have 
\beq
z_1=1+{K_0+1\over 2}\cdot {i\cdot \exp(i\theta)\over \sin(\theta)}\cdot
(1-\exp(i\phi)).
\label{D10}
\eeq
where $\phi\in (0,2\pi-\theta)$ is the angle between the vectors
$z_1-C,1-C$, with $C$ the centre of the circle $D((-K_0,1);\theta)$.
For $z_2\in \partial D(I^\ast;\theta)$ we have a similar expression:
\beq
z_2=A\cdot {\sin(\theta+\gamma)\over \sin\theta}\cdot\exp(i\gamma),
\label{D11}
\eeq
where $A=1.07^{1/\ell}$ and $\gamma\in (0,\pi/\ell)$ is an argument of $z_2$.
Since $z_1$ tends to a real point in $[1,K_0^2]$, it follows from (\ref{D10}),
that $\phi/\theta$ tends to a non-negative finite constant $B$, as $\theta\to
0$, and $1+{K_0+1\over 2}B\le K_0^2$, i.e.,
\beq
0\le B\le 2(K_0-1).
\label{D12}
\eeq
Hence, from the condition $z_1=z_2^\ell$ and from (\ref{D11}),
$\gamma/\theta$ tends to a finite $D\ge 0$. 
Separating now real and imaginary
parts of the equality $z_1=z_2^\ell$, we obtain
the following system for $B$ and $D$:
\beq
1+{K_0+1\over 2}\cdot B=1.07\cdot(1+D)^\ell
\label{D13}
\eeq
\beq
{K_0+1\over 2}\cdot B\cdot (1+{B\over 2})=1.07(1+D)^\ell\cdot \ell\cdot D,
\label{D14}
\eeq
where
$$K_0=1.52.$$
Dividing (\ref{D14}) by (\ref{D13}) and substituting the obtained expression
for $D$ into (\ref{D13}), we come to the equation:
\beq
1+{K_0+1\over 2}\cdot B=A\cdot \{1+{1\over \ell}\cdot {K_0+1\over 2}\cdot
{B\cdot (1+{B\over 2})\over 1+{K_0+1\over 2}\cdot B}\}^\ell,
\label{D15}
\eeq
where
$$K_0=1.52$$
and
$$A=1.07.$$
With these $K_0$ and $A$, this equation (\ref{D15})
has no solutions on the interval
$[0,2(K_0-1)]$ for $\ell=4$,
and, hence, for all $\ell\ge 4$.
In order to see this
we claim that, given $A>1$, this equation has either exactly
two non-negative solutions (maybe one multiple),
or no non-negative solutions at all. Before proving this claim
let us first show that this implies the lemma.

Indeed, if $B=2(K_0-1)$ is a solution,
for some parameter $A_0$, then
$$A_0={K_0^2\over (1+{1\over 4}\cdot {K_0^2-1\over K_0})^4}=1.05835... .$$
On the other hand,
taking the derivative of both sides of (\ref{D15})
(with $\ell=4$) with respect to $B$, we obtain, of course, $D_1=(K_0+1)/2$
on the left hand-side, and
$D_2=(K_0+1)/2\cdot 4(K_0^2-K_0+1)/(K_0^2+4K_0-1)$ on the right-hand side.
Since $D_1>D_2$, it means that $2(K_0-1)$ is the smallest positive solution of
(\ref{D15}) (for $A_0$). Since $A=1.07>A_0=1.05835$,
the smallest positive solution of (\ref{D15}) for
$A=1.07$ is therefore at least $2(K_0-1)$.

So it remains to prove the above
claim. For this it is enough to show that the second derivative
of the right-hand side of (\ref{D15})
(with $\ell=4$) has exactly one positive root.
Let us make a linear change of the variable: define
$x=1+{K_0+1\over 2}B$, so that $1\le x< \infty$. Then
$$\frac{K_0+1}{2} \frac{B(1+\frac{B}{2})}{1+{K_0+1\over 2}B}=1+{x-2\over K_0+1}
-{K_0\over K_0+1}{1\over x}:=T(x)$$
(the latter equality is just notation).
Hence
$$T'(x)={1\over K_0+1}(1+K_0{1\over x^2}),$$
$$T''(x)=-{2K\over K_0+1}{1\over x^3}.$$
And the second derivative of the right-hand side of (\ref{D15})
w.r.t. $x$ is (after calculations):
$$(1+T(x)/4)^2{1\over 4(K_0+1)^2}{1\over x^4}\times
\left\{3x^4+4K_0x^2-2K_0(5K_0+3)x+5K_0^2\right\}.$$
The polynomial in $\{...\}$ has no more than two positive roots
(because the derivative of $\{...\}$ w.r.t. $x$ is an increasing function
of $x\ge 0$). By checking the values of $\{...\}$ at $x=0,1,\infty$
it follows that it does have one positive root between $0$ and $1$ and
at least one root $>1$. So, it has exactly one root greater than $1$
and the claim follows.
\qed

Because of (\ref{C5n}) and (\ref{C5n1}), in order
to conclude that $E_0\circ F^* (\Omega)\subset \Omega$,
we only have to show that
$E_0\circ F^\ast(D(-I^*;\theta))\subset \Omega^*$, for $\theta$ small.
Lemma~\ref{lemmaD1} says that $F^\ast(D(-I^*;\theta))\subset
D((-2.12,0.68);\theta)$. Therefore, by the remark below
Lemma~\ref{lemmaD2}, for this it is enough to check

\begin{lemma}
\label{lemmaD4}
If $\theta$ is small, then
\beq
D((-2.12,0.68);\theta)\subset P_2(D((-1,1);\theta)).
\label{D16}
\eeq
\end{lemma}
\pr
Since $1/0.68=1.47059>1.47$ and $2.12/0.68=3.11765<3.12$,
it is enough to prove (after rescaling) that
\beq
D((-3.12,1);\theta)\subset P_2(D((-(1.47)^{1/2},(1.47)^{1/2});\theta)).
\label{D17}
\eeq
For a possible point $Z$ of intersection of the boundaries, we obtain an
equation 
$$1+{3.12+1\over 2}\cdot {i\exp(i\theta)\over \sin\theta}\cdot (1-\exp(i\phi))=
1.47\cdot(1+{\exp(i(\theta+\alpha))\over \sin^2\theta}\cdot
(2\cos\theta-2\cos(\theta+\alpha))).$$
If $\theta\to 0$ and $\phi,\alpha$ is a solution of this equation, then
$\phi/\theta$ and $\alpha/\theta$ tend
to finite constants $M\ge 0$ and $N\ge 0$ respectively
(proof: $P_2(D((-(1.47)^{1/2},(1.47)^{1/2});\theta))$ certainly contains
$P_2(D((-1,1);\theta))$ while the boundary of the latter domain intersects
$D((-K,1);\theta)$ at a point $Z(K,\theta)$ of an angle $\phi_1$ such
that $\phi_1/\theta$ is bounded as $\theta\to 0$, see Lemma~\ref{lem71}.)
We have the following equations for $M$ and $N$:
\beq
1+{3.12+1\over 2}\cdot M=1.47\cdot (1+4\cdot {N\over 2}\cdot (1+{N\over 2}))
\label{D18}
\eeq
\beq
{3.12+1\over 2}\cdot M\cdot (1+{M\over 2})=1.47\cdot 4\cdot
{N\over 2}\cdot (1+{N\over 2})\cdot (1+N).
\label{D19}
\eeq
This system has no non-negative solutions $(M,N)$. A way to see this is
to reduce the system to a polynomial equation. For this, denote $M+1=x$,
$N+1=y$. Then
\beq
{x-1\over 2}={1.47y^2-1\over 3.12+1}
\label{D20}
\eeq
\beq
{x+1\over 2}={1.47(y^2-1)y\over 1.47y^2-1}.
\label{D21}
\eeq
This implies that $y$ is a zero of the polynomial
\beq
h(y):=y^4-{3.12+1\over 1.47}y^3+{3.12-1\over 1.47}y^2+{3.12+1\over 1.47}y-{3.12\over
1.47^2}.
\label{D22}
\eeq
However, the polynomial $h(y)$ does not have solutions $y\ge 1$.
Indeed, the second derivative of this polynomial $h''(y)$
is a parabola with zero's  at $y=.2000905878$ and at $y=1.201269956$.
If follows that $h'(y)$ is a cubic function with a
local maximum at  $y=.2000905878$ and a local minimum at $y=1.201269956$.
An explicit calculation shows that $h'(1.201269956)>0$
and therefore, one has that $h'(y)> 0$ for each $y\ge 1$.
Hence $h(y)\ge h(1)>0$ for each $y\ge 1$.
\qed

Thus, we have proved that
$E_0\circ F^\ast(\Omega^\ast)$ is inside $\Omega^\ast$.
By symmetry of $\Omega^\ast$, $E_{\ell/2}\circ F^\ast(\Omega^\ast)$ is
inside
$\Omega^\ast$ too, for the same $\theta$ and $\ell$.

Now we consider $E_1\circ F^\ast(\Omega^\ast)$. First note that
the domain $F(\Omega^\ast)$
is contained in $D((-2.12,4/3);\theta)$. So, for given $\ell\ge 4$,
and for $\theta$ sufficiently small,
the domain $P_\ell(D((-1,1);\theta)\cap \{z\st
\pi/\ell\le\arg z\le3\pi/\ell\})$ and its complex conjugate contain
$D((-2.12,4/3);\theta)$ since the diameter of the latter domain grows as
$const/\sin\theta$ while the diameter of the former domain grows as
$const/\sin^\ell\theta$ as $\theta\to 0$.

This is enough
to conclude the proof of the Proposition~\ref{propD1}.
Indeed, each other
$E_i\circ F^\ast(D((-1,1);\theta))$, ($i\not=0,\ell/2$), is contained in
$D((-1,1);\theta)$, since $D((-1,1);\theta)$
is invariant under the rotation $z\mapsto \exp(i\cdot 2\pi/\ell)z$,
for $z$ in the first quadrant.
\qed

As we noted above, Proposition~\ref{propD1} implies Theorem 12.1.

\sect{The proof of the Main Theorem and Theorem A
in the `period doubling case'}
\label{secll}

In this section we shall modify  the proof in
the previous section in order to show that the Main Theorem also
holds in the case 
of an infinitely renormalizable map
of period doubling type (from some renormalization onwards).
This case was not dealt with in Sections~\ref{sec6}
and \ref{sec7} because the space $0.6$ which is used there,
only holds in the case that $f$ is not renormalizable
of both period $s$ and period $s/2$.
In that exceptional case, the space is merely $0.5$, see
Lemma~\ref{str6}  and therefore we can
use the method of round discs as in Section~\ref{sec6} only for
$\ell\ge 8$, see Example 5.1. Therefore we shall use
the ideas of the previous in this case when $\ell<8$. 
These arguments also show that Theorem A holds in this exceptional
case (that $f$ is renormalizable of both periods $s(n)$ and $s(n)/2$).

So let us indicate the differences with the proof
in the previous section.
Of course, the proof of Theorem A already follows from the
previous section if $\ell\ge 4$ and so we have to
take $\ell=2$ in the previous section.
Firstly, define $\Omega$ as before 
with the difference that we take the interval $I=(0,1.09^{1/2})$ now.
Lemma~\ref{lemmaD1} and its proof
go through unchanged (replacing $1.07$ by $1.09$)
because the
actual constants for $\ell=2$ and space $0.5$ are even
better than as in the proof of this lemma.
Lemma~\ref{lemmaD2} is not needed.
In Lemma~\ref{lemmaD3} we have to take $\ell=2$ in the statement.
In the proof we take $K_0=1.4$ and $A=1.09$. The calculations
are slightly different but it is easy to check that
everything works as before. Finally, Lemma~\ref{lemmaD4}
and its proof go through unchanged.
All this concludes the Main Theorem in this case.
Theorem A for this case
follows also in the same way as in Section~\ref{sec7}.
\qed

 %

\sect{Proof of Theorem C}
The proof in this section is an elaboration of
Section 5 of \cite{L3} and Lemmas 14 and 15 in \cite{L5}.
We wish to thank Edson Vargas for several discussions on these
sections.

Let ${\cal E}(T^0)$ be the collection of mappings
$g\colon \cup T_i^1\to T^0$ where $T^0$
is some symmetric interval around $c$ with nice boundary points
and where
$T_i^1$ is a finite collection
of disjoint closed subintervals of $T^0$ for which
\begin{itemize}
\item for $i\ne 0$ the map $g\colon T_i^1\to T^0$
is a diffeomorphism of the form $f^{j(i)}$
and the inverse map $(g|T_i^1)^{-1}$
has a univalent extension to $\cz_{T^0}$;
\item writing $T^1=T_0^1$ we have that
$g|T^1$ is a unimodal map of the form $f^j$
and with $g(\partial T^1)\subset \partial T^0$;
one can write $g|T^1=h\circ f$ where $h^{-1}$ has also a
univalent analytic extension to $\cz_{T^0}$;
\item all iterates of the critical point
$c=0$ under $g$ are in $\cup T_i^1$.
\end{itemize}
Assume there exists a symmetric interval $T^{-1}$
containing $T^0$, so that when writing
as before $g|T^1_i=h_i\circ f$, the map $h_i^{-1}$ has a
univalent analytic extension from $\cz_{T^{-1}}$
into $\cz_{H_i}$ where $H_i$ is some interval containing $f(c)$
such that $f^{-1}(H_i)\cap \rz\subset T^0$.
If this holds then we say that $g\in {\cal E}(T^0,T^{-1})$.

An example of a map $g$ which is of type ${\cal E}(T^0,T^{-1})$
is the first return map
to an interval $W_n=[u_n,\tau(u_n)]$ as in Section 2. More precisely,
since we have assumed that $\omega(c)$ is minimal, we
only consider the finitely many branches which contain points
from $\omega(c)$.  The $f$-image of each branch
can be extended to $W_{n-1}$
(hence the first return map
is in ${\cal E}(W_n,W_{n-1})$).
Indeed, the boundary points of $W_n$ are nice and there are no forward
iterates of $\partial W_n$ in $W_{n-1}\setminus W_n$.
So take a domain $I\cap W_n=\emptyset$
of the first return map $R$ to $W_n$
and a maximal interval $T$ containing $I$
so that h $R|I=f^i$ is monotone.
By maximality of $T$ for each component $T_+$ of
$T\setminus I$ there exists some $j<i$ so that $f^j(T_+)$ contains $c$
in its boundary and since $R$ is the first return map
$f^j(I)\cap W_n=\emptyset$. Hence $f^j(T_+)$ contains
a boundary point of $W_n$ and therefore $f^i(T_+)$ contains
a point of $\partial W_{n-1}$. Since this holds for both components
of $T\setminus I$ this gives $f^i(T)\supset W_{n-1}$.
So on each branch $I\subset W_n$ of $R$ one can write
$R=h\circ f$ and $h$ extends as a diffeomorphism from
some neighbourhood of $f(I)\ni c_1$ onto $W_{n-1}$.

We say that $g\in {\cal E}(T^0)$ has a {\it low return iterate} if
$g(T^1)$ does not contain the critical point $c$.
In this case we define $\R g\in {\cal E}$ as follows.
First define $\R g$ so that
it coincides with $g$ on $\cup_{i\ne 0}T_i^1$.
Define $s_0\ge 2$ to be minimal so that $g^{s_0-1}(c)\in T^0\setminus T^1$
and let $s\ge s_0$ be minimal so that $g^s(c)\in T^1$.
(This means that $g^s(T^1)\cap T^1\ne \emptyset$.)
Because $g$ has a low return iterate and no periodic attractors, 
$s_0$ exists ($s_0=2$ if $g(T^1)\cap T^1=\emptyset$)
and since $\omega(c)$ is minimal the integer
$s$ also exists. Therefore we can define the new central
domain of $T^2$ to
be the component of the domain of $g^s$ containing
$c$. Note that by the choice of $s$ we have
$\R g(T^2)\cap T^1\ne \emptyset$.
Moreover, by the way $s_0$ is chosen we also have
$g(T^1)\cap T^2=\emptyset$.
For $x\in T^1\setminus T^2$ let $s(x)\le s$ be the smallest integer
for which $g^{s(x)}(x)$ and $g^{s(x)}(c)$ are in different
components of $\cup_{i\ne 0}T_i^1$
and define $\R g(x)=g^{s(x)+1}(x)$. The domains of $\R g$
in $T^1\setminus T^2$ map diffeomorphically onto $T^0$.
In fact, we even have $\R g\in {\cal E}$. To show this it suffices
to show that if $\R g|T^2=f^m$ then
there exists an interval $A\supset f(T^2)$ which
is mapped diffeomorphically onto $T^0$ by $f^{m-1}$.
Indeed, if $g|T^1=f^j$ then the Koebe space of $f^{j-1}|f(T^1)$
spreads over $T^0$ and in particular over
one of the intervals connecting $\partial T^0$ to $c$.
Applying $g^{s_0-2}|T^1=f^{(s_0-2)j}$ to this it follows
that the Koebe space of
$g^{s_0-2}\circ f^{j-1}|f(T^1)$ spreads over one of the intervals
connecting $\partial T^0$ to $g^{s_0-2}(c)$
(which by definition is in $T^1$)
and in particular this space contains $\cup_{i\ne 0} T^1_i$.
Since $g$ maps such intervals onto $T^0$ and $s>s_0$
is minimal, the Koebe space of $f^{m-1}|f(T^2)$ spreads also over $T^0$.
This we have proved that $\R g\in {\cal E}$.
If $\R g$ again has a low return then we can define
$\R^2g\in {\cal E}$ and so on.

If $g$ has a {\it high return iterate} (i.e., not a low return iterate)
then let $x$ be the orientation preserving fixed point of $g|T^1$
and $z$ the boundary of $T^0$ on the same side of $c$ as $x$.
Next take preimages of $z_0=z,z_1,z_2,\dots$ of $z$
along this branch of $g|T^1$
(so $[z_1,\tau(z_1)]=T^1$ and
$z_k\to x$ as $k\to \infty$).
Define $U_k=[z_k,\tau(z_k)]$  and choose $k$ minimal so that
$g U_k\supset U_k$. Such an integer $k$ exists
because we have assumed that $f$ is not renormalizable.
The interval $W^1:=U_k$
is the {\it escape interval} associated to a map $g$ with a high
return iterate. For $i=0,1,\dots,k-1$ define
the new map $\widetilde \W g$ on $U_i\setminus U_{i+1}$ as $g^{i+1}\circ g$
and define $\widetilde \W g$ on $U_k$ as the first
return map of $g$ to $U_k$. Note that $\widetilde \W g\notin \E(T^0)$
because the image of the central branch $T^2$ is contained in
the interior of $T^0$ (and does not stretch over to $\partial T^0$).
However, the first return map $\W g$ of $\widetilde \W g$ to $W^1$
is contained in $\E (W^1,T^0)$.
Note that
the domain $W^1$ of $\W g$ is smaller the domain $T^0$ of
$g$, and the extension associated to $\W g$
includes the domain of the original  map $g$.

\begin{lemma}\label{lyno}
Assume that $g\in \E(T^0)$
is a first return map to a symmetric interval
$T^0$ around $c$
with nice boundary bounds. Let $g,\R g,\dots,\R^{k-1}g$ have low returns
and let $T^1,\dots,T^{k+1}$ be the central intervals
corresponding to $g,\R g, \dots, \R^k g$.
Write  $\R^kg|T^{k+1}=f^m$ with $m\in \nz$ and
let $A\supset f(T^{k+1})\ni c_1$
be the interval which is mapped by $f^{m-1}$
diffeomorphically onto $T^0$.
Then
\begin{itemize}
\item there exists a sequence of integers
$k+1\ge p_1>p_2>\dots > p_r \ge 1$
and integers $n(p_1)<n(p_2)<\dots < n(p_r)$
such that $n=n(p_i)<m$ is the largest integer such that
$A^{p_i}:=f^{n(p_i)}(A)\subset T^{p_i-1}\setminus T^{p_i}$.
Then writing, $s(p_i):=n(p_{i+1})-n(p_i)>0$ for each $i=1,\dots,r-1$,
we have that
$f^{s(p_i)}$ maps $A^{p_i}$ onto $A^{p_{i+1}}$.
Moreover, the map $f^{s(p_i)-1}\colon fA^{p_i}\to A^{p_{i+1}}$ has
Koebe space spread over $T^{p_{i+1}-1}$;
\item $p_1=k+1$ or $k$ and
the Koebe space of $f^{n(p_1)}\colon A\to A^{p_1}$ spreads over
$T^{p_1-1}$;
\item $p_r=1$ and the Koebe space of $f^{m-n(p_r)-1}\colon fA^1\to T^0$
spreads over $T^{-1}$.
\end{itemize}
\end{lemma}
\pr
First we observe that the boundary points of each interval $T^i$
are nice. To see this, notice for example that $T^2$ is
the intersection of the branches $T^{2,i}$
of $g^i$ containing $c$
for $i=1,\dots,s$. Here $g^i|T^{2,i}$ is unimodal with $g^i(\partial
T^{2,i})\subset T^0$ because $g(\partial T^{2,i})\subset T^{2,i-1}$. 
By the choice of $s_0$ and $s$ we then
have that $g^i(\partial T^2)$ always remains outside $T^2$
(it is even outside $T^1$ for $s_0\le i\le s$). Hence $\partial T^2$
are nice points. Similarly, $\partial T^i$ are nice points for $i\le k$.
Secondly, from the definition of $s_0$ it follows that
$g(T^1)\cap T^2=\emptyset$ or more generally that
$\R^p g(T^{p+1})\cap T^{p+2}=\emptyset$ for $p=0,1,\dots,k-1$.
Hence consider the first return map $R$ to $T^{p-1}$.
It has a central interval $C^p\subset T^p$ and restricted to $C^p$
one has $R|C^p=\R^{p-1} g$; this holds
because $\R^{p-1}g(T^p)\cap T^{p-1}\ne \emptyset$
by  the choice of $s$ and because, as we just remarked,
$\R^ig(T^{i+1})\cap T^{p-1}=\emptyset$ for $i< p-1$.
Moreover, at each point $x$
its image $R(x)$ is an iterate of $\R^{p-1} g$.
In particular, it follows that $f^iA$ does not intersect
$\partial T^p$ because $\R^{k-1}g(\partial T^p)\subset \partial T^0$
and because $f^{m-1}(A)= T^0$.
All this implies also 
that $f^i(A)\cap T^{p+1}=\emptyset$ for $i< m-1$.
The third observation is that
if $j$ is the smallest integer such that
$f^j(c_1)\in T^{p-1}$ and $H\ni c_1$ is the largest interval
on which $f^j$ is monotone then $f^j(H)\supset T^{p-1}$.
This follows from the fact that $\partial T^{p-1}$ are nice points.
The fourth observation is the following:
consider the first return map $R$ to $T^{p-1}$. It has
a central branch $C^p$ contained inside $T^p$. Now other
branches are mapped by $R$ onto $T^{p-1}$. So if
$J\subset T^{p-1}\setminus C^p$
is so that it remains inside $T^{p-1}\setminus C^p$ for the
first $j$ iterates
of $R$ then $R^j\colon J\to T^{p-1}$ has an extension onto $T^{p-1}$.
The first assertion follows follows from observations 3 and 4.

So let us prove the second assertion of this lemma.
By definition, $\R^{k-1}g(T^k)\cap T^{k-1}\ne \emptyset$.
Hence there exists $n<m$ such that $f^n(A)\subset T^{k-1}$
and by the first observation in the proof the Koebe space
of this map contains $T^{k-1}$. Moreover, $f^i(T^{k+1})\cap T^{k+1}=\emptyset$
for $i=1,\dots,m-2$ and so it follows that $p_1=k$ or
$p_1=k+1$. By the second observation
if $n$ is the largest integer with $f^n(A)\subset A^{p_1}\subset
T^{p_1-1}$
then the Koebe space of this map still contains $T^{p_1-1}$.

The last assertion similarly follows from the
two observations made at the beginning of this proof.
We should observe that $g$ has a low return iterate
and that therefore $p_r=1$.
\qed

\begin{lemma}
There exists $K>0$
and given $\lambda>1$ there exists
$\sigma>0$ and $k_0\in \nz$ with the following property.
Assume that $g\in {\cal E}(T^0,T^{-1})$ as above and that
$|T^{-1}|\ge \lambda |T^0|$.
\begin{enumerate}
\item If all renormalizations $g,\R g, \R^2 g,\dots$
have low return iterates then there exists $k\le k_0$ such that
$\R^k g$ has a polynomial-like extension $G\colon \cup D_i\to D_*(T^0)$;
\item If $k$ is minimal such that
$\R^k g$ does {\bf not} have a low return iterate
then
either $\R^k g$ has a polynomial-like extension $G\colon \cup D_i\to D_*(T^0)$
or
$(|T^0|/|W^{k+1}|)\ge (1-\sigma)^{-1} \cdot (|T^{-1}|/|T^0|)$.
where $W^{k+1}$ is the escape interval associated to $\R^k g$.
\item If $|T^{-1}|/|T^0|>K$ and $k$ is minimal such that
$\R^k g$ does not have a low return iterate
then $\R^k g$ has a polynomial-like extension
$G\colon \cup D_i\to D_*(T^0)$.
\end{enumerate}
\end{lemma}
\pr
Suppose that $g,\dots,\R^{k-1}g$ have low returns
so that $\R^k g\colon \cup I_i\to T^0$ is well defined.
As before, the Schwarz Lemma implies that the
pullbacks of $D_*(T^0)$ under
the extensions of the monotone branches $I_i\to T^0$ of $\R^k g$
fit inside $D_*(T^0)$. So let us consider the pullback
associated to the inverse of the map $\R^kg|T^{k+1}=f^m$
on the central interval $T^{k+1}$.

First we notice that since $g\in \E(T^0,T^{-1})$
there exists because of the Koebe Principle a constant
$\rho>0$ which depends on $\lambda=|T^{-1}|/|T^0|>1$
such that the $\rho$-scaled neighbourhood
of each domain $T_i^1$ is still contained in $T^0$.
Here we use that the interval $H$ from the definition of $\E (T^0,T^{-1})$
is mapped inside $T^0$ by $f^{-1}$.
In particular, $|T^1|<(1+\rho)^{-1}|T^0|$.
In the same way we have
\beq
|T^{i+1}|<(1+\rho)^{-1}|T^i|\mbox{  for each }i=0,1,\dots k.
\label{boaw}
\eeq

Let $A\supset f(T^{k+1})$ be so that
$f^{m-1}$ maps $A$ diffeomorphically onto $T^0$.
Let $R^0=\R^kg(T^{k+1})=f^m(T^{k+1})$ and $I^0=T^0\setminus R^0$.
We want to compare the sizes of the
pullbacks in $I',R'\subset A$ by $f^{m-1}\colon A\to T^0$
of the `real' and `imaginary' pieces $R^0,I^0\subset T^0$.
That is, $R'=f(T^{k+1})$.
Let $p_i$ be as in the previous lemma and let
$R^i=f^{n(p_i)}(R')$ and $I^i=f^{n(p_i)}(I')$
be the partition of $A^i$ corresponding to $I'$ and $R'$.

Write $\mu_{i+1}=\frac{|T_{i+1}|}{|T_i|}\in (0,1)$.
If $R,I$ are two intervals in $T_i\setminus T_{i+1}$
with a unique common point then using the fact
that $f(z)=z^\ell+c_1$ it easy to see that
\beq
\frac{|R|}{|I|}\ge \left[\mu_{i+1}\right]^{\ell-1}
\cdot \frac{|f(R)|}{|f(I)|}
.
\label{impn}
\eeq
Indeed, if $T_i=(-b,b)$ and $T_{i+1}=(-a,a)$ with $0<a<b$
and $I,R$ are contained in $(a,b)$ then
$(|f(R)|/|R|)/(|f(I)|/|I|)$ is maximal when $I=\{a\}$ and $R=(a,b)$.
So (\ref{impn}) follows from the inequality
$$\frac{b^\ell-a^\ell}{b-a}\frac{1}{\ell a^{\ell-1}} \le
\frac{b^{\ell-1}}{a^{\ell-1}}.$$
In fact, it is easy to see that because of (\ref{boaw})
there exists $\tau>0$ so that either
\beq
\frac{|R|}{|I|}\ge
\left[\mu_{i+1}\right]^{\ell-1-\tau}  \cdot   \frac{|f(R)|}{|f(I)|}
\quad \mbox{ or } \quad \frac{|R|}{|I|}\ge 1
.\label{imp}
\eeq

\medskip

Now we will start pulling back the intervals $R^0,I^0$.
Let $K_1\supset f(A^1)$ be the interval which is mapped monotonically
onto $T^{-1}$ by the extension of $f^{m-n(p_r)-1}\colon f(A^1)\to T^0$.
Let $E^0$ be the component of $T^{-1}\setminus T^0$
which is adjacent to $I^0$ (the `extension' in the `imaginary' direction,
see the figure below).

\vskip0.3cm

\hbox to \hsize{\hss\unitlength=1.3mm
\beginpic(70,25)(0,0) \let\ts\textstyle
\put(70,14){$T^{-1}$}
\put(60,11){$T^0$}
\put(8,12){\line(1,0){32}}
\put(0,14){\line(1,0){50}}
\put(8,16){\line(1,0){21}}
\put(30,16){\line(1,0){10}}
\put(41,16){\line(1,0){9}}
\put(17,18.5){$R^0$}\put(16,3){$1$}
\put(34,18.5){$I^0$}\put(26,3){$\frac{1}{t}$}\put(34,3){$1-\frac{1}{t}$}
\put(44,18.5){$E^0$}
\put(24,-2){\line(1,0){26}}\put(54,-2){$t(1-\sigma)$}
\put(24,15.4){$\bullet$}
\put(24,17.8){$c$}
\put(3,-12){{\it The intervals $T^0=I^0\cup R^0$
and $T^{-1}=T^0\cup E^0\cup \tau(E^0)$.}}
\endpic\hss}
\vskip2cm

If $\R^k g$ does not have a low return iterate, then
we have that $R^0=f^m(T^{k+1})\supset W^{k+1}$ where $W^{k+1}$
is the escape interval associated to $\R^k g$.
If the inequality in assertion (2) in the statement of the
lemma does not hold then, writing $t=|T^0|/|W^{k+1}|$
and defining
$\sigma$ as in (2), the relative size of the intervals
$R^0$, $I^0$ and $E^0$ can be estimated as in the figure above.
Therefore,
\beq
C^{-1}([E^0,R^0],I^0)
\ge \frac{\left(1+1/t\right)
\left(t(1-\sigma)-1\right)}{\left(1-1/t\right)\left(1+t(1-\sigma)\right)}
\ge 1-\epsilon(\sigma)
\label{spab}
\eeq
where $\epsilon(\sigma)$ is some function so that $\epsilon(\sigma)\to 0$
as $\sigma\to 0$.
Here we have used that $t\ge |T^0|/|T^1|>1$ is bounded away from $1$
because of (\ref{boaw}).
Hence, using the map which sends an interval $H\supset R^1\cup I^1$
diffeomorphically onto $R^0\cup I^0\cup E^0$ we get
$|f(R^1)|/|f(I^1)|\ge C^{-1}(H,I^1)\ge 1-\epsilon(\sigma)$.
By (\ref{impn}) this gives
\beq
\frac{|R^1|}{|I^1|}\ge 
(1-\epsilon(\sigma))\mu_1^{\ell-1}
=(1-\epsilon(\sigma))\mu_{p_r}^{\ell-1}.
\label{kap1}
\eeq
On the other hand, if $\R^k g$ has a low return then,
because $\R^k g(T^{k+1})=f^m(T^{k+1})$ intersects
$T^k$, we have
by the Koebe Principle
some constant $\delta>0$ which depends on $|T^{-1}|/|T^0|$
such that $|f(R^1)|/|f(I^1)|\ge \delta$ and therefore
\beq
\frac{|R^1|}{|I^1|}\ge \delta \mu_1^{\ell-1}=\delta \mu_{p_r}^{\ell-1}
.
\label{kap2}
\eeq
Let us now compare $|R^{p_{i+1}}|/|I^{p_{i+1}}|$ with $|R^{p_i}|/|I^{p_i}|$.
Here $k+1\ge p_1>\dots>p_r:=1$ are defined as in the previous lemma.
There exists an interval $K_{p_i}\supset f(A_{p_i})$
which is mapped monotonically onto $T^{p_{i+1}-1}$
by $f^{s(p_i)}$.
Let $E^{p_{i+1}}$ be the component of $T^{p_{i+1}-1}\setminus A^{p_{i+1}}$
containing $c$ (again the extension).
Then because $I^{p_{i+1}}\cup R^{p_{i+1}} =A^{p_{i+1}}\subset
T^{p_{i+1}-1}\setminus T^{p_{i+1}}$ and $E^{p_{i+1}}\supset T^{p_{i+1}}$
and because $I^{p_{i+1}}$ is between $c$ and $R^{p_{i+1}}$, we have
$$\frac{|f(R_{p_i})|}{|f(I_{p_i})|}
\ge C^{-1}(K_{p_{i}},f(I^{p_{i}}))\ge$$
$$\ge
C^{-1}([E^{p_{i+1}},R^{p_{i+1}}],I^{p_{i+1}})\ge
\frac{|R^{p_{i+1}}|}{|I^{p_{i+1}}|}\frac{|T^{p_{i+1}}|}{|T^{p_{i+1}-1}|}
\ge \mu_{p_{i+1}} \frac{|R^{p_{i+1}}|}{|I^{p_{i+1}}|}.$$
Hence, using (\ref{impn}),
\beq
\frac{|R^{p_i}|}{|I^{p_{i}}|}
\ge \mu_{p_{i+1}}\mu_{p_i}^{\ell-1} \frac{|R^{p_{i+1}}|}{|I^{p_{i+1}}|}
\label{kapk}
\eeq
If $\R ^kg$ has no low return iterate and
the inequality in assertion (2) in the statement of the
lemma is not satisfied then combining (\ref{kap1}) and (\ref{kapk}) we get
$$\frac{|R^{k+1}|}{|I^{k+1}|}=
\frac{|R^{p_1}|}{|I^{p_1}|}\ge (1-\epsilon(\sigma))
(\mu_{p_r}\cdot \dots \cdot \mu_{{p_2}})^\ell \mu_{p_1}^{\ell-1}$$
and applying an estimate as the one above (\ref{kapk})
to the map $f^{n(p_1)}\colon A\to A^{p_1-1}$ we obtain,
\beq
\frac{|R'|}{|I'|}\ge
(1-\epsilon(\sigma))(\mu_{p_r}\cdot \dots \cdot \mu_{p_1})^\ell \ge
(1-\epsilon(\sigma))\left(\frac{|T^{k+1}|}{|T^0|}\right)^{\ell}
.
\label{chn}
\eeq
In fact we can improve this: using in all the previous
inequalities (\ref{imp}) instead of
(\ref{impn}), we get
\beq
\frac{|R'|}{|I'|}\ge
(1-\epsilon(\sigma))  \cdot
\left(\frac{|T^{k+1}|}{|T^0|}\right)^{\ell-\tau}
\label{delb}
\eeq
because if the second possibility in (\ref{imp}) holds
for $i=j$ with $j$ minimal, then as above (but without using (\ref{kap1})
we have
$|R^k|/|I^k|\ge \mu_{p_{j}} (\mu_{p_{j-1}}\cdot \dots \cdot \mu_1)^{\ell-\tau}$
which gives even a better bound than (\ref{delb}).
Since $|T^{k+1}|/|T^0|\le |T^1|/|T^0|$ is uniformly bounded away from $1$,
see (\ref{boaw}), it follows from (\ref{delb})
that there exists $\kappa>1$ such that provided
$\sigma>0$ is sufficiently small
\beq\frac{|R'|}{|I'|}\ge \kappa \cdot
\left(\frac{|T^{k+1}|}{|T^0|}\right)^{\ell}
\eeq
Since $|R'|=|T^{k+1}|^\ell$ and the pullback under $f^{m-1}$
of $D_*(T^0)$ is inside $D_*([R',I'])$ where $c_1$ is the unique
common point of $R'$ and $I'$,
the last inequality implies that the pullback of the disc $D_*(T^0)$ along
the central branch fits again inside $D_*(T^0)$, showing
that $\R^k g$ has a polynomial-like extension.
This proves assertion (2).

If $\R ^k g$ has a low return then combining (\ref{kap2}), (\ref{kapk})
and also the improved inequality (\ref{imp}), we get
$$\frac{|R'|}{|I'|}\ge
\delta \cdot (\mu_{p_r}\cdot \dots \cdot \mu_{p_1})^{\ell-\tau}
\ge \delta  \left(\frac{|T^{k+1}|}{|T^0|}\right)^{\ell-\tau}.$$
Because of (\ref{boaw}) we have when $k$ is large
that this last term is $\ge 2\left(|T^{k+1}|/|T^0|\right)^{\ell}$
and again the central pullback is mapped inside itself.
From this we get that either there exists $k$ such that
$\R^k g$ does not have a low return iterate or alternatively
$\R^k g$ has a polynomial-like extension. This proves assertion (1).

Let us now prove assertion (3) of the lemma.
Since the last return is high
the expression in (\ref{spab}) can be replaced
by $1\times K/(2+K)$ which tends to $1$
as $K\to \infty$.
Hence
then $|R'|/|I'|\ge K/(2+K) \left(|T^{k+1}|/|T^0|\right)^{\ell-\tau}$
becomes larger than $1$ when $K$ is large because
$|T^{k+1}|/T^0|$ is bounded away from one.
Again we get a polynomial-like mapping.
\qed

{\em Proof of Theorem C:} Let $g\ni \E (T^0,T^{-1})$
be the first return map as in the beginning of this section.
From the previous lemma it follows that $\R^k g$ does not
have a low return iterate for some $k$.
The new map $\W \R^k g$ (as defined above Lemma~\ref{lyno})
is then defined on a smaller domain.
If $\R^k g$ has no polynomial-like extension then
because of the second assertion in the previous lemma,
the corresponding Koebe space
increases by a definite factor $(1-\sigma)^{-1}>1$ (relative
to the size of the new domain).
Applying this idea several times, either one
obtains a polynomial extension at some stage or
the Koebe space becomes arbitrarily large
(compared to the size of the domains).
But from the last assertion
of the previous lemma one then also obtains a polynomial-like
extension.
\qed

Exactly as in Section 12 one has that Theorem C implies the
Main Theorem for each non-renormalizable map with a minimal critical 
point $c$.
If a map is only finitely often renormalizable then again the same
argument can be used: construct the Yoccoz puzzle
associated to the fixed points of the last renormalization
and apply Theorem C also to the last renormalization.
Thus the proof of  the Main Theorem is concluded.
\qed

\newpage

%

\begin{center}
{\huge An Extension and an Erratum}
\end{center}


\section{Theorem A holds for real analytic maps}

Let us first remark that Theorem A in the paper holds
for real analytic maps also. This means that the complex
bounds which Sullivan proved for infinitely renormalizable
Epstein maps of bounded type, even hold for 
arbitrary infinitely renormalizable maps which
are analytic on the dynamical interval. 
This answers a question of W. de Melo
and gives the possibility to extend
certain renormalization results of Sullivan and McMullen
to the class of real analytic maps.

(Let us also note that the generalized polynomial-like map in 
Theorems A-C have the property that the critical point does
not leave the domain of definition under iterates of this 
polynomial-like map.)

\begin{theo} Theorem A holds for a real analytic
unimodal map $f$ which is infinitely renormalizable:
there exists $N(f)$ such that
when $V_n$ is a periodic central interval of $f$
of period $s(n)\ge N(f)$, 
then there exists a polynomial-like
extension $F_n\colon \Omega_n'\to \Omega_n$
of $f^{s(n)}\colon V_n\to V_n$
such that the modulus of $\Omega_n\setminus \Omega_n'$
is universally bounded from below by some
positive number which only depends on $\ell$
and so that the diameter of $\Omega_n$ is of the
same order as that of $V_n$. The number $N(f)$
is uniformly bounded when $f$ runs over a compact space of maps.
\end{theo}
\pr First we prove that the real bounds from Sections 5 and 6 still
hold if $f$ is a of class $C^{1+\mbox{zygmund}}$
with a non-flatness condition at $c$ 
(see Section IV.2.a in [MS]). 
Because of Theorem IV.2.1 from [MS],
this means that if $J\subset T$ are intervals such that $f^s|T$ is
diffeomorphic then
$$C(f^sT,f^sJ)\ge \prod_{i=0}^{s-1} \left(1-o(|f^i(T)|)\right)
C(T,J)$$
where $o(t)$ is some function such that $o(t)\to 0$ as $t\to 0$.
Now in Section 5 of [LS] let $l$ be maximal so that
as before $f^s|l$ is monotone and - this is new - 
$L=f^s(l)$ contains at most $5$ iterates of $V$.
Therefore, in the renormalizable case,
the orbit of $T,\dots, f^k(T)$ in the proof of
Lemma 5.1 has intersection multiplicity bounded by $15$.
Moreover, because the map has no wandering intervals,
one has $\max _{i=0,\dots,k}|f^i(T)|$ tends to zero
if the period tends to infinity. 
(Note also that if $I_n,I_{n+1}$ are consecutive 
central interval of $f$ then 
\beq
|I_{n+1}|\le \lambda |I_n|
\label{eqqm1}
\eeq
where $\lambda<1$ uniformly 
when $f$ runs over a compact space of maps. This follows 
from the extension given by Proposition 7.1
in [LS] and the Koebe Principle.)
In particular, the inequality proved in Lemma 5.1 still
holds with a spoiling factor $O_s$ such that
$O_s\to 1$ as the period $s$ tends to infinity.
Now take in Lemmas 6.2-6.5 also the
same definition for $l$. Then these lemmas still
holds with a spoiling factor $O_s$. In Lemma 6.3
simply note that if $Q_k$ contains more than $5$ iterates
of $V$ then one simply takes $\tilde Q_1\supset f(V)$ so that
$\tilde Q_k=f^{k-1}(\tilde Q_1)$ 
contains $f^k(V)$ and precisely $5$ iterates of $V$. 
Then, because the intersection multiplicity of the
orbit $\tilde Q_1,\dots, f^{k-1}(\tilde Q_1)$ 
is bounded by $15$ and in the same way as before 
we get $C^{-1}(Q_1,f(V))\ge C^{-1}(\tilde Q_1,f(V))
\ge O_s C^{-1}(\tilde Q_k,f^k(V))\ge O_s \cdot 0.6$.
If $Q_k$ contains less than $5$ intervals
then we can obtain $C^{-1}(Q_1,f(V))\ge O_s \cdot 0.6$
in the same cases as before. Only in cases II.b and
II.c we used the interval $Z_1$. But now 
notice that the arguments used there
also apply if replace $Z_1$
by the maximal interval $\tilde Z_1\subset Z_1$ in $H_1$
so that each component of $f^{k-1}(\tilde Z_1)\setminus
Q_k$ contains at most one iterate of $V$. 
Since the intersection 
multiplicity of $\tilde Z_1,\dots,f^{k-1}(\tilde Z_1)$ 
is bounded by $18$ we get that Lemma 6.3 still holds
with a spoiling factor. In Lemmas 6.4 and 6.5 exactly the
same remarks apply. Now in Lemma 8.1 we redefine
$T_1$ as the maximal interval such that $f^{s-1}|T_1$ 
is monotone and such that each component
of $f^{s-1}(T_1)\setminus V$ contains at most $5$
iterates of $V$. So we can still apply Lemma 6.4 in the
proof of Lemma 8.1 to this $T_1$
and so this lemma still holds.

All this implies that the same real 
bounds can be still used in Sections 8 and 9.
Now of course, the Schwarz Lemma 
(that the pullback of some Poincar\'e domain
with angle $\theta$ maps inside a similar region
with the same angle $\theta$) 
which we used
in these sections does not hold anymore,
because the map $f$ is only analytic on a small
neighbourhood of the dynamical interval. However,
in Lemmas VI.5.2 and VI.5.3 of [MS] it is proved that
we can still essentially obtain 
the same inclusion but with a slightly worse
angle. According to Lemmas VI.5.2 and VI.5.3
the loss in the estimate 
tends to zero if the size of the interval tends to zero. 
Therefore we still get the same estimate
in the proof of Theorem A.

The statement that $N(f)$ is uniformly bounded when
$f$ runs over a compact space follows
from (\ref{eqqm1}).
\qed

\section{An erratum}
Firstly, we should point out that the domains of
the polynomial-like
mapping in Theorem C is are disjoint because
the $f$-images of these (near $c_1$) are based on
disjoint intervals in the real line.
Moreover, as Ben Hinkle pointed out, there
is a mistake in the non-renormalizable case when
we prove local connectivity outside the critical point
(on page 42 lines 9-11 it is mistakenly argued that
the sum of the moduli of some annuli 
is infinite in the non-renormalizable case). 
We like to thank Ben Hinkle for this comment and
let us show how to fix the proof.
We shall show that one can argue
as in the proof of the local connectivity of
the Julia set of infinitely renormalizable maps
in Section 8. To do this we have to be a little
careful since the orbit of $c$ enters perhaps
several times in $\Omega\cap \rz$ at times
which do not correspond to iterates of the polynomial-like
mapping $R$. 

So assume that $f$ is non-renormalizable
and that $\omega(c)$ is minimal.
We show that the bounds from Theorem B and C imply
local connectivity.

\begin{prop} 
\label{proer1}
Let $G(j)\colon \cup_i\Omega_i(j)\to \Omega(j)$
be a sequence of polynomial-like mappings associated
to a real polynomial $f(z)=z^\ell+c_1$ such that
the critical point $c=0\in \Omega_0(j)$ does not
escape the domain of $G(j)$ under iterations of $G(j)$. 
(As before, we assume 
$G(j)\colon \Omega_0(j)\to \Omega(j)$ is $\ell$-to-one,
and each other $G\colon \Omega_i(j)\to \Omega(j)$ is an isomorphism.)
Assume moreover, that there exist interval neighbourhoods
$X(j)$ of $\Omega_0(j)\cap \rz$ so that
when $x,f^i(x)\in X(j)$ then $f^i(x)$
is an iterate of $x$ under $G(j)$
(we call this {\rm the first return condition}) 
and so that the modulus of
the annuli $\cz_{X(j)}\setminus \Omega_0(j)$
is uniformly bounded away from zero. 
Then the Julia set of $f$ is locally connected.
\end{prop}
\pr
If $z$ is in the Julia set but $\omega(z)$ does not contain $c$ then 
the Julia set is locally connected at $z$ because of
the contraction principle. So choose
a point $z$ from the Julia set of $f$ so that $\omega(z)\ni c$.
Let $P_n$ be an open piece of the Yoccoz puzzle (corresponding to $\Omega$)
based on two preimages $v,-v$  
of the orientation reversing fixed point of $f$
so that $\Omega\cap \rz$ is
either equal to $[-v,v]$ or to a small neighbourhood of this interval.
There exists a large integer $N$ such that the full preimage $G^{-N}(P_n)$
is inside the domain of definition $\cup_i\Omega_i$ of $G$,
see Section 3.
Note that $G^{-N}(P_n)$ consists of finitely many (open) Yoccoz pieces.
Let us consider the pieces of $G^{-N}(P_n)$ 
inside the central domain
$\Omega_0$, i.e., 
$$P'_n= G^{-N}(P_n)\cap \Omega_0.$$
Since $\omega(z)\ni c$,
there exists a {\it minimal} $k$ such that $f^k(z)\in P'_n$.
In particular, the point $f^k(z)$ 
belongs to one of the Yoccoz pieces inside 
$\Omega_0$.
Let $F$ be the branch of $f^{-k}$ which maps
a neighborhood of $f^k(z)$ to a neighbourhood of $z$. 
Let $X$ be as in the statement of this lemma.

\noindent
{\it Claim 1.} The map $F$ extends to a holomorphic
map in the domain $\cz_{X}$.
\noindent
{\it Proof of the claim.} Assume the contrary.
We then get that for some minimal $r< k$ 
that $f^{-r}(\cz_X)$ (along the same
orbit) meets the critical value
$c_1$.  
This means that the branch $f^{-r}$ follows
the points $c_{r+1}=f^r(c_1)\in \cz_X$, $c_r=f^{r-1}(c_1)$,...,
$c_2=f(c_1)$, $c_1$.
Among these iterations of $c_1$, 
let us mark all those $c_{j_1}, c_{j_2},\dots, c_{j_m}$, 
where $j_1<j_2<\dots<j_m$, which hit the domain $\cz_X$
(i.e., are in $X$).
Because of the first return assumption 
there exists integers $k(1)<k(2)<\dots$ such that
$c_{j_1}=G^{k(1)}(c)$,
$c_{j_2}=G^{k(2)-k(1)}(c_{j_1})=G^{k(2)}(c)$,...,
$c_{j_m}=G^{k(m)}(c)$. It follows, that
$f^{-r}=f^{-(s-1)}\circ G^{-(k(m)-1)}$,
where
$f^{-(s-1)}$ is the branch from $V$ to $\hat U$ corresponding to
the restriction of $G$ on $\Omega_0$ (so $G|\Omega_0=f^{s-1}\circ f$). 
Hence, $f^{-(r+1)}(\cz_X\cap \Omega)\subset
G^{-r(m)}(\Omega)\subset \Omega_0$ and $f^{k-r-1}(z)\in f^{-(r+1)}(P'_n)=
(G|_{\Omega_0})^{-1}\circ G^{-k(m)+1}(P'_n)\subset P'_n$.
This contradicts the minimality of $k$ and proves
the claim. 

Now apply the claim to a sequence
of maps $G(j)\colon \cup \Omega_i (j)\to \Omega(j)$.
This gives a sequence of annuli $\cz_{X(j)}\setminus \Omega_0(j)$
of modulus $\ge \delta$ such that some iterate 
of $z$ maps to a puzzle piece inside $\Omega_0(j)$.
Since the diameter of $\Omega_0(j)$ shrinks to zero
one completes the proof exactly like 
in the infinitely renormalizable case.
\qed

\begin{corr} If $f$ has infinitely many high first
returns, then the Julia set is locally connected.
\end{corr}
\pr
For each high return, we have a polynomial-like
mapping $G:\cup_i\Omega_i\to \Omega$ constructed in Sections 11-12,
such that $G|_V$ is the first return map. Remember that
$\Omega$ here is a definite
complex neighborhood of the interval $V$ so that $\Omega_0$ is inside a
definite neighborhood of the interval $U=[-u,u]\subset V$ and so
that $G|_{\Omega_0}=f^s$. Denote by $l'$ a maximal interval outside $U$
with a common boundary point, such that $f^s|_{l'}$ is monotone,
and let $f^{ks}(u)$ be the first moment when it leaves $l'$. Then any
$f$-iterate of $c$ in the interval $X=[-f^{ks}(u), f^{ks}(u)]$ is, in fact,
an iterate of $c$ under $G$. Hence, we can apply the Proposition 2.1
(note that the gaps $X\setminus \Omega\cap \rz$ are not small because of
Proposition 7.1 [LS]. 
\qed

Let $W_n$ be the sequence of intervals as in Section 2
and let $R_n$ be the corresponding first return maps.
We say that the return is {\it low} if the image of $R_n$ of 
the central component $W_{n+1}$ does not contain $c$
and it is called {\it central} if $R_n(c)\subset
W_{n+1}$.

\begin{lemma} There exists
a universal number $\lambda>1$ (only depending
on $\ell$) the following property. 
Assume that we are in one of the following situations:
\noindent
1) either $R_{n-1}$ or $R_n$ has a non-central low return;
\noindent
2) the return of $R_{n-1}$ is non-central high.
\noindent
Then $|W_n|\ge \lambda |W_{n+1}|$.
\end{lemma}
\pr
If $R_n$ has a non-central low return 
then $|W_n|\ge \lambda |W_{n+1}|$ according 
to the corollary on page 345 in [MS].
In the same way, if the return to $W_{n-1}$ is non-central low then 
again $|W_{n-1}|\ge \lambda |W_n|$.
Now let the first return map to $W_n$ restricted
to $W_{n+1}$ be equal to $f^s$.
There exists an interval neighbourhood $T$ of $f(W_{n+1})$
such that $f^{s-1}$ maps $T$ diffeomorphically onto
$W_{n-1}$ and so that $f^{-1}(T)\cap \rz\subset W_n$.
Hence, using the Koebe Principle it follows that
$f^{-1}(T)\cap \rz \subset W_n$
is a definite neighbourhood of $W_{n+1}$.
So in this case we are done also.

If the return map $R_{n-1}$ to $W_{n-1}$ is high
then according to part 1 of Lemma 1.2 on page 342 in [MS]
(or Proposition 7.1 of the present paper), 
the map on the
central domain $W_n$ is a composition of $f$ and a map 
with bounded distortion. From this and the fact that 
$R_{n-1}$ has a non-central high return it follows that
$W_{n+1}$ has to be a definite factor smaller than $W_n$
also.
\qed

Of course, there are infinitely many integers
$n$ for which the map $R_{n-2}$ has a non-central return.
Choose such an $n$ and write $T^0=W_n$, $T^{-1}=W_{n-1}$
and study the situation as in Section 14.
So take a first return map
$g$ as in Section 14 and define $T^{2,i-1}$ to be the component
of $g^i$ containing $c$. 
Let us begin by remarking that on page 59 line -8
one better defines $s$ minimal so that 
$g^s(T^{2,s-1})\cap T^1\ne\emptyset$ 
and on line -3 we should have  defined $s(x)$ 
as the smallest integer for which
$g^{s(x)}(x)$ and $g^{s(x)}(c)$ are in different components
of $\cup_i T_i^1$ (so we have to also allow $T_0^1$.)
Moreover on page 60 line 17 one should read $g^i\circ g$
instead of $g^{i+1}\circ g$.
First we need the following proposition.

From the lemma above, there exists a universal
$\lambda>1$ (only depending on $\ell$) such that 
\beq
|T^{-1}|\ge \lambda |T^0|. 
\label{eqq1}
\eeq
Now consider the return map $g$ to $T^0$ with
central domain $T^1=W_{n+1}$. 

\begin{lemma} 
\label{lemer2}
Suppose that $g,\R g, \dots, \R^k g$ exist.
Then there exists a universal $\kappa>0$ (only 
depending on $\ell$) and a $\kappa$-scaled neighbourhood
$W^{k+1}$ of $T^{k+1}$ such that 
each iterate of $c$ inside $W^{k+1}$ is 
an iterate of $c$ under the map $\R^k g$. 
\end{lemma}
\pr
First notice that there exists for
each component $I$ of $g$ an integer $k$ 
such that $g|I=f^k$. Moreover, there exists $U\supset f(I)$
such that $f^{k-1}$ maps $U$ diffeomorphically onto $T^{-1}$
and so that $f^{-1}(U)\cap \rz\subset T^0$. 
In particular, some neighbourhood $U$ of $f(T^1)$
is mapped diffeomorphically onto $T^{-1}$
and because of (\ref{eqq1}) and by Koebe this implies that
$U$ contains a definite neighbourhood of $f(T^1)$
and since $f^{-1}(U)\cap \rz\subset T^0$
one gets that some definite neighbourhood of $T^1$ is contained
in $T^0$. Hence there exists a universal $\lambda'>1$ (only
depending on $\ell$)
such that 
\beq
|T^0|\ge \lambda' |T^1|.
\label{eqq2}
\eeq
Similarly each component $I$ of the domain of $g$ in $T^0\setminus T^1$
has adjacent to it an interval $I'\subset T^0$ 
with one point in common with $I$ `further away from $c$'
(i.e., so that $I$ lies between $I'$ and $c$) 
so that $g|(I\cup I')$ is monotone and
\beq
|I'|\ge \rho |I|.
\label{eqq3}
\eeq
 
For $k=0$ the result is obvious.
Let us first show the lemma for $k=1$.
Let us first consider the case
that $g$ has a non-central low return. 
We claim that if $x, g^i(x)\in T^1$ then $g^i(x)$ is an iterate
of $x$ under $\R g$. To see this we first remark that
\beq
\R g |(T^0\setminus T^1)=g. 
\label{eqq0}
\eeq
Next take $x\in T^1$. 
If $g(x)$ is contained in a component of $\cup_{i\ne 0} T^1_i$
which is entirely contained in $g(T^1)$ 
(which is in $T^0\setminus T^1$ since the return
was assumed to be low) then we have
$\R g(x)=g^2(x)$. So if $g^2(x)\in T^1$ then the required statement
holds for $x$ and if $g^2(x)\in T^0\setminus T^1$ 
then by (\ref{eqq0}) again the required statement
holds for $x$. If $g(x)$ is not contained in such a component
then $x$ is contained in a symmetric interval $T^{2,1}\subset T^1$ 
such that $g(T^{2,1})$ is inside one of the components
of $\cup_{i\ne 0} T^1_i$ (in fact, $T^{2,1}$ is the component
of $g^{-1}(\cup_{i\ne 0}T_i)$ containing $c$). 
In particular $g(x)\notin T^1$.
So if $x$ is not as before
consider $g^2(x)$. If $g^2(x)$ is contained in a
component of $\cup_{i\ne 0}T^1_i$ which is entirely 
contained in $g(T^{2,1})$ then $\R g(x)=g^3(x)$
and by (\ref{eqq0}) again the required statement holds for $x$. 
If $g^2(x)$ is not contained in such a component
then again $x$ is contained in a symmetric
interval $T^{2,2}\subset T^{2,1}$ such that 
$g^2(T^{2,2})$ is inside one of the components of
$\cup_{i\ne 0}T^1_i$ and in particular $g^2(x)\notin T^1$.
In this way one proves the claim inductively.
Now we set $W^2=T^1$ 
and it suffices to show that $T^1$ is a definite amount
larger then $T^2$. To see this write $g|T^1=f^m$.
The map $f^{m-1}$ maps some neighbourhood of $f(T^1)$
diffeomorphically onto $T^0$ (in fact onto $T^{-1}$ we do not
need this anymore), and because the Koebe Principle
and because of (\ref{eqq3}) it follows that a definite
piece of $T^1$ is mapped by $g$ outside any given
component in $T^0\setminus T^1$ of $g$. Hence 
$|T^1|>(1+\kappa)|T^2|$ for some universal number $\kappa>0$.

If $g$ has a central low return then define $s_0$
as before and consider a (shrinking) nested
sequence of intervals $T^{2,s_0-2}\subset \dots T^{2,0}:=T^1$ 
such that $g(\partial T^{2,i})\subset \partial T^{2,i-1}$
for $i=1,\dots,s_0-2$
(so these intervals are associated to the saddle-cascade of the
central branch; they are symmetric around $c$
and their endpoints are preimages of $\partial T^1$ under
the central branch of $g$).

Since we now assume that one has a central return
we have $s_0>2$. Because of the corollary
on page 345 of [MS], one has a universal $\lambda>1$ (only
depending on $\ell$) such that 
$|T^{2,s_0-2}|\ge \lambda |T^{2,s_0-1}|$
and because of (\ref{eqq2}) and the Koebe Principle 
we get in the same way also 
\beq
|T^{2,0}|\ge \lambda |T^{2,1}|.
\label{eqq4}
\eeq 
(If one has a long saddle-cascade 
a similar uniform comparison between $T^{2,i}$ with $T^{2,i-1}$
is certainly not true.)
Now $s_0$ is by definition the minimal integer such that
$g^{s_0-1}(c)\notin T^1$. As in (\ref{eqq4}) one has
that $|g^{s_0-1}(c)-g^{s_0-2}(c)|$ is comparible
to size of $T^{2,0}\setminus T^{2,1}$, i.e., to
the size to $T^0\setminus T^1$. Write $g^{s_0-2}|T^{2,s_0-2}=f^m$.
Then it follows from this
that some neighbourhood of $f(T^{2,s_0-2})$
is mapped by $f^{m-1}$ onto a definite neighbourhood 
of $f^m(T^{2,s_0-2})$. So 
\beq
f^{m-1}|f(T^{2,s_0-2})\mbox{ has uniformly bounded distortion.}
\label{eqq5}
\eeq
Now define $W^2\supset T^2$ to be equal to $T^{2,s_0-2}$.
One has by construction that $g(T^1)\cap T^{2,s_0-2}=\emptyset$
and by definition of $\R g$ if $x,g^i(x)\in T^{2,s_0-2}$
then $g^i(x)$ is an iterate of $x$ under $\R g$.
So it suffices to show that $V^2$ is a definite 
amount larger than $T^2$.
But this follows from (\ref{eqq3}) and (\ref{eqq5}).

Now if $g,\dots,\R^k g$ are defined then we get 
$|T^i|\ge \lambda |T^{i+1}|$ for $i=0,\dots,k$.
If $\R^{k-1} g$ has a non-central low return
then and this we can set 
$V^{k+1}=T^k$ and we are done.
If $\R^{k-1} g$ has a central low return
then we argue as above with 
intervals $T^{k+1,i}$.
\qed

If $g,\R g,\R^2 g, \dots$ all exist
then there exists $\R ^k g$ with central domain $T^{k+1}$
which has a polynomial-like with central piece 
$\Omega_0(k)\supset T^{k+1}$
extension and Lemma~\ref{lemer2}
gives an interval $W^{k+1}\supset T^{k+1}$
such that the modulus of $\cz_{W^{k+1}}\setminus
\Omega_0(k)$ is uniformly bounded from below
and a first return condition is satisfied.
Because of Proposition~\ref{proer1} above
the Julia set is locally connected if this case
happens infinitely often.
If $g,\dots, \R^k g$ exist but $\R^k g$ has a high return then
according to Lemma 14.2 either 
$\R^k g$ has a polynomial-like extension and we can again use
Lemma~\ref{lemer2} above or one considers
$g_1=\W \R^k g$ which has a better extension domain.
Remark that $\W \R^k g$ is by construction
a first return map to its domain.
So we can apply Lemma 14.2 again when
$g_1,\dots,\R^{k_1} g$ are well defined
(with again the first return condition and
even better extension scale).
By assertions 2 and 3 of
Lemma 14.2 we must reach a situation which has
a polynomial-like extension. 
In this way we get
eventually a map which a polynomial-like extension
and satisfies the required conditions of Proposition~\ref{proer1} above.
Thus we get

\begin{corr}
The Julia set of $f$ is locally connected.
\end{corr}



\newpage
%

 %


\begin{thebibliography}{NvS2}
\bibitem[BH]{BH}
        B. Branner and J.H.Hubbard, {\em
        The iteration of cubic polynomials I,} Acta Math. {\bf 160}, (1988),
        143-206, {\em The iteration of cubic polynomials II,}
        Acta Math. {\bf 169}, (1992),  229-325.
\bibitem[BL]{BL}
        A.M. Blokh and M.Yu. Lyubich, {\em Measurable
        dynamics of S-unimodal maps of the interval},
        Ann. Sc. E.N.S. $4e$ s\'erie, {\bf 24}, 545--573, (1991).
\bibitem[BKNS]{BKNS}
        H. Bruin, G. Keller, T. Nowicki and
        S. van Strien, {\em Absorbing Cantor sets in dynamical systems:
        Fibonacci maps}, Stonybrook IMS preprint 1994/2. To appear
        in the Annals of Math.
\bibitem[DH]{DH}
         A. Douady and J.H. Hubbard, {\em
         On the dynamics of polynomial-like mappings,}
         Ann.Sc.E.N.S.4e s\'erie, {\bf 18}, (1985), 287--343
\bibitem[DH1]{DH1}
	A. Douady and  J.H. Hubbard, {\em Etude dynamique des polynomes
        complexes I, II,} Publication Mathematiques D'Orsay, 
        no. 84-02 (1984), no. 85-04 (1985).
\bibitem[GS]{GS}
         J. Graczyk and G. \'Swi\accent'30atek, 
         {\em Polynomial-like property for real quadratic polynomials},
         February 3, 1995
\bibitem[HJ]{HJ}
        J. Hu and Y. Jiang, {\em The Julia set of the Feigenbaum quadratic
        polynomial is locally connected,} Preprint 1993.
\bibitem[Ji1]{Ji1}
         Y. Jiang, {\em Infinitely renormalizable quadratic
        Julia sets}, Preprint ETH, Z\"urich (1993).
\bibitem[Ji2]{Ji2}
         Y. Jiang, {\em Renormalization in  quadratic-like mappings},
        Preprint November 1994.
\bibitem[KN]{fibo}
        G. Keller and T. Nowicki,
        {\em Fibonacci maps re(a$\ell$)visited,}
        Preprint (1992).
\bibitem[Ly1]{L1}
        M.Yu. Lyubich, {\em Ergodic theory for smooth one-dimensional
        dynamical systems}.
        Stonybrook IMS Preprint 1991/11.
\bibitem[Ly2]{L2}
        M.Yu. Lyubich, {\em On the Lebesgue measure of the
        Julia set of a quadratic polynomial}.
        Stonybrook IMS Preprint 1991/10.
\bibitem[Ly3]{L3}
        M.Yu. Lyubich, {\em Combinatorics, geometry and attractors
        of quasi-quadratic maps}. Annals of Math. {\bf 140}, (1994),
        347-404
\bibitem[Ly4]{L4}
        M.Yu. Lyubich, {\em Milnor's attractors, persistent
        recurrence and renormalization}. In:
        ``Topological Methods in Modern Mathematics, A Symposium in Honor
        of John Milnor's 60th Birthday'', Publish or Perish (1992).
\bibitem[Ly5]{L5}
        M.Yu. Lyubich, {\em Geometry of quadratic polynomials: moduli, rigidity
        and local connectivity,} Stonybrook IMS Preprint 1993/9.
\bibitem[LM]{LM}
        M. Lyubich and J. Milnor, {\em The unimodal Fibonacci map,}
        Journal of the A.M.S. {\bf 6}, 425-457 (1993).
\bibitem[LY]{LY}
        M. Lyubich and M. Yampolsky, {\em Dynamics of quadratic polynomials:
        complex bounds for real maps}, draft Februari 22, 1995.
\bibitem[Mar]{Mar}
        M. Martens, {\em Interval dynamics,}  Thesis, Delft, (1990).
\bibitem[MS]{MS}
        W. de Melo and S. van Strien, {\em One-dimensional dynamics,}
        Ergebnisse Series {\bf 25}, Springer Verlag, (1993).
\bibitem[MMS]{MMS}
        M. Martens, W. de Melo and S. van Strien,
        {\em Julia-Fatou-Sullivan theory for real one-dimensional
        dynamics},  Acta Math. {\bf 168}, 273-318 (1992).
\bibitem[McM]{McM}
        C. McMullen, {\em Complex dynamics and renormalization,}
        Princeton University Press, to appear.
\bibitem[Mil]{Mil}
        J. Milnor, {\em Local connectivity of Julia sets; expository
        lectures}. Stonybrook IMS Preprint 1992/11.
\bibitem[Mil1]{Mil1}
        J. Milnor, Questions in: {\em Problems in Holomorphic
        Dynamics}, Stony Brook IMS Preprint 1992/7, (editors:
        B. Bielefeld and M. Lyubich).
\bibitem[Pe]{Pe}
        C. Petersen, {\em Local connectivity of some Julia sets
        containing a circle with an irrational rotation}, 
        Preprint I.H.E.S. /M/94/26 (1994). 
\bibitem[Sh]{Sh}  M. Shishikura, {\em Unpublished}.
\bibitem[SN]{SN}
        S. van Strien and T. Nowicki, {\em Polynomial maps with
        a Julia set of positive Lebesgue measure: Fibonacci maps}
        Stonybrook IMS Preprint 1994/3.
\bibitem[Sul]{Sul}
        D. Sullivan, {\em Bounds, quadratic differentials,
        and renormalization conjectures, 1990,}
        In AMS Centennial Publications. {\bf 2}:
        Mathematics into Twenty-first Century.
\bibitem[Y]{Y}
        J.C. Yoccoz, {\em MLC}, Manuscript.
\end{thebibliography}
\end{document}